\documentclass[11pt,a4paper,reqno]{amsart}

\usepackage{amsmath} 
\usepackage{amssymb}  
\usepackage{amsthm}
\usepackage{euscript}
\usepackage[T1]{fontenc}
\usepackage{mathtools}
\usepackage{relsize}
\usepackage{amsfonts}
\usepackage{xcolor}
\usepackage{graphicx}
\usepackage{mathrsfs}
\usepackage{tikz}
\usepackage{calc}
\usepackage[numbers]{natbib}
\usepackage[letterpaper,top=2cm,bottom=2cm,left=3cm,right=3cm,marginparwidth=1.75cm]{geometry}
\usepackage[colorlinks=true, allcolors=blue]{hyperref}

\pgfdeclarelayer{edgelayer}
\pgfdeclarelayer{nodelayer}
\pgfsetlayers{edgelayer,nodelayer,main}

\tikzstyle{none}=[inner sep=0pt]

\newtheorem{thm}{Theorem}[section]
\newtheorem{cor}[thm]{Corollary}
\newtheorem{lem}[thm]{Lemma}
\newtheorem{ass}[thm]{Assumption}
\newtheorem{prop}[thm]{Proposition}
\theoremstyle{definition}
\newtheorem{defn}[thm]{Definition}
\newtheorem{rem}[thm]{Remark}

\numberwithin{equation}{section}

\DeclareMathOperator{\gr}{gr}

\DeclareFontFamily{U}{mathx}{}
\DeclareFontShape{U}{mathx}{m}{n}{<-> mathx10}{}
\DeclareSymbolFont{mathx}{U}{mathx}{m}{n}
\DeclareMathAccent{\widehat}{0}{mathx}{"70}
\DeclareMathAccent{\widecheck}{0}{mathx}{"71}

\newlength{\widthofexpr}
\renewcommand{\Im}{\operatorname{Im}}

\newcommand{\dom}{\operatorname{dom}}

\newcommand{\ran}{\operatorname{ran}}

\newcommand{\mform}[3]{\left[\!\left[ #1, #2 \right]\!\right]_{#3}^-}

\newcommand{\mperp}{\left[\!\left[ \perp \right]\!\right]^-}

\newcommand{\mrperp}{\langle \langle \cdot , \cdot \rangle \rangle_- }

\newcommand{\R}{\ensuremath{\mathbb R}}    
\newcommand{\C}{\ensuremath{\mathbb C}}    
\newcommand{\N}{\ensuremath{\mathbb N}}    

\newcommand{\calB}{\mathcal B}

\newcommand{\calD}{\mathcal D}
\newcommand{\calE}{\mathcal E}
\newcommand{\calF}{\mathcal F}
\newcommand{\calG}{\mathcal G}

\newcommand{\calJ}{\mathcal J}

\newcommand{\calL}{\mathcal L}

\newcommand{\calP}{\mathcal P}

\newcommand{\calS}{\mathcal S}

\newcommand{\calW}{\mathcal W}
\newcommand{\calX}{\mathcal X}

\newcommand{\calZ}{\mathcal Z}

\newcommand{\euA}{\EuScript A}  
\newcommand{\euB}{\EuScript B}

\newcommand{\euJ}{\EuScript J}

\newcommand{\vphi}{\varphi}

\newcommand{\rank}{\operatorname{rk}}
\newcommand{\mul}{\operatorname{mul}}

\newcommand{\dif}[1]{\frac{\partial^{#1}}{\partial \xi^{#1}}}
\newcommand{\tdif}[1]{\tfrac{\partial^{#1}}{\partial \xi^{#1}}}

\DeclareMathOperator{\re}{Re}

\makeatletter
\let\@evenhead\relax
\let\@oddhead\relax
\makeatother

\title[]{Extension theory via boundary triplets for infinite-dimensional implicit port-Hamiltonian systems}

\author{H. Gernandt$^1$}
\thanks{$^1$IMACM Work Group Functional Analysis, School of Mathematics and Natural Sciences, University of Wuppertal, Germany, Mail: \textsc{gernandt@uni-wuppertal.de}}
\author{F. Philipp$^2$}
\thanks{$^2$Optimization-based Control Group, Institute of Mathematics, Technische Universität Ilmenau, Germany, Mail: \textsc{friedrich.philipp@tu-ilmenau.de}}
\author{T. Preuster$^3$}
\thanks{$^3$Junior Professorship Numerical Mathematics, Faculty of Mathematics, Chemnitz University of Technology, Germany, Mail: \textsc{\{till.preuster,manuel.schaller\}@math.tu-chemnitz.de}}
\author{M. Schaller$^3$}

\thanks{F.P. was funded by the Carl Zeiss Foundation within the project {\it DeepTurb--Deep Learning in and from Turbulence}. He was further supported by the free state of Thuringia and the German Federal Ministry of Education and Research (BMBF) within the project {\it THInKI--Th\"uringer Hochschulinitiative für KI im Studium} and he gratefully acknowledges funding by the German Research Foundation (DFG) – Project-ID 519323897. \\ 
\indent H.G., T.P. and M.S. acknowledge funding by the Deutsche Forschungsgemeinschaft (DFG, German Research Foundation) – Project-ID 531152215 – CRC 1701.}

\begin{document}

\begin{abstract}
    The solution of constrained linear partial-differential equations can be described via parametric representations of linear relations. To study these representations, we provide a novel definition of boundary triplets for linear relations in range representations where the associated boundary map is defined on the domain of the parameterizing operators rather than the relation itself. This allows us to characterize all boundary conditions such that the underlying dynamics is represented by a self-adjoint, skew-adjoint or maximally dissipative relation. The theoretical results are applied to a class of implicit port-Hamiltonian systems on one-dimensional spatial domains. More precisely, we explicitly construct a boundary triplet which solely depends on the coefficient matrices of the involved matrix differential operators and we derive the associated Lagrangian subspace. We exemplify our approach by means of the Dzektser equation, the biharmonic wave equation, and an elastic rod with non-local elasticity condition.
\end{abstract}

\maketitle

\section{Introduction}
The time evolution of various physical phenomena may be formulated as
\begin{align}\label{eq:eq1}
     \tfrac{\mathrm{d}}{\mathrm{d}t} \calP x(t) = \calS x(t), \qquad t \geq 0,
\end{align}
 where $\calP$ and $\calS$ are (possibly unbounded) linear operators, e.g.\ differential operators. Equivalently, we may study the inclusion $(z(t), \tfrac{\mathrm{d}}{\mathrm{d}t} z(t)) \in \calL$ with the linear relation $\calL = \ran \left[\begin{smallmatrix} \calP \\ \calS \end{smallmatrix}\right]$.
 Such multi-valued linear operators $\calL$, or simply subspaces, naturally occur in the port-Hamiltonian modeling formalism (therein called \emph{Lagrangian subspaces}) if are self-adjoint in the sense of linear relations. 
 
 Both, for port-Hamiltonian modeling and well-posedness analysis of the degenerate Cauchy problem \eqref{eq:eq1}, self-adjointness or (maximal) dissipativity of $\calL$ are central aspects, cf. \cite{knuckles1994remarks, arendt2023semigroups}. In the finite-dimensional case, i.e., if $\calP, \calS \in \C^{n \times n}$, it is easy to determine whether $\calL$ is self-adjoint or (maximally) dissipative. In this case, there are symmetry and maximality (i.e., rank) conditions on the matrices $\calP, \calS$ that are straightforward to check for particular applications, cf. \cite{GernHR21}. In infinite dimensions, determining the (possibly multiple) domains $\dom \left[\begin{smallmatrix} \calP \\ \calS \end{smallmatrix}\right]$ 
 for which $\calL$ is self-adjoint or (maximally) dissipative, respectively, is generally not straightforward. In this paper we address the problem of characterizing these domains using extension theory for linear relations and the notion of boundary triplets having its foundations in \cite{gorbachuk1989extension, derkach1991generalized}. Roughly speaking, a boundary triplet $(\calG, \Gamma_0, \Gamma_1)$ consists of a Hilbert space $\calG$ and two maps $\Gamma_0, \Gamma_1: \calL \to \calG \times \calG$ that are intertwined with the (multi-valued) operator $\calL$ by means of an abstract Green identity, cf. \cite{BehrHdS20}. The boundary maps $\Gamma_0, \Gamma_1$ then may be directly used to formulate suitable boundary conditions of~\eqref{eq:eq1} or to enable boundary control and observation. 
 
 Our first main contribution is to 
 provide a novel definition of boundary triplets where the boundary maps $\Gamma_0, \Gamma_1$ are defined on $\dom \left[\begin{smallmatrix} \calP \\ \calS \end{smallmatrix}\right]$ instead of the whole relation $\calL$. Having found such a boundary triplet we are able to characterize all self-adjoint and maximal dissipative realizations of $\calL=\ran\left[\begin{smallmatrix} \calP \\ \calS \end{smallmatrix}\right]$. Although the classical boundary triplet framework from \cite{BehrHdS20, behrndt2007boundary} is also applicable for $\calL$, this requires a definition of the boundary evaluation maps directly on $\calL=\ran \left[\begin{smallmatrix}
       \calP \\ \calS
   \end{smallmatrix}\right]$. This means for \eqref{eq:pHs} that the imposed boundary conditions are expressed in terms of $\calP \lambda $ and $ \calS \lambda$ instead of $\lambda$ itself, which makes the formulation cumbersome. Our suggested distinct notion of boundary triplets in which the boundary maps are formulated directly on the domain of the combined operator $\left[\begin{smallmatrix}
       \calP \\ \calS
   \end{smallmatrix}\right]$ allows us to overcome this issue. \par
As our second result we apply the abstract extension theory to a class of implicit port-Hamiltonian systems (pHs) on one-dimensional spatial domains. 
In this paper, we focus on the implicit PDE system 
\begin{align}\label{eq:pHs}
    \frac{\mathrm{d}}{\mathrm{d}t} \underbrace{\sum\limits_{k=0}^N \tdif{k} P_k \tdif{k} }_{=\calP} x(t, \xi) =\displaystyle \underbrace{\sum\limits_{k=0}^M J_k \tdif{k} }_{=\calJ} \Big( \underbrace{\sum\limits_{l=0}^N \tdif{l} S_l \tdif{l} }_{=\calS} x(t,\xi) \Big) 
\end{align}
where $x(t,\xi) \in \C^n$ for any time $t \geq 0$ and spatial variable $\xi \in (a,b)\subset \mathbb{R}$ and where $P_k, J_k, S_l \in \C^{n\times n}$. In \cite{MascvdSc23} the class \eqref{eq:pHs} was approached by the inclusions (often called the \emph{geometric notion} of port-Hamiltonian systems)
\begin{align}\label{eq:pHs_dyn}
    (z(t), e(t)) \in \ran \left[\begin{smallmatrix}
       \calP \\ \calS
   \end{smallmatrix}\right]\subset \calX \times \calX, \qquad (e(t), \tfrac{\mathrm{d}}{\mathrm{d}t} z(t)) \in \gr \calJ \subset \calX \times \calX,
\end{align}
where $\calX=L^2(a,b;\C^n)$. Here, $\calL=\ran \left[\begin{smallmatrix}
       \calP \\ \calS
   \end{smallmatrix}\right]$ is called \emph{Lagrangian subspace} if it is self-adjoint and $\calD=\gr \calJ$ is a \emph{Dirac structure} if $\calJ$ is skew-adjoint. In this case we call the pair $(\calD, \calL)$ a generalized pHs, see also \cite{van2018generalized, bendimerad2024stokes}. 
   We define a boundary triplet for both $\calD$ and $\calL$ to characterize the domains $\dom \calJ$ and $\dom \left[\begin{smallmatrix}
       \calP \\ \calS
   \end{smallmatrix}\right]$ for which $(\calD, \calL)$ is a generalized pHs. Therein, the boundary maps may be explicitly derived and only depend on the coefficient matrices $P_k, J_k$ and $S_l$. Our general framework allows non-invertible leading coefficients in each of the differential operators. This generalizes previous characterizations of 
    maximally dissipative extensions from \cite{le2005dirac,JacoZwar12,villegas2007port}.
   As the governing operator $\calJ$ in the pHs formulation \eqref{eq:pHs} is assumed to be skew-symmetric, we recall the notion of boundary triplets for skew-symmetric operators, which were recently considered in \cite{wegner2017boundary, trostorff2023m,skrepek2021linear}.  
   
The implicit PDEs \eqref{eq:pHs} can be used to model various physical systems, such as the Dzektser equation, the wave equation, the Euler-Bernoulli beam, and non-local elasticity phenomena, cf. \cite{bendimerad2024stokes,heidari2019port,roetman1967biharmonic,MascvdSc23,JacoZwar12}. It serves as a natural extension to existing port-Hamiltonian formulations of PDEs on one-dimensional domain. The case of $N=0$ and $P_0=I$ \eqref{eq:pHs} was thoroughly studied in \cite{villegas2007port,JacoZwar12} by means of a semigroup approach. An extension to higher-order spatial derivatives was considered in \cite{augner2016stabilisation}, higher-dimensional spatial domains were included in \cite{kurula2015linear,skrepek2021linear}, a formulation via abstract dissipative block operators was given in \cite{gernandt2024stability}, and a system node approach was introduced in~\cite{PhilReis23}. Since the operator on the left-hand side is not assumed to be invertible, \eqref{eq:pHs} corresponds to a linear partial differential-algebraic equation (PDAE). Recent works on PDAEs, in particular in the context of pHs include \cite{jacob2022solvability, erbay2024weierstra,gernandtreis2023} focusing particularly on solution theory. Beyond that, dissipative Hamiltonian DAEs which are closely related to maximally dissipative relations were studied in \cite{mehrmann2023abstract,mehl2024spectral}. Approaches leveraging extension theory and boundary triplets for pHs include \cite{gernandt2024stability} showing stability and passivity properties that are preserved under Kirchhoff-type interconnections and \cite{kurula2010dirac} where the authors gave a relationship between Dirac structures, unitary operators in Kre\u{\i}n spaces and boundary triplets.
   \par
   This paper is organized as follows. After recalling some notations and preliminaries on linear relations and function spaces, we develop the extension theory for range representations of linear relations via boundary triplets in~\ref{sec:ext_theory}. The construction of the boundary triplet is exemplified by means of the Dzektser equation and we provide a one-to-one correspondence between our approach and so-called Lagrangian subspaces, that are frequently used in geometrical port-Hamiltonian systems theory.  In Section~\ref{sec:pHs} we construct boundary triplets to range representations given by the matrix differential operators that define implicit pH systems~\eqref{eq:pHs}. We conclude this work in Section~\ref{sec:applications} with applications and the examples of the biharmonic wave equation and the elastic rod with non-local elasticity condition.
\par 
\noindent \textbf{Notation and preliminaries.} This section on notation is divided into two parts. The first part introduces the concept of linear relations, primarily following the definitions and conventions outlined in \cite{BehrHdS20}. The second part focuses on the function spaces that arise in the context of differential operators on one-dimensional spatial domains. \par  
Let $\calX_1, \calX_2,  \calX$ be Hilbert spaces over $\C$. A \textit{linear relation} from $\calX_1$ to $\calX_2$ is a linear subspace $H$ of the product space $\calX_1 \times \calX_2$. A linear relation from $\calX$ to $\calX$ is called linear relation in $\calX$. The product space $\calX_1 \times \calX_2$ is endowed with the canonical inner product $\langle \cdot , \cdot \rangle_{\calX_1\times \calX_2} := \langle \cdot , \cdot \rangle_{\calX_1} + \langle \cdot , \cdot \rangle_{\calX_2}$ with the corresponding norm $\|\cdot\|_{\calX_1 \times \calX_2} \coloneqq \langle \cdot , \cdot \rangle^{1/2}_{\calX_1 \times \calX_2}$. The \textit{domain, range, kernel} and \textit{multivalued part} of a relation $H \subset \calX_1 \times \calX_2$ are defined by
\begin{equation*}
    \begin{array}{rcl}
        \dom H & \coloneqq & \left\lbrace h_1 \in \calX_1 \, \middle| \, (h_1,h_2) \in H  \ \text{for some} \ h_2 \in \calX_2 \right\rbrace, \\ 
        \ran H & \coloneqq & \left\lbrace h_2 \in \calX_2 \, \middle| \, (h_1,h_2) \in H  \ \text{for some} \ h_1 \in \calX_2 \right\rbrace, \\ 
        \ker H & \coloneqq & \left\lbrace h_1 \in \calX_1 \, \middle| \, (h_1, 0) \in H \right\rbrace, \\ 
        \mul H & \coloneqq & \left\lbrace h_2 \in \calX_2 \, \middle| \, (0, h_2) \in H  \right\rbrace.
    \end{array}
\end{equation*}
The \textit{inverse} of the relation $H$ from $\calX_1$ to $\calX_2$ is
\begin{equation*}
    H^{-1} \coloneqq  \left\lbrace (h_2, h_1) \in \calX_2 \times \calX_1 \, \middle| \, (h_1, h_2) \in H \right\rbrace.
\end{equation*}
For $\lambda \in \C$ we set
\begin{equation*}
    \lambda H \coloneqq  \left\lbrace (h_1, \lambda h_2) \in \calX_1 \times \calX_2 \, \middle| \, (h_1, h_2) \in H \right\rbrace.
\end{equation*}
The \textit{product} of two linear relations $H,K$ in $\calX$ is 
\begin{equation*}
    HK=\left\lbrace (f,h) \in \calX \times \calX \, \middle| \, \exists \, g \in \calX \, \text{s.t.} \,  (f,g) \in K, (g,h) \in H \right\rbrace.
\end{equation*}
The \textit{sum} of relations $H,K$ in $\calX$ is 
\begin{align*}
    H+K \coloneqq \left\lbrace (f,h+k) \in \calX \times \calX \, \middle| \, (f,h) \in H, (f,k) \in K \right\rbrace
\end{align*}
and we abbreviate
\begin{align*}
    H+\mu \coloneqq H+\mu I= \left\lbrace (f,h+\mu f) \in \calX \times \calX \, \middle| \, (f,h) \in H \right\rbrace
\end{align*}
for $\mu \in \C$. The \textit{componentwise sum} of relations $H,K$ in $\calX$ is 
\begin{align*}
    H \, \widehat{+} \, K \coloneqq \left\lbrace (f+g,h+k) \in \calX \times \calX \, \middle| \, (f,h) \in H, (g,k) \in K \right\rbrace.
\end{align*}
The \textit{adjoint relation} of a linear relation $H$ from $\calX_1$ to $\calX_2$ is given by
\begin{equation*}
    H^* \coloneqq   \left\lbrace (g_1,  g_2) \in \calX_2 \times \calX_1 \, \middle| \, \langle g_2 , h_1 \rangle_{\calX_1}=\langle g_1 , h_2 \rangle_{\calX_2}   \ \text{for all} \  (h_1, h_2) \in H \right\rbrace
\end{equation*}
where it is straightforwardly verified that $(H^*)^{-1}=(H^{-1})^*\eqqcolon H^{-*}$.

We define the sesquilinear form 
\begin{equation*}
    \begin{array}{rccl}
   \mform{\cdot}{\cdot}{\calX}:& (\calX \times \calX) \times (\calX \times \calX) & \to & \C, \\ & ((h_1,h_2)
    ,  (g_1,g_2))
    & \mapsto&  \langle h_1 , g_2 \rangle_{\calX}-\langle h_2 , g_1 \rangle_{\calX}  .
    \end{array}
\end{equation*}
For a linear relation $H \subset \calX \times \calX$ the orthogonal complement with respect to the indefinite inner product $\mform{\cdot}{\cdot}{\calX}$ is defined by
\begin{equation*}
    H^{\mperp} = \left\lbrace (g_1,  g_2) \in \calX \times \calX \, \middle| \, \mform{(h_1,h_2)}{(g_1,g_2)}{\calX} =0  \ \text{for all} \  (h_1, h_2) \in H \right\rbrace.
\end{equation*}
Again, a simple calculation yields
\begin{equation}\label{eq:mperp_eq_star}
    H^{\mperp}=H^*.
\end{equation}
We call a relation $H$ in $\calX$ \textit{symmetric} if $H \subset H^*$ and \textit{self-adjoint} if $H = H^*$. Correspondingly, a relation $H$ is \textit{skew-symmetric} if $H \subset -H^*$ and \textit{skew-adjoint} if $H = -H^*$. Let $H$ be a (skew-) symmetric linear relation, i.e., $\pm H \subset H^*$. A linear relation $K$ in $\calX$ is an \textit{intermediate extension} of $H$ if $\pm H \subset K \subset H^*$. Clearly, if $K$ is an intermediate extension of $H$, then also $\pm K^*$ is.
A relation $H$ in $\calX$ is \textit{dissipative} if $\re\langle g,f\rangle_\calX \leq 0$ holds for all $(f,g)\in H$. Moreover, $H$ is \textit{maximally dissipative} if $H$ is dissipative and for all dissipative relations $K$ in $\calX$ with $H \subset K$ it holds that $H=K$. Note that the name dissipativity is not consistently used in the literature, e.g.\ in \cite{BehrHdS20} dissipativity requires the imaginary part to be positive. 
Moreover, a relation $H$ from $\calX_1$ to $\calX_2$ is \textit{closed} if $H$ is a closed subspace of $\calX_1 \times \calX_2$, i.e., for all sequences $(x_n)_{n \in \N} \subset H$, the convergence $x_n \to x$ in $\calX_1 \times \calX_2$ implies $x \in H$. The \textit{closure} $\overline{H}$ of a relation $H$ is the set of all limit points $x$ of sequences $(x_n)_{n \in \N} \subset H$. By definition, we have $\overline{H}=H^{**}$. 

A bounded linear operator $R: \calX_1 \to \calX_2$ is called \textit{unitary} if $R$ satisfies $\ran R = \calX_2$ and $\langle Rx,Ry \rangle_{\calX_2}=\langle x,y \rangle_{\calX_1}$ for all $x,y \in \calX_1$. For a bounded operator (or matrix) $R: \calX_1 \to \calX_2$ a \textit{singular value} $\sigma_i(R)=\sqrt{\lambda_i(R^*R)}$ is the square root of an eigenvalue of $R^*R$. We denote by $\sigma_{\min}(R)$ the \textit{smallest singular value} of $R$. For an unbounded operator $A: \calX_1 \supset \dom A \to \calX_2$ we call $\dom A$ the \textit{domain} of $A$ and $\ker A=\left\lbrace x \in \dom A \, \middle| \, Ax =0 \right\rbrace$ the \textit{kernel} of $A$. We note that $\dom A$ and $\ker A$ coincide with the domain and kernel of the linear relation 
\begin{align*}
    \gr A \coloneqq \left\lbrace (x,Ax) \, \middle| \, x \in \dom A \right\rbrace.
\end{align*}
When it is clear from the context, we will sometimes write $A$ for the linear relation $\gr A$ and we set
\begin{align*}
    \ran A \coloneqq \left\lbrace Ax \, \middle| \, x \in \dom A \right\rbrace \subset \calX_2.
\end{align*}

We utilize the following standard function spaces, cf.~\cite{adams2003sobolev}. 
For an open interval $(a,b) \subset \R$, $L^2(a,b;\C^n)$ is the space of square integrable functions from $(a,b)$ with values in $\C^n$ which is a Hilbert space when endowed with the standard inner product. Moreover, the space of smooth compactly supported function with values in $\C^n$ is denoted by $C_c^\infty(a,b;\C^n)$ and the space of $N$-times weakly differentiable functions is $H^N(a,b;\C^n)$ endowed with the usual norm. We let $H_0^N(a,b;\C^n)$ be the closure of $C_c^\infty(a,b;\C^n)$ with respect to $||\cdot||_{H^N(a,b;\C^n)}$.

\section{Extension theory for parametric representations of linear relations}\label{sec:ext_theory}
In this section, we introduce a novel class of boundary triplets for symmetric relations in range representations and use this to characterize the set of intermediate extensions of the considered relation. Then,  in Section~\ref{subsec:dzektser}, the approach is applied to the Dzektser equation. In Section~\ref{subsec:Lagrange} a one-to-one correspondence is given between the suggested boundary triplet definition and the so-called Stokes-Lagrange subspaces. In Section~\ref{subsec:bdd_triplets_skew} we conclude with a short recapitulation of boundary triplets for skew-symmetric operators.

Let $\calX$ be a Hilbert space. 
Throughout this section, we consider two
densely defined linear operators 
\begin{equation}
\label{def:Q0T0}
    \calP_0: \calX \supset \dom \calP_0 \to \calX, \qquad \calS_0: \calX \supset \dom \calS_0 \to \calX
\end{equation}
with adjoints $\calP \coloneqq \calP_0^*$ and $\calS \coloneqq \calS_0^*$ and the column operator 
\begin{align}
\label{def:PS}
    \begin{bmatrix}
       \calP \\ \calS
   \end{bmatrix} : \calX \supset \dom \begin{bmatrix}
       \calP \\ \calS
   \end{bmatrix} \to \calX \times \calX, \qquad \begin{bmatrix}
       \calP \\ \calS
   \end{bmatrix}  x = \begin{bmatrix}
       \calP x \\ \calS x
   \end{bmatrix} 
\end{align}
for all  $x \in  \dom 
\left[\begin{smallmatrix}
        \calP \\ \calS
\end{smallmatrix}\right] \coloneqq \dom \calP \, \cap \, \dom \calS$. 

In \eqref{eq:pHs_dyn} it was observed that the range of the operator $\left[\begin{smallmatrix}
        \calP \\ \calS
\end{smallmatrix}\right]$ naturally occurs in the class of implicit port-Hamiltonian systems, cf. \cite{MascvdSc23}. In particular, the above operator range can be used to formulate degenerate Cauchy problems and to investigate their existence of solutions, \cite{knuckles1994remarks, arendt2023semigroups}.

\subsection{Main result: Self-adjointness and dissipativity via boundary triplets}\label{subsec:bdd_triplets}

\begin{ass}\label{ass:ass1}
We assume that the operators $\calP_0$ and $\calS_0$ given by \eqref{def:Q0T0} are such that $\left[\begin{smallmatrix}
       \calP \\ \calS
   \end{smallmatrix}\right]$ given by \eqref{def:PS} is \emph{coercive}, meaning that for some $c>0$ it holds  
   \begin{align}\label{eq:coercivity}
       \left\|\begin{bmatrix}
       \calP \\ \calS
   \end{bmatrix}x \right\|_{\calX \times \calX} \geq c \|x\|_\calX\quad \text{for all} \ x \in \dom\begin{bmatrix}
       \calP \\ \calS
   \end{bmatrix}.
   \end{align}
\end{ass}
\noindent It is clear that Assumption~\ref{ass:ass1} is fulfilled as soon as $\mathcal{P}$ or $\mathcal{S}$ is boundedly invertible. However, we will show by means of the Dzektser equation treated in Subsection~\ref{subsec:dzektser}, how this coercivity of the combined operator may be verified for a broader class of (truly implicit) examples without invertibility of one of the operators.

We introduce our central definition of a boundary triplet for $\ran \left[\begin{smallmatrix} \calP \\ \calS \end{smallmatrix}\right]$.
\begin{defn}
\label{def:range_triplet}
    We define a boundary triplet for $\ran \left[\begin{smallmatrix}
        \calP \\ \calS
    \end{smallmatrix}\right]$ by a triplet $(\calG, \Gamma_0, \Gamma_1)$ consisting of a Hilbert space $\calG$ and a linear and surjective map $\Gamma = \left[ \begin{smallmatrix}
        \Gamma_0 \\ \Gamma_1
    \end{smallmatrix} \right]: \dom \left[\begin{smallmatrix}
        \calP \\ \calS
    \end{smallmatrix}\right] \to \calG \times \calG$ such that the \textit{abstract Green identity} 
    \begin{equation}\label{eq:Green}
        \left\langle \calS x , \calP y \right\rangle_\calX - \left\langle \calP x , \calS y \right\rangle_\calX = \left\langle \Gamma_1 x, \Gamma_0 y \right\rangle_\calG -\left\langle \Gamma_0 x, \Gamma_1 y\right\rangle_\calG 
    \end{equation}
    holds for all $x,y \in \dom \left[\begin{smallmatrix}
        \calP \\ \calS
    \end{smallmatrix}\right]$. The space $\calG$ is called the \textit{boundary space}. 
\end{defn}
We comment on the relation to (ordinary) boundary triplets in Remark~\ref{rem:relation_to_jussi}. \par
The first main result of this work is the following characterization of crucial properties of the operator by means of the relation in the boundary space.

 \begin{thm}\label{thm:self-adj_extension}
        Let Assumption~\ref{ass:ass1} hold, let $(\calG, \Gamma_0, \Gamma_1)$ be a boundary triplet for $\ran \left[\begin{smallmatrix}
        \calP \\ \calS
    \end{smallmatrix}\right]$, let $\Theta \subset \calG \times \calG$, and define an operator $A_\Theta: \calX \supset \dom A_\Theta \to \calX \times \calX$ by
    \begin{equation}
            A_\Theta = \left[\begin{smallmatrix}
        \calP \\ \calS
    \end{smallmatrix}\right] , \qquad \dom A_\Theta = \left\lbrace  x \in \dom \left[\begin{smallmatrix}
        \calP \\ \calS
    \end{smallmatrix}\right] \, \middle| \, \Gamma x \in \Theta \right\rbrace.
    \end{equation}
    Then, $H_\Theta = \ran A_\Theta \subset \calX \times \calX$ is self-adjoint (maximally dissipative) if and only if $\Theta \subset \calG \times \calG$ is self-adjoint (maximally dissipative).
    \end{thm}
The previous corollary enables us to derive the set of boundary conditions of \eqref{eq:pHs_dyn} that lead to a port-Hamiltonian system representation. The proof of Theorem \ref{thm:self-adj_extension} is the main objective of the subsequent Subsection~\ref{subsec:proof}. Before proceeding, we briefly comment on the relation to the existing notion of boundary triplets for linear relations.
\begin{rem}\label{rem:relation_to_jussi}
    The definition of a boundary triplet provided in Definition~\ref{def:range_triplet} is closely related to the usual definition of a boundary triplet for the adjoint of a linear relation, cf. \cite[Def. 2.1.1]{BehrHdS20}. Normally, such an ordinary boundary triplet for the adjoint $H^*$ of a symmetric relation $H$ is defined via $(\calG, \Gamma_0, \Gamma_1)$ with $\calG$ being a Hilbert space and $\Gamma =  \left[ \begin{smallmatrix}
        \Gamma_0 \\ \Gamma_1
    \end{smallmatrix} \right]: H^* \to \calG \times \calG$ linear, surjective and satisfying the abstract Green identity
    \begin{equation}\label{eq:Green_for_relations}
        \left\langle f' , g \right\rangle_\calX - \left\langle f , g' \right\rangle_\calX = \left\langle \Gamma_1 \left[ \begin{smallmatrix}
            f \\ f'
        \end{smallmatrix} \right], \Gamma_0\left[ \begin{smallmatrix}
             g \\ g'
        \end{smallmatrix} \right] \right\rangle_\calG -\left\langle \Gamma_0 \left[ \begin{smallmatrix}
            f \\ f'
        \end{smallmatrix} \right], \Gamma_1 \left[ \begin{smallmatrix}
            g \\ g'
        \end{smallmatrix} \right]\right\rangle_\calG 
    \end{equation}
    for all $\left[ \begin{smallmatrix}
            f \\ f'
        \end{smallmatrix} \right], \left[ \begin{smallmatrix}
            g \\ g'
        \end{smallmatrix} \right] \in H^*$. Instead of defining the boundary map $\Gamma$ on the whole relation $H^*$, we define it on $\dom \left[\begin{smallmatrix}
        \calP \\ \calS
    \end{smallmatrix}\right]$. When $\calP=I$, Definition~\ref{def:range_triplet} aligns precisely with the standard definition of a boundary triplet for the adjoint $\mathcal{S}$ of a symmetric operator. We present a further remark on this topic later, where we provide more details.
\end{rem}

\subsection{Proof of Theorem~\ref{thm:self-adj_extension}}\label{subsec:proof}

   The following lemma gathers properties of $\left[\begin{smallmatrix}
       \calP \\ \calS
   \end{smallmatrix}\right]$ that will be central in our subsequent analysis.
   \begin{lem}\label{lem:injective_and_closed_range}
       Let Assumption \ref{ass:ass1} hold. Then the operator $\left[\begin{smallmatrix}
       \calP \\ \calS
   \end{smallmatrix}\right]$ is closed, injective and $\ran \left[\begin{smallmatrix}
       \calP \\ \calS
   \end{smallmatrix}\right]$ is closed in $\calX \times \calX$.
   \end{lem}
   \begin{proof}
According to~\cite[Prop. 3.1]{moller2008adjoints}, the adjoint of the densely defined operator $\begin{bmatrix}
            \calP_0 & \calS_0
        \end{bmatrix}$ is the combined operator
        \begin{equation*}
        \left[\begin{smallmatrix}
         \calP \\\calS
    \end{smallmatrix}\right]: \calX \supset \dom\calP \cap \dom\calS \to \calX \times \calX.
    \end{equation*} 
    Therefore, $\left[\begin{smallmatrix}
       \calP \\ \calS
   \end{smallmatrix}\right]$ is an adjoint operator and therefore closed.
    The coercivity condition~\eqref{eq:coercivity} immediately shows that $\left[\begin{smallmatrix}
       \calP \\ \calS
   \end{smallmatrix}\right]$ is injective with closed range by \cite[Thm.~IV.5.2]{kato2013perturbation}. 
   \end{proof}
   The following lemma is helpful for characterizing adjoint relations.
    \begin{lem}\label{lem:adjoint_rel}
        Let $A: \calX \supset \dom A \to \calX, B: \calX \supset \dom B \to \calX$ two densely defined linear operators and define their graphs as $K=\gr A, L=\gr B$. Then,
        \begin{align*}
            (K^{-1}L)^*=L^*K^{-*}
        \end{align*}
        if and only if $\ran  \left[ \begin{smallmatrix}
            A & B
    \end{smallmatrix}\right]$ is closed.
    \end{lem}
    \begin{proof}
        First of all, we provide a small auxiliary result. To this end, consider the row operator
    \begin{equation*}
        \begin{bmatrix}
            -B & A
        \end{bmatrix}: \calX \times \calX \supset  \dom B\times \dom A \to \calX 
    \end{equation*}
    defined by $\begin{bmatrix}
            -B & A
        \end{bmatrix} \left[ \begin{smallmatrix}
            f \\ g
        \end{smallmatrix}\right]=-B f + A g$ and note that $K^{-1}L=\ker \begin{bmatrix}
            -B & A
        \end{bmatrix}$. Then, by \cite[Prop. 3.1]{moller2008adjoints} the Hilbert space adjoint of $\begin{bmatrix}
            -B & A
        \end{bmatrix}$ is the combined operator
        \begin{equation*}
        \left[\begin{smallmatrix}
         -B^* \\ A^*
    \end{smallmatrix}\right]: \calX \supset \dom A^* \cap \dom B^* \to \calX \times \calX.
    \end{equation*}
    Assume that $\ran  \left[ \begin{smallmatrix}
            A & B
    \end{smallmatrix}\right] $ is closed. Following \cite[Def. 1.3.1]{BehrHdS20}, we observe that $(K^{-1}L)^*=J(K^{-1}L)^\perp$ where $J: \calX^2 \to \calX^2$ denotes the \textit{flip-flop operator} defined by $J\left[ \begin{smallmatrix}
            f \\ g
        \end{smallmatrix}\right]=\left[ \begin{smallmatrix}
            g \\ -f
        \end{smallmatrix}\right]$. Hence,
    \begin{equation*}
        (K^{-1}L)^*
        = J \left( \ker \left[ \begin{smallmatrix}
            -B & A
    \end{smallmatrix}\right]\right)^\perp \stackrel{(*)}{=} J \ran \left[\begin{smallmatrix}
         -B^* \\ A^*
    \end{smallmatrix}\right]=\ran \left[\begin{smallmatrix}
         A^* \\ B^*
    \end{smallmatrix}\right]=(\gr B)^*(\gr A)^{-*},
    \end{equation*}
    where in $(*)$ we applied the closed range theorem, \cite[Chap. 7.5]{yosida2012functional} and the representation of the adjoint of $\left[\begin{smallmatrix}
            -B & A
         \end{smallmatrix}\right]$ from the beginning of this proof. To prove the converse direction we see that $(K^{-1}L)^*=\ran \left[\begin{smallmatrix}
         A^* \\ B^*
    \end{smallmatrix}\right]$ is closed as it is an adjoint relation. The closed range theorem and the fact that $\left[ \begin{smallmatrix}
            A & B
    \end{smallmatrix}\right]^*=\left[\begin{smallmatrix}
         A^* \\ B^*
    \end{smallmatrix}\right]$ yields that $\ran \left[ \begin{smallmatrix}
            A & B
    \end{smallmatrix}\right] $ is closed in $\calX$. 
    \end{proof}
    We now define linear relations given in kernel and range representation involving the previously defined operators.
   \begin{cor}\label{cor:H}
    Let Assumption~\ref{ass:ass1} hold and define the linear relation
\begin{equation}\label{eq:H}
    H \coloneqq (\gr \calP_0)^{-1} (\gr \calS_0)=
    \left\lbrace \left[ \begin{smallmatrix}
        f_1 \\ f_2
    \end{smallmatrix} \right] \in \dom \calS_0 \times \dom \calP_0 \, \middle| \, 
        \calS_0 f_1= \calP_0 f_2
    \right\rbrace  .
\end{equation}
Then, the adjoint relation $H^*$ is given by
    \begin{equation}\label{eq:H_star}
        H^*= (\gr \calS)(\gr \calP)^{-1} =\left\lbrace  \begin{bsmallmatrix}
        \calP \\ \calS 
    \end{bsmallmatrix} \lambda  \, \middle| \, \lambda \in  \dom \left[ \begin{smallmatrix}
        \calP \\ \calS
    \end{smallmatrix}\right] \right\rbrace=\ran \left[ \begin{smallmatrix}
        \calP \\ \calS
    \end{smallmatrix} \right]. 
    \end{equation}
\end{cor}
\begin{proof}
    As Lemma \ref{lem:injective_and_closed_range} yields the closedness of $\ran \left[ \begin{smallmatrix}
        \calP \\ \calS
    \end{smallmatrix} \right]$, we apply Lemma \ref{lem:adjoint_rel} with $A=\calP_0, B=\calS_0$.
\end{proof}
In the subsequent analysis, we let $H$ be given as in \eqref{eq:H}. 

\medskip\noindent
\textbf{Short interlude.} We briefly investigate the symmetry behavior of $H$. At first glance, it not easy to check $H \subset H^*$ in general because $H$ is defined as a kernel of an unbounded row operator. 
In this part, we derive a more verifiable characterization of symmetry of $H$ which will be used in the application to port-Hamiltonian systems in Section~\ref{sec:pHs}. First, for a given $x \in \dom H$ we define $H(x)\coloneqq \left\lbrace y \in \calX \, \middle| \, \left( \begin{smallmatrix}
       x \\ y
    \end{smallmatrix} \right) \in H \right\rbrace$. The linear relation $H$ is uniquely determined by its graph, also denoted by
    \begin{align*}
        H =  \left\lbrace \left( \begin{smallmatrix}
       x \\ y
    \end{smallmatrix} \right) \in \calX \times \calX \, \middle| \, x \in \dom H , \,  y \in H(x) \right\rbrace.
    \end{align*}
    By identifying $H$ with its graph it is clear that $H$ is symmetric, i.e., $H \subset H^*$, if and only if $\dom H \subset \dom H^*$ and $H(x) \subset H^*(x)$ for all $x \in \dom H$. 
    \par 
    To proceed, we define operators $N_0 = I_{\ran \calP_0} \calS_0   \calP$ and $M=\calP_0   \calS I_{\dom \calP}$. In detail, this means that
    \begin{equation*}
    \begin{array}{rclrcl}
         \dom N_0 &\coloneqq & \calP^{-1}\calS_0^{-1}\ran \calP_0  \qquad & N_0 x & =& \calS_0 \calP x ,\\
        \dom M &\coloneqq &   \calS^{-1}\dom \calP_0 \cap \dom \calP  \qquad & M x &=& \calP_0 \calS x  .
    \end{array}
    \end{equation*}
\begin{prop}\label{prop:H_symm}
    The linear relation $H$ defined in \eqref{eq:H} is symmetric if and only if the following two conditions are satisfied:
    \begin{itemize}
        \item[(i)] $\calS_0^{-1} \ran \calP_0 \subset \calP \dom \calS$;
        \item[(ii)] it holds that 
        \begin{align}\label{eq:symm_prop}
            \gr N_0 \subset \gr M \, \widehat{+} \, \left( \ker \calP \times \{ 0 \} \right).
        \end{align} 
    \end{itemize}
\end{prop}
Observe that \eqref{eq:symm_prop} means that for all $\lambda \in \dom N_0$ we have $\lambda = z + n$ where $z \in \dom M, n \in \ker \calP$ and it holds that
\begin{align}\label{eq:symm_prop_2}
     \calS_0  \calP z = \calP_0  \calS z.
\end{align}
\begin{proof}
    Before we start with the proof, observe that
\begin{equation*}
   \begin{array}{rclrcl}
       \dom H &=& \calS_0^{-1} \ran \calP_0 , \qquad  \qquad & \dom H^* &=& \calP \dom \calS, \\
       \mul H &=& \ker \calP_0  ,& \mul H^* &=& \calS \ker \calP.
   \end{array}
\end{equation*}
We first show that $(i)$ and $(ii)$ imply $H \subset H^*$. As previously discussed, it is enough to show $\dom H \subset \dom H^*$ and $H(x) \subset H^*(x)$ for all $x \in \dom H$. Clearly, $(i)$ is equivalent to $\dom H \subset \dom H^*$. 
Thus, it remains to show that $(ii)$ implies $H(x) \subset H^*(x)$ for all $x \in  \dom H$. To this end, let $y \in H(x)$ for an arbitrary but fixed $x \in \dom H$, i.e., $\calS_0 x = \calP_0y$. It follows from $(i)$ that $x=\calP \lambda$ for some $\lambda \in \dom \calP \cap \dom \calS$ which implies that $\lambda \in  \calP^{-1}\dom H = \dom N_0$. By $(ii)$, we have $\lambda = z + n$ where $z \in \dom M, n \in \ker \calP$ and it holds \eqref{eq:symm_prop_2}. Hence, 
\begin{align*}
     \calP_0 y =  \calS_0 x= \calS_0 \calP \lambda = \calP_0 \calS z.
\end{align*}
In other words, we have $y - \calS z \in \ker \calP_0$. Now, we conclude from \cite[Prop. 1.3.2]{BehrHdS20} that $(\dom H)^\perp = \mul H^*$ and $(\dom H^*)^\perp = \mul \overline{H}$ which implies
\begin{align*}
    \ker \calP_0 = \mul H \subset \mul \overline{H} = (\dom H^*)^\perp \subset (\dom H)^\perp = \mul H^* = \calS \ker \calP.
\end{align*}
Consequently, $y-\calS z \in \calS \ker \calP$ which means that $y-\calS z= \calS \eta$ for a $\eta \in \ker \calP$, that is $y = \calS (z+ \eta)$. Thus 
$$y \in H^*(x)=\left\lbrace \calS \xi \, \middle| \, x=\calP \xi \, \text{for some} \,  \xi \in \dom  \left[ \begin{smallmatrix}
        \calP \\ \calS
    \end{smallmatrix}\right]\right\rbrace$$ as $\eta \in \ker \calP$. \par
Conversely, let $H \subset H^*$ which clearly implies $(i)$. We directly obtain from the definition of products of linear relations that
\begin{align}\label{eq:relations_incl}
    (\gr \calP_0)H (\gr \calP) \subset (\gr \calP_0)H^* (\gr \calP).
\end{align}
Since $H=(\gr \calP_0)^{-1} (\gr \calS_0)$ and $(\gr \calP_0)(\gr \calP_0)^{-1}=\gr I_{\ran \calP_0}$, by \cite[Lem. 1.1.4]{BehrHdS20} we have 
\begin{align*}
    (\gr \calP_0)H (\gr \calP) &= (\gr \calP_0)(\gr \calP_0)^{-1} (\gr \calS_0)(\gr \calP ) \\
        &=(\gr I_{\ran \calP_0})(\gr \calS_0 \calP) \\
        &= \gr N_0.
\end{align*}
Again by \cite[Lem. 1.1.4]{BehrHdS20} we get
\begin{align*}
        (\gr  \calP)^{-1}(\gr  \calP) = \gr I_{\dom \calP} \, \widehat{+} \, \left( \{ 0 \} \times \ker \calP ) \right)
\end{align*}
Moreover, the fact $H^*=(\gr \calS)(\gr \calP)^{-1}$ yields
\begin{align*}
    (\gr \calP_0)H^* (\gr \calP )  &= (\gr \calP_0)(\gr \calS)(\gr \calP)^{-1} (\gr \calP ) \\
        &= (\gr \calP_0 \calS)\left(\gr I_{\dom 
        \calP } \, \widehat{+} \, \left( \{ 0 \} \times \ker \calP ) \right) \right) \\
        &= \left\lbrace \left[ \begin{smallmatrix}
       x \\ \calP_0 \calS y
    \end{smallmatrix}\right] \, \middle| \, x \in \dom \calP , y= x+ n, n \in \ker \calP \right\rbrace \\
        &= \gr M \, \widehat{+} \, \left( \ker \calP \times  \{ 0 \}  \right).
\end{align*}
Combining the last equalities with \eqref{eq:relations_incl} we obtain \eqref{eq:symm_prop} which shows the claim. 
\end{proof}
\par 
With this brief detour complete, we now return to our main discussion. In the following lemma, we equip $\dom \left[ \begin{smallmatrix}
        \calP \\ \calS
    \end{smallmatrix}\right]$ with a (natural) inner product. 
\begin{lem}\label{lem:dom_HS}
    Let Assumption~\ref{ass:ass1} hold, then the space $\dom  \left[ \begin{smallmatrix}
        \calP \\ \calS
    \end{smallmatrix}\right]$ endowed with
    \begin{align}\label{eq:dom_inner_product}
        \langle x,y \rangle_{\dom  \left[ \begin{smallmatrix}
        \calP \\ \calS
    \end{smallmatrix}\right]} = \langle \calP x, \calP y \rangle_\calX+ \langle \calS x, \calS y \rangle_\calX
    \end{align}
    is a Hilbert space.
\end{lem}
\begin{proof}
We define a sesquilinear form $\mathfrak{t}:\dom\mathfrak{t} \times \dom \mathfrak{t} \to \C$ via
\begin{align}
    \mathfrak{t}(x,y)=\left\langle Tx,Ty\right\rangle_\calZ
\end{align}
where $\dom \mathfrak{t}=\dom \left[ \begin{smallmatrix}
        \calP \\ \calS
    \end{smallmatrix}\right]$, $T=\left[ \begin{smallmatrix}
        \calP \\ \calS
    \end{smallmatrix}\right]$ and $\calZ=\calX \times \calX$. As $T=\left[ \begin{smallmatrix}
        \calP \\ \calS
    \end{smallmatrix}\right]$ is closed due to Lemma~\ref{lem:injective_and_closed_range}, \cite[Lem. 5.1.21]{BehrHdS20} implies that $\mathfrak{t}$ is a closed form. Moreover, Assumption~\ref{ass:ass1} yields the existence of $c>0$ such that $\mathfrak{t}(x,x)=\|Tx\|_\calZ^2=\left|\left| \left[\begin{smallmatrix}
       \calP \\ \calS
   \end{smallmatrix}\right]x \right|\right|_{\calX \times \calX} \geq c^2\|x\|_\calX^2$, that is, $\mathfrak{t}$ is bounded from below with a positive lower bound. This implies that $(\dom t, t)$ is a Hilbert space by \cite[Lem. 5.1.9]{BehrHdS20} and thus the claimed result.
\end{proof}
To prepare the following steps we note that the abstract Green identity \eqref{eq:Green} is equivalent to 
\begin{align}\label{eq:imaginary_green}
    \pm i\mform{ \left[ \begin{smallmatrix}
        \calP  \\ \calS 
    \end{smallmatrix}\right] x }{ \left[ \begin{smallmatrix}
        \calP  \\ \calS 
    \end{smallmatrix}\right] y }{\calX}  = \pm i \mform{ \left[ \begin{smallmatrix}
        \Gamma_0  \\ \Gamma_1
    \end{smallmatrix}\right]x}{ \left[ \begin{smallmatrix}
         \Gamma_0  \\ \Gamma_1
    \end{smallmatrix}\right] y }{\calG} 
\end{align}
for all $x,y \in \dom \left[\begin{smallmatrix}
        \calP \\ \calS
    \end{smallmatrix}\right]$ by premultiplication with $\pm i$. 
The following lemma relates the inner product $ \langle \cdot ,\cdot \rangle_{\dom \left[\begin{smallmatrix}
         \calP \\ \calS
     \end{smallmatrix}\right]} $ with the indefinite inner product $\mform{\cdot}{\cdot}{\calX}$. 
\begin{lem}\label{lem:defect_mform}
    Let $x \in \dom \left[\begin{smallmatrix}
         \calP \\ \calS
     \end{smallmatrix}\right]$ and $y \in \ker (\calS - (\pm i \calP))$. Then,
     \begin{align*}
         \langle x,y \rangle_{\dom \left[\begin{smallmatrix}
         \calP \\ \calS
     \end{smallmatrix}\right]} 
    = \pm i \mform{ \left[ \begin{smallmatrix}
        \calP \\ \calS 
    \end{smallmatrix}\right]x}{ \left[ \begin{smallmatrix}
        \calP  \\ \calS 
    \end{smallmatrix}\right]  y}{\calX}
     \end{align*}
\end{lem}
\begin{proof}
    Note that $y \in \ker (\calS - (\pm i \calP))$ implies $\calS y = \pm i \calP y$, which yields
     \begin{align*}
         \langle x,y \rangle_{\dom \left[\begin{smallmatrix}
         \calP \\ \calS
     \end{smallmatrix}\right]} 
     & = \langle \calP x, \calP y \rangle_\calX + \langle \calS x,  \calS y \rangle_\calX   \\
     & = \langle \calP x, \mp i\calS y \rangle_\calX + \langle \calS x, \pm i \calP y \rangle_\calX    \\
     & = \pm i \langle \calP x, \calS y \rangle_\calX  \mp i \langle \calS x,  \calP y \rangle_\calX \\
     & = \pm i \left( \langle \calP x, \calS y \rangle_\calX - \langle \calS x,  \calP y \rangle_\calX \right) \\
    & = \pm i \mform{ \left[ \begin{smallmatrix}
        \calP  \\ \calS 
    \end{smallmatrix}\right]x}{ \left[ \begin{smallmatrix}
        \calP  \\ \calS 
    \end{smallmatrix}\right]  y}{\calX}
     \end{align*}
    for all $x \in \dom \left[\begin{smallmatrix}
         \calP \\ \calS
     \end{smallmatrix}\right]$ and $y \in \ker (\calS - (\pm i \calP))$.
\end{proof}
In the subsequent lemma, we derive central properties of  boundary triplets for $\ran \left[\begin{smallmatrix}
        \calP \\ \calS
    \end{smallmatrix}\right]$ by adapting techniques from \cite[Prop. 2.1.2]{BehrHdS20} to our situation. The main difference is that our boundary map $\Gamma$ is defined directly on $\dom \left[\begin{smallmatrix}   \calP \\ \calS
    \end{smallmatrix}\right]$ instead of $H^*$.  
\begin{lem}\label{lem:Gamma_properties}
    Let Assumption~\ref{ass:ass1} hold, $H \subset H^*$ and let $(\calG, \Gamma_0, \Gamma_1)$ be a boundary triplet for $H^*=\ran \left[\begin{smallmatrix}
        \calP \\ \calS
    \end{smallmatrix}\right]$. Then, the following hold:
    \begin{enumerate}
        \item[(i)] $\Gamma = \left[ \begin{smallmatrix}
        \Gamma_0 \\ \Gamma_1
    \end{smallmatrix}\right]: \left( \dom \left[\begin{smallmatrix}
        \calP \\ \calS
    \end{smallmatrix}\right] , \left\| \cdot \right\|_{ \dom \left[\begin{smallmatrix}
        \calP \\ \calS
    \end{smallmatrix}\right] }\right)  \to \calG \times \calG$ is continuous;
   \item[(ii)]  $\left[\begin{smallmatrix}
        \calP \\ \calS
    \end{smallmatrix}\right] \ker \Gamma=\overline{H}$.
    \end{enumerate}
\end{lem}
\begin{proof}
\begin{itemize}
   \item[(i)] Recall that $\dom \left[\begin{smallmatrix}
        \calP \\ \calS
    \end{smallmatrix}\right] $ endowed with $\|\cdot\|_{\dom \left[\begin{smallmatrix}
        \calP \\ \calS
    \end{smallmatrix}\right]} $ is a Hilbert space by Lemma~\ref{lem:dom_HS}. We use the closed graph theorem, \cite[Chap. 7.5]{yosida2012functional}, to show the boundedness. Therefore, let $(x_n)_{n \in \N} \subset \dom \left[\begin{smallmatrix}
        \calP \\ \calS
    \end{smallmatrix}\right] $ such that
    \begin{equation*}
        x_n \stackrel{n \to \infty}{\to} x \in \dom \left[\begin{smallmatrix}
        \calP \\ \calS
    \end{smallmatrix}\right], \qquad  
    \Gamma x_n \stackrel{n \to \infty}{\to} \begin{bmatrix}
                \varphi_1 \\ \varphi_2
            \end{bmatrix} \in \calG \times \calG.
    \end{equation*}
    By surjectivity of $\Gamma$ we find for any $(\psi_1, \psi_2) \in \calG \times \calG$ an element $y \in \dom  \left[\begin{smallmatrix}
        \calP \\ \calS
    \end{smallmatrix}\right] $ with $\left[\begin{smallmatrix}
                \psi_1 \\ \psi_2
    \end{smallmatrix}\right] = \left[\begin{smallmatrix}
        -\Gamma_1 \\ \Gamma_0
    \end{smallmatrix}\right]y$.  Then, formula \eqref{eq:Green} implies
    \begin{align*}
            \left\langle \begin{bmatrix}
                \varphi_1 \\ \varphi_2
    \end{bmatrix} ,\begin{bmatrix}
                \psi_1 \\ \psi_2
    \end{bmatrix} \right\rangle_{\calG\times \calG}
            &= \lim\limits_{n \to \infty} \left( \left\langle \Gamma_1 x_n, \Gamma_0 y \right\rangle_\calG -\left\langle \Gamma_0 x_n, \Gamma_1 y\right\rangle_\calG \right) \\
            &= \lim\limits_{n \to \infty} \left( \left\langle \calS x_n , \calP y \right\rangle_\calX - \left\langle \calP x_n , \calS y \right\rangle_\calX  \right) \\
            &=  \left\langle \calS x , \calP y \right\rangle_\calX - \left\langle \calP x , \calS y \right\rangle_\calX \\
            &= \left\langle \Gamma_1 x, \Gamma_0 y \right\rangle_\calG -\left\langle \Gamma_0 x, \Gamma_1 y\right\rangle_\calG \\
            &= \left\langle \Gamma x ,  \begin{bmatrix}
                \psi_1 \\ \psi_2
    \end{bmatrix}\right\rangle_{\calG\times \calG}.
    \end{align*}
    Since $\psi_1$ and $\psi_2$ were chosen arbitrarily, we conclude 
$\Gamma x = \left[\begin{smallmatrix}
                \varphi_1 \\ \varphi_2
           \end{smallmatrix}\right]$, which yields the closedness of $\Gamma$ and thus its continuity. 
    \item[(ii)] We show both subset inclusions. Let $x \in \ker \Gamma$ and observe that for all $y \in \dom  \left[\begin{smallmatrix}
    	\calP \\ \calS
    \end{smallmatrix}\right] $ we have
    \begin{equation*}
    	0=\left\langle \Gamma_1 x, \Gamma_0 y \right\rangle_\calG -\left\langle \Gamma_0 x, \Gamma_1 y\right\rangle_\calG =  \left\langle \calS x , \calP y \right\rangle_\calX - \left\langle \calP x , \calS y \right\rangle_\calX .
    \end{equation*}
   Hence, $\left[\begin{smallmatrix}
   	\calP \\ \calS
   \end{smallmatrix}\right] x \in (H^{*})^{\mperp}=H^{**}=\overline{H}$. Conversely, let $x \in \left[\begin{smallmatrix}
   \calP \\ \calS
   \end{smallmatrix}\right]^{-1} \overline{H}$. Then, by surjectivity of $\Gamma $ for all $\left(\begin{smallmatrix}
   \varphi_1 \\ \varphi_2
   \end{smallmatrix}\right) \in \calG \times \calG$ there exists $y \in \dom \left[\begin{smallmatrix}
   \calP \\ \calS
   \end{smallmatrix}\right]$ such that $\left[ 
   	\begin{smallmatrix}
   	    \varphi_1 \\ \varphi_2    
   	\end{smallmatrix} \right]
   	 = \left[ \begin{smallmatrix}
   		-\Gamma_1 \\ \Gamma_0
   	\end{smallmatrix}\right]y$. This yields
   \begin{equation*}
   	\left\langle \Gamma x ,  \left[\begin{smallmatrix}
   		\varphi_1 \\ \varphi_2
   	\end{smallmatrix}\right] \right\rangle_{\calG^2} = \left\langle \Gamma_1 x, \Gamma_0 y \right\rangle_\calG -\left\langle \Gamma_0 x, \Gamma_1 y\right\rangle_\calG =  \left\langle \calS x , \calP y \right\rangle_\calX - \left\langle \calP x , \calS y \right\rangle_\calX =0
   \end{equation*}
   since $\left[\begin{smallmatrix}
   	\calP \\ \calS
   \end{smallmatrix}\right] y \in \ran \left[\begin{smallmatrix}
   \calP \\ \calS
   \end{smallmatrix}\right]= H^{*}$ and $\left[\begin{smallmatrix}
   \calP \\ \calS
   \end{smallmatrix}\right] x \in \overline{H}=(H^{*})^{\mperp}$. Consequently, $\Gamma x=0$ and hence $x \in \ker \Gamma$. \endproof
   \end{itemize}
\end{proof}
\noindent In the following, we will rely on a decomposition of $\dom \left[\begin{smallmatrix}
         \calP \\ \calS
     \end{smallmatrix}\right]$.
To this end, we introduce following subspaces associated with the eigenspace of a linear relation $T \subset \calX \times \calX$ at $\mu \in \C$
    \begin{equation*}
        \mathfrak{N}_\mu (T) = \ker (T-\mu), \qquad \widehat{\mathfrak{N}}_\mu(T)=\left\lbrace \begin{bmatrix}
                f_\mu \\ \mu f_\mu
            \end{bmatrix}  \, \middle| \, f_\mu \in \mathfrak{N}_\mu(T) \right\rbrace.
    \end{equation*}
The following result provides an explicit representation of the \emph{defect subspace} $\widehat{\mathfrak{N}}_\mu(H^*)$ at $\mu \in \C \setminus \R$.
\begin{lem}\label{lem:defect_representation}
    Consider the linear relation $H^*= \ran \left[\begin{smallmatrix}
         \calP \\ \calS
     \end{smallmatrix}\right]$. It holds that 
     \begin{align*}
     \mathfrak{N}_\mu (H^*)&=\calP\ker(\calS-\mu\calP) ,\\
         \widehat{\mathfrak{N}}_\mu(H^*) &= \left\lbrace \left[\begin{smallmatrix}
         \calP \\ \calS
     \end{smallmatrix}\right] x  \, \middle| \, x\in \ker (\calS - \mu  \calP) \right\rbrace = \ran \left. \left[\begin{smallmatrix}
         \calP \\ \calS
     \end{smallmatrix}\right] \right|_{\ker(\calS - \mu  \calP)}.
     \end{align*}
\end{lem}
\begin{proof}
   Let $\mu \in \C$ such that $\Im (\mu)\neq 0$ and consider
    \begin{equation*}
        \begin{array}{rcl}
            \mathfrak{N}_\mu (H^*) &=& \ker (H^*-\mu) \\
            &=& \ker \left\lbrace \left[ \begin{smallmatrix}
                h \\ h' - \mu h
            \end{smallmatrix} \right] \, \middle| \,\left[ \begin{smallmatrix}
                h \\ h' 
            \end{smallmatrix} \right] \in H^*=\ran \left[ \begin{smallmatrix}
        \calP \\ \calS
    \end{smallmatrix}\right] \right\rbrace  \\[1.5ex]
    &=& \ker \left\lbrace \left[ \begin{smallmatrix}
                \calP g \\ \calS g - \mu \calP g
            \end{smallmatrix} \right] \, \middle| \, g \in \dom  \left[ \begin{smallmatrix}
        \calP \\ \calS
    \end{smallmatrix}\right]  \right\rbrace 
         \\[1.5ex]
    &=& \calP\ker(\calS-\mu\calP).
        \end{array} 
    \end{equation*}
    Since the condition $x \in \ker(\calS - \mu \calP)$ is equivalent to $\calS x = \mu \calP x$, we may also write 
    \begin{equation*}
        \mathfrak{N}_\mu (H^*) = \left\lbrace \mu^{-1} \calS x  \, \middle| \, x \in \ker(\calS - \mu \calP) \right\rbrace = \ran \left. \mu^{-1} \calS  \right|_{\ker(\calS - \mu \calP)}.
    \end{equation*}
    However, this means
    \begin{align*}
        \widehat{\mathfrak{N}}_\mu(H^*)& = \left\lbrace \left[ \begin{smallmatrix}
                f_\mu \\ \mu f_\mu
            \end{smallmatrix} \right] \, \middle| \, f_\mu \in \mathfrak{N}_\mu(H^*) \right\rbrace =  \left\lbrace \left[ \begin{smallmatrix}
                \calP \\ \mu \mu^{-1} \calS 
            \end{smallmatrix} \right] \lambda \, \middle| \, \lambda \in \ker (\calS - \mu \calP) \right\rbrace \\
            &= \ran \left. \begin{bmatrix}
         \calP \\ \calS
     \end{bmatrix}\right|_{\ker(\calS - \mu  \calP)}.
    \end{align*}
\end{proof}
A consequence of Lemma~\ref{lem:defect_representation} is a \emph{von Neumann-type} decomposition. 
\begin{lem}\label{lem:orthogonal_decomposition}
Let Assumption~\ref{ass:ass1} hold, $H \subset H^*$ and let $(\calG, \Gamma_0, \Gamma_1)$ be a boundary triplet for $H^*=\ran \left[\begin{smallmatrix}
        \calP \\ \calS
    \end{smallmatrix}\right]$. If we endow $\dom \left[\begin{smallmatrix}
         \calP \\ \calS
     \end{smallmatrix}\right] $ with the inner product \eqref{eq:dom_inner_product}, then we have the following orthogonal decomposition
\begin{equation}
    \dom \left[\begin{smallmatrix}
         \calP \\ \calS
     \end{smallmatrix}\right] = \ker \Gamma \oplus \ker (\calS - i \calP) \oplus \ker (\calS + i \calP).
\end{equation}
\end{lem}
\begin{proof}
    Following \cite[Thm. 1.7.11]{BehrHdS20}, we conclude that 
    \begin{equation*}
        H^*=\overline{H} \oplus \widehat{\mathfrak{N}}_{i}(H^*) \oplus \widehat{\mathfrak{N}}_{-i}(H^*),
    \end{equation*}
    which can be rewritten in range representation with Lemma \ref{lem:Gamma_properties} (ii) and Lemma \ref{lem:defect_representation} as
    \begin{equation*}
        \ran \left[\begin{smallmatrix}
         \calP \\ \calS
     \end{smallmatrix}\right]= \ran \left. \left[\begin{smallmatrix}
         \calP \\ \calS
     \end{smallmatrix}\right]\right|_{ \ker \Gamma }\oplus \ran
     \left. \left[\begin{smallmatrix}
         \calP \\ \calS
     \end{smallmatrix}\right] \right|_{ \ker   (\calS - i \calP)}  \oplus \ran
     \left.  \left[\begin{smallmatrix}
         \calP \\ \calS
     \end{smallmatrix}\right]\right|_{  \ker   (\calS + i \calP) }.
    \end{equation*}
    Hence, for any $z\in \dom \left[\begin{smallmatrix}
         \calP \\ \calS
     \end{smallmatrix}\right]$, there are unique $x\in \ker \Gamma$ and $y_{\pm} \in \ker (\calS - (\pm i \calP))$ such that $$\left[\begin{smallmatrix}
         \calP \\ \calS
     \end{smallmatrix}\right]z = \left[\begin{smallmatrix}
         \calP \\ \calS
     \end{smallmatrix}\right]x + \left[\begin{smallmatrix}
         \calP \\ \calS
     \end{smallmatrix}\right]y_++ \left[\begin{smallmatrix}
         \calP \\ \calS
     \end{smallmatrix}\right]y_-.$$ 
     Moreover, due to the injectivity of $\left[\begin{smallmatrix}
         \calP \\ \calS
     \end{smallmatrix}\right]$ these elements $x$ and $y_\pm$ are unique. Again due to injectivity, $z=x+y_++y_-$ such that we can decompose 
     \begin{align}\label{eq:decomp_dom}
              \dom \left[\begin{smallmatrix}
         \calP \\ \calS
     \end{smallmatrix}\right] = \ker \Gamma \oplus \ker (\calS - i \calP) \oplus \ker (\calS + i \calP),
     \end{align}
     where the orthogonality remains to be shown. To this end, let $x \in \ker \Gamma$ and $y_{\pm} \in \ker (\calS - (\pm i \calP))$. As a consequence, $\calS y_\pm = \pm i \calP y_\pm$ holds and
     \begin{align*} 
         \langle x,y_\pm \rangle_{\dom \left[\begin{smallmatrix}
         \calP \\ \calS
     \end{smallmatrix}\right]} 
     = \pm i \mform{ \left[ \begin{smallmatrix}
        \calP \\ \calS 
    \end{smallmatrix}\right]x}{ \left[ \begin{smallmatrix}
        \calP  \\ \calS 
    \end{smallmatrix}\right]  y_\pm}{\calX}
      = \pm i \mform{ \left[ \begin{smallmatrix}
        \Gamma_0  \\ \Gamma_1
    \end{smallmatrix}\right]x}{ \left[ \begin{smallmatrix}
         \Gamma_0  \\ \Gamma_1
    \end{smallmatrix}\right] y }{\calG}  = 0,
     \end{align*}
     where we applied Lemma \ref{lem:defect_mform} and the identity \eqref{eq:imaginary_green}. This shows $\ker \Gamma \perp
     \ker( \calS - (\pm i \calP)) $. Furthermore, consider
     \begin{align*}
          \langle y_+,y_- \rangle_{\dom \left[\begin{smallmatrix}
         \calP \\ \calS
     \end{smallmatrix}\right]} 
     & = \langle \calP y_+, \calP y_- \rangle_\calX + \langle \calS y_+,  \calS y_- \rangle_\calX   \\
     & = \langle \calP y_+, \calP y_- \rangle_\calX + \langle i \calP y_+,  -i\calP y_- \rangle_\calX \\
     & = \langle \calP y_+, \calP y_- \rangle_\calX + i^2 \langle \calP y_+,  \calP y_- \rangle_\calX \\ & =0
     \end{align*}
     and thus the claimed orthogonality.
\end{proof}
The subsequent lemma shows the bijectivity of a restriction of the boundary map. The proof is based on an adaptation of \cite[Lem. 2.4.7]{skrepek2021linear} which is a similar result for skew-symmetric operators.
\begin{lem}\label{lem:bdd_map_invertibility}
Let Assumption~\ref{ass:ass1} be satisfied, let $H$ be symmetric and let $(\calG, \Gamma_0, \Gamma_1)$ be a boundary triplet for $H^*=\ran \left[\begin{smallmatrix}
        \calP \\ \calS
    \end{smallmatrix}\right]$. Then, the restriction $\Gamma'=\left. \Gamma \right|_{ \ker (\calS - i \calP) \oplus \ker (\calS + i \calP)}$ is bijective from $ \ker (\calS - i \calP) \oplus \ker (\calS + i \calP)$ to $\calG \times \calG$ and boundedly invertible.
\end{lem}
\begin{proof}
     By Lemma \ref{lem:orthogonal_decomposition} we decompose 
     \begin{align*}
         \dom \left[\begin{smallmatrix}
         \calP \\ \calS
     \end{smallmatrix}\right] = \ker \Gamma \oplus \ker (\calS - i \calP) \oplus \ker (\calS + i \calP).
     \end{align*}
     This shows immediately that $\ker \Gamma'=\ker \Gamma \cap (\ker (\calS - i \calP) \oplus \ker (\calS + i \calP))= {0}$. Further,
    \begin{align*}
        \ran \Gamma &= \left\lbrace \Gamma z \,\middle|\, z\in \dom  \left[\begin{smallmatrix}
         \calP \\ \calS
     \end{smallmatrix}\right]\right \rbrace \\ &=  \left\lbrace \Gamma x + \Gamma y \,  \middle| \, x \in \ker \Gamma, \, y \in \ker (\calS - i \calP) \oplus \ker (\calS + i \calP) \right\rbrace \\
     &= \left\lbrace \Gamma y \,  \middle|  \, y \in \ker (\calS - i \calP) \oplus \ker (\calS + i \calP)\right\rbrace = \ran \Gamma'
    \end{align*}
    such that $\Gamma'$ inherits the surjectivity of $\Gamma$. Thus, $\Gamma'$ maps $ \ker (\calS - i \calP) \oplus \ker (\calS + i \calP)$ isomorphic to $\calG \times \calG$. \\
    In the sequel, we will show that $\Gamma'^{-1}$ is bounded. To this end, let $\left( \left[ \begin{smallmatrix}
        u_n \\ v_n
    \end{smallmatrix}\right]\right)_{n \in \N}$ be a sequence in $\ran \left. \Gamma \right|_{\ker (\calS - i \calP)} \subset \calG \times \calG$ such that $\left( \left[ \begin{smallmatrix}
        u_n \\ v_n
    \end{smallmatrix}\right]\right)_{n \in \N} \stackrel{n \to \infty}{\to}  \left[ \begin{smallmatrix}
        u \\ v
    \end{smallmatrix}\right] \in \calG \times \calG$. We define
    \begin{align*}
        x_n = \Gamma'^{-1} \left[ \begin{smallmatrix}
        u_n \\ v_n
    \end{smallmatrix}\right], \qquad  x = \Gamma'^{-1} \left[ \begin{smallmatrix}
        u \\ v
    \end{smallmatrix}\right] \in \ker (\calS - i \calP) \oplus \ker (\calS + i \calP)
    \end{align*}
    which defines a sequence $(x_n)_{n \in \N} \subset \ker (\calS - i \calP)$. We decompose $x=x_++x_-$ uniquely and orthogonally with $x_\pm \in \ker (\calS -(\pm  i \calP))$. This yields together with Lemma \ref{lem:defect_mform}
    \begin{align*}
        & ||x_n-x||^2_{\dom \left[\begin{smallmatrix}
        \calP \\ \calS
    \end{smallmatrix}\right]} \\ &= \langle x_n-x, x_n - x_+ \rangle_{\dom \left[\begin{smallmatrix}
        \calP \\ \calS
    \end{smallmatrix}\right]} - \langle x_n-x, x_- \rangle_{\dom \left[\begin{smallmatrix}
        \calP \\ \calS
    \end{smallmatrix}\right]} \\
    & = i \mform{ \left[ \begin{smallmatrix}
        \calP \\ \calS 
    \end{smallmatrix}\right](x_n -x)}{ \left[ \begin{smallmatrix}
        \calP  \\ \calS 
    \end{smallmatrix}\right]  (x_n - x_+)}{\calX} + i \mform{ \left[ \begin{smallmatrix}
        \calP \\ \calS 
    \end{smallmatrix}\right](x_n -x)}{ \left[ \begin{smallmatrix}
        \calP  \\ \calS 
    \end{smallmatrix}\right] x_-}{\calX}  \\
    & = i \mform{ \left[ \begin{smallmatrix}
        \Gamma_0 \\ \Gamma_1 
    \end{smallmatrix}\right](x_n -x)}{ \left[ \begin{smallmatrix}
        \Gamma_0 \\ \Gamma_1
    \end{smallmatrix}\right]  (x_n - x_+)}{\calG} + i \mform{ \left[ \begin{smallmatrix}
        \Gamma_0 \\ \Gamma_1
    \end{smallmatrix}\right](x_n -x)}{ \left[ \begin{smallmatrix}
        \Gamma_0 \\ \Gamma_1
    \end{smallmatrix}\right] x_-}{\calG}  \\
    & = i \mform{ \left[ \begin{smallmatrix}
        u_n \\ v_n
    \end{smallmatrix}\right]- \left[ \begin{smallmatrix}
        u \\ v
    \end{smallmatrix}\right]}{\left[ \begin{smallmatrix}
        u_n \\ v_n
    \end{smallmatrix}\right] - \left[ \begin{smallmatrix}
        \Gamma_0 \\ \Gamma_1
    \end{smallmatrix}\right] x_+}{\calG} + i \mform{  \left[ \begin{smallmatrix}
        u_n \\ v_n
    \end{smallmatrix}\right]- \left[ \begin{smallmatrix}
        u \\ v
    \end{smallmatrix}\right]}{ \left[ \begin{smallmatrix}
        \Gamma_0 \\ \Gamma_1
    \end{smallmatrix}\right] x_-}{\calX} \stackrel{n \to \infty}{\to} 0
    \end{align*}
since $ \left[ \begin{smallmatrix}
        u_n \\ v_n
    \end{smallmatrix}\right]$ is bounded and converges to $\left[ \begin{smallmatrix}
        u \\ v
    \end{smallmatrix}\right]$. As a consequence, $x_n \stackrel{n \to \infty}{\to} x$ with respect to 
    $\|\cdot\|_{\dom \left[\begin{smallmatrix}
        \calP \\ \calS
    \end{smallmatrix}\right]}$ and $x \in \ker (\calS - i \calP)$ by the closedness of $\ker (\calS - i \calP)$ and we conclude that $\left. \Gamma ^{-1}\right|_{\ker (\calS - i \calP)} $ is bounded. Furthermore, this also implies the closedness of $\ran \left. \Gamma \right|_{\ker (\calS - i \calP)}$.
    By analogue argumentation, we see the boundedness of $\left. \Gamma ^{-1}\right|_{\ker (\calS + i \calP)}$ and hence $\Gamma'^{-1}$ is bounded.
\end{proof}
In the following theorem we show, that boundary triplets can be used to parametrize the set of all intermediate extensions between $\left. \left[\begin{smallmatrix}
        \calP \\ \calS
    \end{smallmatrix}\right]\right|_{\ker \Gamma} $ and $ \left[\begin{smallmatrix}
    \calP \\ \calS
    \end{smallmatrix}\right]$ by forcing the values of the boundary map to be contained in a linear relation in the boundary space.
\begin{thm}\label{thm:intermediate_extensions}
Let Assumption~\ref{ass:ass1} hold, $H \subset H^*$ and let  $(\calG, \Gamma_0, \Gamma_1)$ be a boundary triplet for $H^*=\ran \left[\begin{smallmatrix}
        \calP \\ \calS
    \end{smallmatrix}\right]$ and for $\Theta \subset \calG \times \calG$ define $A_\Theta: \calX \supset \dom A_\Theta \to \calX \times \calX$ by 
    \begin{equation}\label{eq:operator_theta}
            A_\Theta = \left[\begin{smallmatrix}
        \calP \\ \calS
    \end{smallmatrix}\right] , \qquad \dom A_\Theta = \left\lbrace  x \in \dom \left[\begin{smallmatrix}
        \calP \\ \calS
    \end{smallmatrix}\right] \, \middle| \, \Gamma x \in \Theta \right\rbrace. 
        \end{equation}
        Then the following holds:
    \begin{enumerate}
    \item[(i)] For any linear operator $A :\calX \supset \dom A \to \calX \times \calX$ with $ \ran \left. \left[\begin{smallmatrix}
        \calP \\ \calS
    \end{smallmatrix}\right]\right|_{\ker \Gamma} \subset \ran A \subset  \ran \left[\begin{smallmatrix}
    \calP \\ \calS
    \end{smallmatrix}\right]$ there is a unique $\Theta \subset \calG \times \calG$ such that $A=A_\Theta$. 
    \item[(ii)] $\overline{\ran A_\Theta}=\ran A_{\overline{\Theta}}$ for all $\Theta \subset \calG \times \calG$;
    \item[(iii)] $(\ran A_\Theta)^*=\ran A_{\Theta^*}$ for all $\Theta \subset \calG \times \calG$;
    \item[(iv)] $\ran A_\Theta \subset \ran A_{\Theta'}$ if and only if $\Theta \subset \Theta'$ for $\Theta, \Theta' \subset \calG \times \calG$.
    \end{enumerate}
\end{thm}
\begin{proof}
\begin{itemize}
\item[(i)]  Following Lemma \ref{lem:orthogonal_decomposition}, we orthogonally decompose
    \begin{equation}
    \dom \left[\begin{smallmatrix}
         \calP \\ \calS
     \end{smallmatrix}\right] = \ker \Gamma \oplus \ker (\calS - i \calP) \oplus \ker (\calS + i \calP).
\end{equation}
    By Lemma \ref{lem:bdd_map_invertibility}, we equivalently rewrite this orthogonal decomposition via
    \begin{equation}
    \dom \left[\begin{smallmatrix}
         \calP \\ \calS
     \end{smallmatrix}\right] = \ker \Gamma \oplus \ran \left. \Gamma \right|_{ \ker (\calS - i \calP) \oplus \ker (\calS + i \calP)}^{-1}.
\end{equation}
Thus, any $A: \calX \supset \dom A \to \calX \times \calX$ with $\overline{H}=\ran \left. \left[\begin{smallmatrix}
        \calP \\ \calS
    \end{smallmatrix}\right] \right|_{\ker \Gamma}  \subset \ran A \subset \ran \left[\begin{smallmatrix}
        \calP \\ \calS
    \end{smallmatrix}\right]  =H^*$ satisfies
    \begin{align*}
        \dom A = \Gamma^{-1}\{0\} \oplus \left. \Gamma \right|_{ \ker (\calS - i \calP) \oplus \ker (\calS + i \calP)}^{-1} \Theta = \dom A_\Theta
    \end{align*}
     for a unique $\Theta \subset \calG \times \calG$.  
\item[(ii)]  Clearly, by Lemma \ref{lem:bdd_map_invertibility}, we have 
    \begin{align*}
        \overline{\left. \Gamma \right|_{ \ker (\calS - i \calP) \oplus \ker (\calS + i \calP)}^{-1} \Theta }=\left. \Gamma \right|_{ \ker (\calS - i \calP) \oplus \ker (\calS + i \calP)}^{-1} \overline{\Theta}.
    \end{align*}
    We define the embedding 
    \begin{equation*}
        \iota: \dom \left[\begin{smallmatrix}
        \calP \\ \calS
    \end{smallmatrix}\right] \to \ran \left[\begin{smallmatrix}
        \calP \\ \calS
    \end{smallmatrix}\right], \qquad x \mapsto \left[\begin{smallmatrix}
        \calP \\ \calS
    \end{smallmatrix}\right] x
    \end{equation*}
    which is continuous by Lemma \ref{lem:injective_and_closed_range}. As a consequence, 
    \begin{align*}
        \overline{\ran A_\Theta} & =\overline{\overline{H} \ \oplus \ \iota \left. \Gamma \right|_{ \ker (\calS - i \calP) \oplus \ker (\calS + i \calP)}^{-1} \Theta}= \overline{H} \ \oplus \ \iota \left. \Gamma \right|_{ \ker (\calS - i \calP) \oplus \ker (\calS + i \calP)}^{-1}\overline{\Theta} \\
        &= \ran A_{\overline{\Theta}}.
    \end{align*}
    \item[(iii)] Before we start, recall that for any intermediate extension $B$ of $\overline{H}$, the adjoint $B^*$ is also an intermediate extension of $\overline{H}$. At first, we show $(\ran A_\Theta)^* \subset \ran A_{\Theta^*}$. Therefore, let $\left[\begin{smallmatrix}
        g_1 \\ g_2
    \end{smallmatrix}\right] \in (\ran A_\Theta)^*=(\ran A_\Theta)^{\mperp}\subset H^*$. In particular, this means that there is a $y \in \dom \left[\begin{smallmatrix}
        \calP \\ \calS
    \end{smallmatrix}\right]$ with $\left[\begin{smallmatrix}
        g_1 \\ g_2
    \end{smallmatrix}\right]=\left[\begin{smallmatrix}
        \calP \\ \calS
    \end{smallmatrix}\right]y$. Moreover, let $(\varphi_1, \varphi_2) \in \Theta$ and choose $x \in \dom A_\Theta$ with $\Gamma x = \left[\begin{smallmatrix}
            \varphi_1 \\ \varphi_2
        \end{smallmatrix}\right]$, which is always possible by surjectivity of $\Gamma$. Then,
    \begin{align*}
         \left\langle \varphi_2, \Gamma_0 y \right\rangle_\calG -\left\langle \varphi_1, \Gamma_1 y \right\rangle_\calG & = \left\langle \Gamma_1 x, \Gamma_0 y \right\rangle_\calG -\left\langle \Gamma_0 x, \Gamma_1 y \right\rangle_\calG \\
         &= \left\langle \calS x , \calP y \right\rangle_\calX - \left\langle \calP x , \calS y \right\rangle_\calX=0,
    \end{align*}
    and thus $\Gamma y \in \Theta^{\mperp}=\Theta^*$ which shows $(g,g') \in \ran A_{\Theta^*}$. \\
    Conversely, let $x \in \dom A_{\Theta^*}$, i.e., $\Gamma x \in \Theta^*=\Theta^{\mperp}$. Further, let $(\varphi_1, \varphi_2) \in \Theta$ and choose $y \in \dom A_\Theta$ such that $\Gamma y = \left[\begin{smallmatrix}
        \varphi_1 \\ \varphi_2
    \end{smallmatrix} \right] $. Thus,
    \begin{align*}
          0 &=
        \left\langle \varphi_2, \Gamma_0 x \right\rangle_\calG -\left\langle \varphi_1, \Gamma_1 x \right\rangle_\calG  =   \left\langle \Gamma_1 y, \Gamma_0 x \right\rangle_\calG -\left\langle \Gamma_0 y, \Gamma_1 x \right\rangle_\calG \\
        &= \left\langle \calS y , \calP x \right\rangle_\calX - \left\langle \calP y , \calS x \right\rangle_\calX,
    \end{align*}
    i.e., $ \left[\begin{smallmatrix}
        \calP \\ \calS
    \end{smallmatrix}\right]x \in (\ran A_{\Theta})^{\mperp}=(\ran A_\Theta)^*$. 
\item[(iv)] The assertion follows directly from the definition of $A_{\Theta}$.
    \end{itemize}
    \end{proof}
As a consequence of Theorem \ref{thm:intermediate_extensions} (iii), (iv), combined with \cite[Cor. 2.1.4]{BehrHdS20} we immediately see that Theorem \ref{thm:self-adj_extension} is true. \par
In the following remark, we discuss the relation of the boundary triplet approach for $H^*$ suggested in Definition~\ref{def:range_triplet} with the concept of ordinary boundary triplets, cf. \cite[Def. 2.1.1]{BehrHdS20}.
    \begin{rem}\label{rem:relation_between_triplets}
    Under Assumption~\ref{ass:ass1}, the boundary triplet for $H^*=\ran \left[\begin{smallmatrix}
        \calP \\ \calS
    \end{smallmatrix}\right]$ can be related to an ordinary boundary triplet for $H^*$. Recall, that $\left[\begin{smallmatrix}
        \calP \\ \calS
    \end{smallmatrix}\right]: \calX \supset \dom \left[\begin{smallmatrix}
        \calP \\ \calS
    \end{smallmatrix}\right] \to \calX \times \calX$ is injective by Lemma \ref{lem:injective_and_closed_range}, and hence,  there exists a bounded inverse $\left[\begin{smallmatrix}
        \calP \\ \calS
    \end{smallmatrix}\right]^{-1}: \ran \left[\begin{smallmatrix}
        \calP \\ \calS
    \end{smallmatrix}\right] \to \dom \left[\begin{smallmatrix}
        \calP \\ \calS
    \end{smallmatrix}\right]$.
    Since $\left[\begin{smallmatrix}
        \calP \\ \calS
    \end{smallmatrix}\right]^{-1}$ is surjective, $\Gamma \left[\begin{smallmatrix}
        \calP \\ \calS
    \end{smallmatrix}\right]^{-1}$ is surjective if and only if $\Gamma$ is surjective. Moreover, to show that the abstract Green identity \eqref{eq:Green} holds if and only if \eqref{eq:Green_for_relations} holds, we define $\left[\begin{smallmatrix}
        f \\ f'
    \end{smallmatrix}\right]\coloneqq \left[\begin{smallmatrix}
        \calP \\ \calS
    \end{smallmatrix}\right] x, \left[\begin{smallmatrix}
        g \\ g'
    \end{smallmatrix}\right]\coloneqq \left[\begin{smallmatrix}
        \calP \\ \calS
    \end{smallmatrix}\right] y$ and obtain 
    \begin{align*}
         \left\langle \calS y , \calP x \right\rangle_\calX - \left\langle \calP y , \calS x \right\rangle_\calX = \left\langle \Gamma_1 y, \Gamma_0 x \right\rangle_\calG -\left\langle \Gamma_0 y, \Gamma_1 x \right\rangle_\calG
    \end{align*}
    if and only if
    \begin{align*}
        \left\langle f' , g \right\rangle_\calX - \left\langle f , g' \right\rangle_\calX 
         &= \left\langle \Gamma_1  \left[\begin{smallmatrix}
        \calP \\ \calS
    \end{smallmatrix}\right]^{-1} \left[ \begin{smallmatrix}
            f \\ f'
        \end{smallmatrix} \right], \Gamma_0  \left[\begin{smallmatrix}
        \calP \\ \calS
    \end{smallmatrix}\right]^{-1} \left[ \begin{smallmatrix}
            g \\ g'
        \end{smallmatrix} \right] \right\rangle_\calG \\
        & \hphantom{=} -\left\langle \Gamma_0  \left[\begin{smallmatrix}
        \calP \\ \calS
    \end{smallmatrix}\right]^{-1} \left[ \begin{smallmatrix}
            f \\ f'
        \end{smallmatrix} \right], \Gamma_1  \left[\begin{smallmatrix}
        \calP \\ \calS
    \end{smallmatrix}\right]^{-1}\left[ \begin{smallmatrix}
            g \\ g'
        \end{smallmatrix} \right]\right\rangle_\calG .
    \end{align*}
    Consequently, $(\calG, \Gamma_0, \Gamma_1)$ is a boundary triplet for $H^*=\ran \left[\begin{smallmatrix}
        \calP \\ \calS
    \end{smallmatrix}\right]$ if and only if the triplet $(\calG, \Gamma_0 \left[\begin{smallmatrix}
        \calP \\ \calS
    \end{smallmatrix}\right]^{-1}, \Gamma_1 \left[\begin{smallmatrix}
        \calP \\ \calS
    \end{smallmatrix}\right]^{-1})$ where $\Gamma \left[\begin{smallmatrix}
        \calP \\ \calS
    \end{smallmatrix}\right]^{-1} =  \left[ \begin{smallmatrix}
        \Gamma_0 \\ \Gamma_1
    \end{smallmatrix} \right]\left[\begin{smallmatrix}
        \calP \\ \calS
    \end{smallmatrix}\right]^{-1}: H^* \to \calG \times \calG$ is an ordinary boundary triplet for $H^*$. However, for most applications it is not obvious how to compute $\left[\begin{smallmatrix}
        \calP \\ \calS
    \end{smallmatrix}\right]^{-1}$.
    \end{rem}
    The following theorem provides a necessary and sufficient condition for the existence of a boundary triplet for $H^*$. 
    \begin{thm}\label{thm:existence}
        Let Assumption~\ref{ass:ass1} hold and assume that $H$ is symmetric. There exists a boundary triplet for $H^*=\ran \left[\begin{smallmatrix}
        \calP \\ \calS
    \end{smallmatrix}\right]$ if and only if 
        \begin{align}
            \dim \ker  (\calS - i \calP) = \dim \ker (\calS + i \calP).
        \end{align}
    \end{thm}
    \begin{proof}
        As shown in \cite[Chap. 2.4]{BehrHdS20}, an ordinary boundary triplet $(\calG', \Gamma_0', \Gamma_1')$ for $H^*=\ran \left[\begin{smallmatrix}
        \calP \\ \calS
    \end{smallmatrix}\right]$ exists if and only if
    \begin{align}
        \dim \mathfrak{N}_\mu (H^*) =\dim \mathfrak{N}_{\overline{\mu}} (H^*)
    \end{align}
    for some and hence for all $\mu \in \C \setminus \R$. Moreover, by \cite[Thm. 2.5.1]{BehrHdS20} there is a unique bounded $\calW : \calG \times \calG \to \calG' \times \calG'$ such that $\Gamma'=\calW \Gamma  \left[\begin{smallmatrix}
        \calP \\ \calS
    \end{smallmatrix}\right]^{-1}$ and that $(\calG, \Gamma_0 \left[\begin{smallmatrix}
        \calP \\ \calS
    \end{smallmatrix}\right]^{-1}, \Gamma_1 \left[\begin{smallmatrix}
        \calP \\ \calS
    \end{smallmatrix}\right]^{-1})$ is an ordinary boundary triplet for $H^*=\ran \left[\begin{smallmatrix}
        \calP \\ \calS
    \end{smallmatrix}\right]$. In particular, this is equivalent to the existence of a boundary triplet for $H^*=\ran \left[\begin{smallmatrix}
        \calP \\ \calS
    \end{smallmatrix}\right]$ by Remark~\ref{rem:relation_between_triplets}. Specifically, this implies that
    \begin{align}\label{eq:dim_defect_spaces}
        \dim \calP (\ker (\calS-i\calP)) = \dim \calP (\ker (\calS+i\calP)) 
    \end{align}
    is a necessary and sufficient condition, where we used 
    \begin{align*}
        \mathfrak{N}_\mu (H^*) = \ker (H^* -\mu)= \calP (\ker (\calS-\mu\calP))
    \end{align*}
    shown in Lemma \ref{lem:defect_representation}. Moreover, as it was shown in \cite{gernandt2020invariance} for operator pencils with bounded coefficients, the operator $\left. \calP\right|_{\ker (\calS-(\pm i\calP))}$ maps $\ker (\calS-(\pm i\calP))$ bijectively to $\calP(\ker (\calS-(\pm i \calP)))$. Indeed, this can be seen by
    \begin{align}
        \ker \left. \calP\right|_{\ker (\calS-(\pm i\calP)) }= \left\lbrace x \in \dom  \left[\begin{smallmatrix}
        \calP \\ \calS
    \end{smallmatrix}\right] \, \middle| \mp i \calS x=\calP x =0\right\rbrace=\{0\}
    \end{align}
    since $\ker \calP \cap \ker \calS = \{0\}$ by Lemma \ref{lem:injective_and_closed_range}. Clearly, any operator maps surjectively onto its range. Consequently, \eqref{eq:dim_defect_spaces} can be rewritten by
    \begin{align*}
        \dim \ker (\calS-i\calP) = \dim \ker (\calS+i\calP)
    \end{align*}
    which shows the claim. 
    \end{proof}
    \par
    The following lemma provides a way to construct a boundary triplet by verifying the abstract Green identity \eqref{eq:Green} on a dense subspaces. 
    \begin{lem}\label{lem:extension_by_density}
        Let $D$ be a dense subspace of $\dom \left[\begin{smallmatrix}
        \calP \\ \calS
    \end{smallmatrix}\right]$ with respect to the topology induced by \eqref{eq:dom_inner_product}. Assume that $H \subset H^*$ and that there exist a Hilbert space $\calG$ and linear, surjective and bounded maps $\Gamma_0, \Gamma_1: D \to \calG$ such that 
    \begin{equation}
        \left\langle \calS x , \calP y \right\rangle_\calX - \left\langle \calP x , \calS y \right\rangle_\calX = \left\langle \Gamma_1 x, \Gamma_0 y \right\rangle_\calG -\left\langle \Gamma_0 x, \Gamma_1 y\right\rangle_\calG\label{eq:green_on_dense_subspace} 
    \end{equation}
    holds for all $x,y \in D$. Then, $(\calG, \widetilde{\Gamma_0}, \widetilde{\Gamma_1})$ is a boundary triplet for $H^*=\ran \left[\begin{smallmatrix}
        \calP \\ \calS
    \end{smallmatrix}\right]$, where $\widetilde{\Gamma_0}, \widetilde{\Gamma_1}: \dom \left[\begin{smallmatrix}
        \calP \\ \calS
    \end{smallmatrix}\right] \to \calG$ denote the unique extensions of $\Gamma_0, \Gamma_1$ to $\dom \left[\begin{smallmatrix}
        \calP \\ \calS
    \end{smallmatrix}\right]$, respectively.
    \end{lem}
    \begin{proof}
        Following \cite[Thm. 1.9.1]{megginson2012introduction}, we conclude that $\Gamma_0, \Gamma_1: D \to \calG$ extend uniquely to linear and bounded operators $\widetilde{\Gamma_0}, \widetilde{\Gamma_1}: \dom \left[\begin{smallmatrix}
        \calP \\ \calS
    \end{smallmatrix}\right] \to \calG$. In particular, this means that $\left[\begin{smallmatrix}
        \widetilde{\Gamma_0} \\ \widetilde{\Gamma_1}
    \end{smallmatrix}\right]x = \lim\limits_{n \to \infty}  \left[\begin{smallmatrix}
        \Gamma_0 \\ \Gamma_1
    \end{smallmatrix}\right]x_n$ for a sequence $(x_n)_{n \in \N} \subset D$ with $x_n \to x$ as $n \to \infty$. Obviously, $\widetilde{\Gamma_0}, \widetilde{\Gamma_1}: \dom \left[\begin{smallmatrix}
        \calP \\ \calS
    \end{smallmatrix}\right] \to \calG$ are surjective because the restrictions to $D$ is already surjective, respectively. To show that the abstract Green identity extends to all functions in $\dom \left[\begin{smallmatrix}
        \calP \\ \calS
    \end{smallmatrix}\right]$ we let $(x_n)_{n \in \N}, (y_n)_{n \in \N}  \subset D$ with $x_n \to x, y_n \to y$ as $n \to \infty$. This yields
    \begin{equation*}
        \begin{array}{rcl}
             \left\langle \calS x , \calP y \right\rangle_\calX - \left\langle \calP x , \calS y \right\rangle_\calX 
             &=& \lim\limits_{n \to \infty} \left\langle \calS x_n , \calP y_n \right\rangle_\calX - \left\langle \calP x_n , \calS y_n \right\rangle_\calX \\
             &=& \lim\limits_{n \to \infty}  \left\langle \Gamma_1 x_n, \Gamma_0 y_n \right\rangle_\calG -\left\langle \Gamma_0 x_n, \Gamma_1 y_n\right\rangle_\calG \\
             &=&  \left\langle \widetilde{\Gamma_1} x, \widetilde{\Gamma_0} y \right\rangle_\calG -\left\langle \widetilde{\Gamma_0} x, \widetilde{\Gamma_1}y\right\rangle_\calG 
        \end{array}
    \end{equation*}
    for all $x,y \in \dom \left[\begin{smallmatrix}
        \calP \\ \calS
    \end{smallmatrix}\right]$ by using the closedness of the operators $\calS, \calP$ and the Green identity \eqref{eq:green_on_dense_subspace} on the dense subspace $D$. Eventually, this shows that $(\calG, \widetilde{\Gamma_0}, \widetilde{\Gamma_1})$ is a boundary triplet for $\ran \left[\begin{smallmatrix}
        \calP \\ \calS
    \end{smallmatrix}\right]$.
    \end{proof}
\subsection{An Example: Dzektser equation}\label{subsec:dzektser}
In this section, we illustrate our approach by means of the Dzektser equation \cite{Dexter,bendimerad2024stokes}. More precisely, we will construct a boundary triplet in the sense of Definition~\ref{def:operator_bdd_triplet}, and, consequently, provide a full characterization of self-adjointness of the corresponding relation by applying our main result Theorem~\ref{thm:self-adj_extension}. 
In our example, we consider a modification of the original Dzektser equation from \cite{jacob2022solvability} that was studied in the context of solution theory of dissipative partial differential algebraic equations.
    \begin{equation}\label{eq:Dzektser}
        \frac{\mathrm{d}}{\mathrm{d} t} \left(1+ \dif{2}\right)x(t,\xi)=\left(\dif{2} + 2\dif{4}\right)x(t,\xi), \quad (t,\xi)\in[0,\infty)\times[0,\pi]. 
    \end{equation}
    Our goal is to characterize all boundary conditions such that \eqref{eq:Dzektser} admits a self-adjoint representation in the spirit of the previous Section \ref{subsec:bdd_triplets}. To this end, we consider the Hilbert space $\calX = L^2(0,\pi)$ and the operators $\calP_{D,0}, \calS_{D,0}$ given by 
        \begin{equation*}
        \begin{array}{rclrl}
            \calP_{D,0} x&=&x+\dif{2}x, &\qquad \dom \calP_{D,0} = &H_0^2(0,\pi),\\
            \calS_{D,0} x&=&\dif{2}x + 2\dif{4}x, &\qquad \dom \calS_{D,0} =& H_0^4(0,\pi).
        \end{array}
        \end{equation*}
        The adjoint operators of $\calP_{D,0}$ and $\calS_{D,0}$ are given by
        \begin{equation*}
        \begin{array}{rclrl}
            \calP_D x&=&x +\dif{2}x, &\qquad \dom \calP_D = &H^2(0,\pi),\\
            \calS_D x&=&\dif{2}x + 2\dif{4}x, &\qquad \dom \calS_D =& H^4(0,\pi),
        \end{array}
        \end{equation*}
        respectively. Accordingly, we define the combined operator by $\left[\begin{smallmatrix}
        \calP_D \\ \calS_D
    \end{smallmatrix}\right]:  \calX \supset H^4(0,\pi) \to \calX \times \calX$. 
    In the following, we verify that the associated sesquilinear form $\mathfrak{t}_D: H^4(0,\pi) \times H^4(0,\pi) \to \C$ via
    \begin{align}
    \label{def:t_D}
            \mathfrak{t}_D(x,y)=\langle \calP_D x, \calP_D y\rangle_{\calX}+\langle \calS_D x, \calS_D y\rangle_{\calX}
        \end{align}
    fulfills the coercivity condition
    \begin{align}
    \label{eq:t_d_coercive}
        \mathfrak{t}_D(x,x)\geq \|x\|_{\calX}^2\quad \text{for all $x\in H^4(0,\pi)$}.
        \end{align}
To prove this, one has to observe that it is sufficient to verify \eqref{eq:t_d_coercive} for all elements of an orthonormal basis. In particular, we choose the orthonormal basis  $(\varphi_k)_{k=1}^\infty$ that is given by the eigenfunctions $\varphi_k(\cdot)=\sqrt{\tfrac{2}{\pi}}\sin(k\cdot)$ of the Dirichlet Laplacian for the eigenvalue $\lambda_k=-k^2$. This results in the following estimate for $k \geq 1$ 
\begin{align*}
\mathfrak{t}_D(\varphi_k,\varphi_k)&\geq \langle \calP_D \varphi_k, \calP_D \varphi_k\rangle_{\calX}+\langle \calS_D \varphi_k, \calS_D \varphi_k\rangle_{\calX}\\&\geq \langle (1-k^2)\varphi_k, (1-k^2) \varphi_k\rangle_{\calX}+\langle (2k^4-k^2) \varphi_k, (2k^4-k^2) \varphi_k\rangle_{\calX}\\
&=(1-k^2)^2+(2k^4-k^2)^2\geq 1 = \|\varphi_k\|^2,
\end{align*}
i.e., \eqref{eq:t_d_coercive}. As a consequence, Assumption~\ref{ass:ass1} is fulfilled for operators associated with the Dzektser equation \eqref{eq:Dzektser}. 
Hence, $(H^4(0,\pi), \mathfrak{t}_D)$ defines a Hilbert space. Obviously, $H=\ker \left[\begin{smallmatrix}
    - \calS_{D,0} & \calP_{D,0} 
\end{smallmatrix}\right]$ is symmetric. We now construct a boundary triplet for $H^*=\ran \left[\begin{smallmatrix}
        \calP_D \\ \calS_D
    \end{smallmatrix}\right]$ in terms of boundary evaluations of functions in $H^4(0,\pi)$.
The following result is taken from \cite[Thm. 7.8.1]{aubin2000applied} and follows by straightforward interpolation.
\begin{prop}\label{prop:trace_surjective}
    Let $k \in \N$. Then the trace operator $\gamma: H^k(a,b; \C^n) \to \C^{2nk}$ given by
    \begin{equation}\label{eq:trace}
        \gamma x = \begin{bmatrix}
            x(b) ~ x'(b) & \ldots & x^{(k-1)}(b) ~ x(a) ~ x'(a) & \ldots & x^{(k-1)}(a)
        \end{bmatrix}^\top
    \end{equation}
 is linear, bounded and surjective and satisfies $\ker \gamma=H_0^k(a,b; \C^n)$.
\end{prop}
\noindent Define the boundary maps $\Gamma =  \begin{bmatrix}
        \Gamma_0 \\ \Gamma_1
    \end{bmatrix} : H^{4}(0, \pi) \to \C^{4}\times \C^{4}$
\begin{align}
    \label{eq:bt_dexter}
\Gamma x = \begin{bmatrix}
        \Gamma_0 \\ \Gamma_1
    \end{bmatrix} x = \tfrac{1}{\sqrt{2}} \begin{bmatrix}
        A & -A \\ -I & -I
    \end{bmatrix} \gamma x,\quad 
         A= \begin{bmatrix}
                \phantom{-}0&1&\phantom{-}0&2 \\ -1&0&-2&0 \\ \phantom{-}0&2&\phantom{-}0&2 \\ -2&0&-2&0
            \end{bmatrix} \in \R^{4 \times 4},
    \end{align}
    where $\gamma: H^{4}(0,\pi) \to \C^{8}$ denotes the trace operator \eqref{eq:trace}. 
\begin{thm}\label{thm:Dzektser_bdd_triplet}
    We have that $(\C^4, \Gamma_0, \Gamma_1)$ given by \eqref{eq:bt_dexter} defines a boundary triplet for $\ran \left[\begin{smallmatrix}
        \calP_D \\ \calS_D
    \end{smallmatrix}\right]$.
\end{thm}
\begin{proof}
    Observe that
    \begin{align*}
        & \left\langle \calS_D x , \calP_D y \right\rangle_\calX \\
        &= \langle x'', y \rangle_\calX + \langle x'',  y'' \rangle_\calX + \langle 2x^{(4)}, y \rangle_\calX + \langle 2x^{(4)}, y'' \rangle_\calX \\
        &=  \left[  x'y - xy' \right]_0^\pi + 2\left[  x^{(3)}y - x''y' + x'y'' - xy^{(3)} + x^{(3)}y'' - x''y^{(3)}\right]_0^\pi \\
        &\hspace{2.59ex} + \underbrace{\langle x, y'' \rangle_\calX + \langle x'', y'' \rangle_\calX + \langle x, 2y^{(4)} \rangle_\calX + \langle x'', 2y^{(4)} \rangle_\calX}_{=\left\langle \calP_D x , \calS_D y \right\rangle_\calX }.
    \end{align*}
    for all $x,y \in H^4(0,\pi)$. This is equivalent to
    \begin{align*}
        \left\langle \calS_D x , \calP_D y \right\rangle_\calX - \left\langle \calP_D x , \calS_D y \right\rangle_\calX & = (\gamma x)^H  \left[ \begin{smallmatrix}
        A^H &0 \\ 0 & -A^H
    \end{smallmatrix} \right] \gamma y = \left\langle \Gamma x, \left[ \begin{smallmatrix}
       0 & -I \\ I &0
    \end{smallmatrix} \right]  \Gamma y \right\rangle_{\C^8},
    \end{align*}
    where we used the fact that 
    \begin{align*}
        \frac{1}{\sqrt{2}} \begin{bmatrix}
        A & -A \\ -I & -I
    \end{bmatrix}^H  \begin{bmatrix}
       0 & -I \\ I &0
    \end{bmatrix}  \begin{bmatrix}
        A & -A \\ -I & -I
    \end{bmatrix}\frac{1}{\sqrt{2}} =  \begin{bmatrix}
        A^H &0 \\  0& -A^H
    \end{bmatrix}.
    \end{align*}
    Clearly, $\C^4$ is a Hilbert space and $\Gamma$ is surjective by Proposition~\ref{prop:trace_surjective} and as the composition of surjective maps is surjective.
\end{proof}
Observe that $\ker \Gamma = \ker \gamma = H_0^4(0,\pi)$, and consequently, Lemma \ref{lem:orthogonal_decomposition} leads to the orthogonal decomposition with respect to $\mathfrak{t}_D$
    \begin{align*}
        H^4(0,\pi) = H_0^4(0, \pi)  \oplus \ker (\tdif{4} - \tfrac{1+i}{2} \tdif{2} - \tfrac{i}{2}) \oplus \ker (\tdif{4} + \tfrac{1+i}{2} \tdif{2} + \tfrac{i}{2}).
    \end{align*}
As a consequence of by Theorem~\ref{thm:self-adj_extension}, the relation $\ran  A_\Theta$ with 
    \begin{align*}
         A_\Theta = \begin{bmatrix}
        \calP_D \\ \calS_D
    \end{bmatrix} , \qquad \dom A_\Theta = \left\lbrace  x \in H^4(0,\pi) \, \middle| \, \Gamma x \in \Theta \right\rbrace,
    \end{align*}
    is self-adjoint if and only if $\Theta \subset \C^4 \times \C^4$ is self-adjoint.

\subsection{Correspondence of boundary triplets and Lagrangian subspaces}\label{subsec:Lagrange}
Lagrangian subspaces are a central topic in port-Hamiltonian systems, see \cite{MascvdSc18, MascvdSc23, bendimerad2024stokes}. In this part, we analyze the connection between the concept of Lagrangian subspaces defined on the so-called \emph{bond space} and our proposed extension theory of linear relations. When dealing with differential operators, these two domains are closely related via the notion of a boundary triplet. We note that the correspondence of skew-symmetric and skew-adjoint relations to Dirac structures has been established neatly in \cite{kurula2010dirac}. 

In the following, let $\calF$ be a complex Hilbert space, the so-called \textit{flow space} and its Hilbert space dual $\calE=\calF^*$ which we call \textit{effort space}. Let $U:\calE \to \calF$ be a unitary operator. As usual the product space $\calF \times \calE$ is endowed with the natural inner product
\begin{equation*}
    \left\langle \begin{bmatrix}
        f_1 \\ e_1
    \end{bmatrix}, \begin{bmatrix}
        f_2 \\ e_2
    \end{bmatrix} \right\rangle_{\calF \times \calE} = \langle f_1, f_2 \rangle_\calF + \langle e_1, e_2 \rangle_\calE.
\end{equation*}
Furthermore, we define the \textit{bond space} $\calB=\calF \times \calE$ endowed with the indefinite inner product $\langle\langle  \cdot , \cdot \rangle\rangle_- : \calB \times \calB \to \C$ defined by
    \begin{equation*}
        \left\langle \left\langle  \left[\begin{smallmatrix}
            f_1 \\ e_1
        \end{smallmatrix} \right], \left[\begin{smallmatrix}
            f_2 \\ e_2
        \end{smallmatrix} \right] \right\rangle \right\rangle_- = \left\langle  \left[\begin{smallmatrix}
            f_1 \\ e_1
        \end{smallmatrix} \right], \left[ \begin{smallmatrix}
           0 & -U \\ U^* &
0        \end{smallmatrix} \right] \left[\begin{smallmatrix}
            f_2 \\ e_2
        \end{smallmatrix} \right] \right\rangle_{\calF \times \calE}= \langle e_1 , U^*f_2 \rangle_\calE -  \langle f_1 , Ue_2 \rangle_\calF.
    \end{equation*}
    For a linear subspace $\calL \subset \calB$, we define 
    \begin{equation*}
        \calL^{\mrperp} = \left\lbrace \left[\begin{smallmatrix}
            f_2 \\ e_2
        \end{smallmatrix} \right]\in \calB \,  \middle| \, \left\langle \left\langle  \left[\begin{smallmatrix}
            f_1 \\ e_1
        \end{smallmatrix} \right], \left[\begin{smallmatrix}
            f_2 \\ e_2
        \end{smallmatrix} \right] \right\rangle \right\rangle_- =0  \ \text{for all}  \left[\begin{smallmatrix}
            f_1 \\ e_1
        \end{smallmatrix} \right] \in \calL  \right\rbrace.
    \end{equation*}
    By definition, we have
    \begin{equation*}
        \calL^{\mrperp}=\left[ \begin{smallmatrix}
            0& -U \\ U^* &0
        \end{smallmatrix} \right] \calL^\perp.
    \end{equation*}
    \begin{defn}
        A linear subspace $\calL \subset \calB$ is a \textit{Lagrangian subspace} if $\calL=\calL^{\mrperp}$.
    \end{defn}
    \begin{cor}
        We identify the flow space $\calF$ with its dual $\calE$ by means of the Riesz isomorphism. Further, let $U=I:\calF \to \calF$. Then, for a subspace $\calL \subset \calB$ we have $\calL^{\mrperp}=\calL^*$. In this case, a linear relation $\calL$ is Lagrangian if and only if it is self-adjoint. 
    \end{cor}
    \begin{rem}
        In the port-Hamiltonian literature, cf. \cite{MascvdSc23}, a Lagrangian subspace defined on an infinite-dimensional state space is also called \textit{Stokes-Lagrange subspace}. For example, for $\tdif{2}: L^2(0,1) \supset H^2(0,1) \to L^2(0,1)$ we have 
        \begin{align*}
            \langle x'', y \rangle_{L^2(0,1)} - \langle x, y'' \rangle_{L^2(0,1)} &=x'(1)y(1)  -  x'(0)y(0)  - x(1)y'(1)  + x(0)y'(0)  \\
            &= \left\langle \hspace{-0.75mm} \begin{bmatrix}
                -x'(0) \\ \hphantom{-}x'(1)
            \end{bmatrix}\hspace{-0.5mm}, \hspace{-0.5mm}\begin{bmatrix}
                y(0) \\ y(1)
            \end{bmatrix} \hspace{-0.75mm}\right\rangle_{\C^2} - \left\langle\hspace{-0.75mm} \begin{bmatrix}
                x(0) \\ x(1)
            \end{bmatrix} \hspace{-0.5mm}, \hspace{-0.5mm}\begin{bmatrix}
                -y'(0) \\ \hphantom{-}y'(1)
            \end{bmatrix} \hspace{-0.75mm}\right\rangle_{\C^2}
        \end{align*}
        for all $x,y \in H^2(0,1)$ by the application of the integration by parts formula. By defining the linear and surjective \textit{boundary port} maps 
        \begin{align*}
            \chi_\partial,\varepsilon_\partial  & : H^2(0,1) \to \C^2 , \qquad \chi_\partial x = \left[ \begin{smallmatrix}
                x(0) \\ x(1)
            \end{smallmatrix} \right] ,\quad  \varepsilon_\partial x = \left[ \begin{smallmatrix}
                -x'(0) \\ \hphantom{-}x'(1)
            \end{smallmatrix} \right]
        \end{align*}
        this may be rewritten as
        \begin{align*}
            \left\langle  \left[\begin{smallmatrix}
             x \\ \chi_\partial x \\ x'' \\ \varepsilon_\partial x
        \end{smallmatrix} \right], \left[ \begin{smallmatrix}
           0 & -U \\ U^* &0
        \end{smallmatrix} \right] \left[\begin{smallmatrix}
            y \\ \chi_\partial y \\ y'' \\ \varepsilon_\partial y
        \end{smallmatrix} \right] \right\rangle_{L^2(0,1) \times \C^2}=0
        \end{align*}
        for all $x,y \in H^2(0,1)$ where $U=\left[\begin{smallmatrix}
            I_{L^2(0,1)} &0 \\ 0& -I_{\C^2}
        \end{smallmatrix} \right]$. 
        Moreover, we define the \textit{Stokes-Lagrange subspace}
        \begin{align*}
            \calL_{2} = \left\lbrace \left(\left[\begin{smallmatrix}
            f \\ \chi
        \end{smallmatrix} \right], \left[\begin{smallmatrix}
            e \\ \varepsilon
        \end{smallmatrix} \right]\right) \in (L^2(0,1) \times \C^2)^2 \, \middle| \, f \in H^2(0,1): e=f'', \, \left[\begin{smallmatrix}
            \chi \\ \varepsilon 
        \end{smallmatrix} \right] = \left[\begin{smallmatrix}
        \chi_\partial \\ \varepsilon_\partial 
    \end{smallmatrix}\right] f \right\rbrace
        \end{align*}
        which allows us to rewrite the partial differential equation
        \begin{equation*}
        \begin{array}{rcll}
            \frac{\partial}{\partial t} x(t,\xi) & =& \dif{2} x(t,\xi) , \qquad & \text{for} \ (t,\xi)\in(0, \infty) \times (0,1) \\
            x(t,\xi) &=& \dif{} x(t,\xi)=0 , \qquad & \text{for} \ (t,\xi)\in(0, \infty) \times \{0,1\}
        \end{array}
        \end{equation*}
        implicitly by
        \begin{align*}
            ( x, 0, \tfrac{\mathrm{d}}{\mathrm{d}t}x, 0) \in \calL_2. 
        \end{align*}
        More generally, it holds that $\calL_2=\calL_2^{\mrperp}$ as the following theorem shows. This connects extension theory for linear
relations with the geometric description of systems via Lagrangian subspaces. For a similar characterization for Dirac structures, we refer the reader to \cite{kurula2010dirac}.
    \end{rem}
\begin{thm}
\label{thm:lagrange}
    Let $\calF=\calX \times \calG$ for Hilbert spaces $\calX, \calG$, let $U=\left[\begin{smallmatrix}
            I_\calX &0 \\ 0& -I_\calG
        \end{smallmatrix} \right]:\calF \to \calF$ and define $\calB=\calF \times \calF$. Further, let $\left[\begin{smallmatrix}
        \Gamma_0 \\ \Gamma_1
    \end{smallmatrix}\right]: \dom \left[\begin{smallmatrix}
        \calP \\ \calS
    \end{smallmatrix}\right] \to \calG \times \calG$ be a map. Then, $(\calG, \Gamma_0, \Gamma_1)$ is a boundary triplet for $H^*=\ran \left[\begin{smallmatrix}
        \calP \\ \calS
    \end{smallmatrix}\right]$ if and only if
    \begin{equation*}
        \calL=\left\lbrace \left(\left[\begin{smallmatrix}
            f \\ \chi_\partial
        \end{smallmatrix} \right], \left[\begin{smallmatrix}
            e \\ \varepsilon_\partial
        \end{smallmatrix} \right]\right) \in \calB \, \middle| \, \exists \, x \in \dom \left[\begin{smallmatrix}
        \calP \\ \calS
    \end{smallmatrix}\right]: \left[\begin{smallmatrix}
            f \\ e
        \end{smallmatrix} \right]=\left[\begin{smallmatrix}
        \calP \\ \calS
    \end{smallmatrix}\right] x, \, \left[\begin{smallmatrix}
            \chi_\partial \\ \varepsilon_\partial 
        \end{smallmatrix} \right] = \left[\begin{smallmatrix}
        \Gamma_0 \\ \Gamma_1
    \end{smallmatrix}\right] x \right\rbrace
    \end{equation*}
    satisfies $\calL=\calL^{\mrperp}$.
\end{thm}
\begin{proof}
    We successively verify both implications. To this end, let $(\calG, \Gamma_0, \Gamma_1)$ be a boundary triplet for $\ran \left[\begin{smallmatrix}
        \calP \\ \calS
    \end{smallmatrix}\right]$. As a consequence, we have 
    \begin{align*}
             & \left\langle \hspace{-1.5mm} \left\langle  \left[\begin{smallmatrix}
            f_1\\ \chi_\partial^1 \\ e^1 \\ \varepsilon_\partial^1
        \end{smallmatrix} \right], \left[\begin{smallmatrix}
            f_2 \\ \chi_\partial^2 \\ e^2 \\ \varepsilon_\partial^2
        \end{smallmatrix} \right] \right\rangle \hspace{-1.5mm}\right\rangle_- 
        = \left\langle \hspace{-1.5mm} \left\langle   \left[\begin{smallmatrix}
            \calP \\ \Gamma_0  \\ \calS  \\ \Gamma_1 
        \end{smallmatrix} \right] x,  \left[\begin{smallmatrix}
            \calP \\ \Gamma_0  \\ \calS  \\ \Gamma_1 
        \end{smallmatrix} \right] y \right\rangle \hspace{-1.5mm}\right\rangle_- \\
        &= \langle \calS x , \calP y \rangle_\calX - \langle \calP x , \calS y \rangle_\calX  - (\langle \Gamma_1 x , \Gamma_0 y \rangle_\calG - \langle \Gamma_0 x, \Gamma_1 y \rangle_\calG) \\
        &=0
    \end{align*}
    for all $(f^1,\chi_\partial^1,e^1,\varepsilon_\partial^1), (f^2,\chi_\partial^2,e^2,\varepsilon_\partial^2) \in \calL$ which yields $\calL \subset \calL^{\mperp}$. Additionally, let an element $(f^2,\chi_\partial^2,e^2,\varepsilon_\partial^2) \in \calL^{\mperp}$, i.e., 
    \begin{equation*}
         \left\langle \hspace{-1.5mm} \left\langle  \left[\begin{smallmatrix}
            f_1\\ \chi_\partial^1 \\ e^1 \\ \varepsilon_\partial^1
        \end{smallmatrix} \right], \left[\begin{smallmatrix}
            f_2 \\ \chi_\partial^2 \\ e^2 \\ \varepsilon_\partial^2
        \end{smallmatrix} \right] \right\rangle \hspace{-1.5mm}\right\rangle_- =0
    \end{equation*}
    for all $(f_1,\chi_\partial^1,e^1,\varepsilon_\partial^1) \in \calL$. In particular, this implies
    \begin{equation*}
        \left\langle  \left[\begin{smallmatrix}
            \calP \\ \calS
        \end{smallmatrix} \right] x , \left[ \begin{smallmatrix}
            0& -I \\ I &0
        \end{smallmatrix} \right] \left[\begin{smallmatrix}
            f^2 \\ e^2
        \end{smallmatrix} \right] \right\rangle_{\calX \times \calX}=
        \langle \calS x , f^2 \rangle_\calX - \langle \calP x , e^2 \rangle_\calX  =0
    \end{equation*}
    for all $x \in \ker \Gamma$ which is equivalent to \[\left[\begin{smallmatrix}
            f^2 \\ e^2
        \end{smallmatrix} \right] \in (\overline{H})^{\mperp}=H^*=\ran  \left[\begin{smallmatrix}
            \calP \\ \calS
        \end{smallmatrix} \right].
        \]
        Furthermore, there exists $y \in \dom  \left[\begin{smallmatrix}
            \calP \\ \calS
        \end{smallmatrix} \right]$ such that $\left[\begin{smallmatrix}
            f^2 \\ e^2
        \end{smallmatrix} \right] = \left[\begin{smallmatrix}
            \calP \\ \calS
        \end{smallmatrix} \right]y$ and hence
        \begin{align*}
             &   \left\langle \hspace{-1.5mm} \left\langle  \left[\begin{smallmatrix}
            \calP \\ \Gamma_0  \\ \calS  \\ \Gamma_1 
        \end{smallmatrix} \right] x, \left[\begin{smallmatrix}
            \calP y\\ \chi_\partial^2  \\ \calS y  \\ \varepsilon_\partial^2 
        \end{smallmatrix} \right]  \right\rangle \hspace{-1.5mm}\right\rangle_- \\
        &= \langle \calS x , \calP y \rangle_\calX - \langle \calP x , \calS y \rangle_\calX - (\langle \Gamma_1 x , \chi_\partial^2 \rangle_\calG - \langle \Gamma_0 x, \varepsilon_\partial^2 \rangle_\calG) \\
        &=\langle \Gamma_1 x , \Gamma_0 y \rangle_\calG - \langle \Gamma_0 x, \Gamma_1 y \rangle_\calG - (\langle \Gamma_1 x , \chi_\partial^2 \rangle_\calG - \langle \Gamma_0 x, \varepsilon_\partial^2 \rangle_\calG) \\
        &= \langle \Gamma_1 x , \Gamma_0 y - \chi_\partial^2 \rangle_\calG - \langle \Gamma_0 x, \Gamma_1 y - \varepsilon_\partial^2 \rangle_\calG 
        \end{align*}
        for all $x \in \dom \left[\begin{smallmatrix}
            \calP \\ \calS
        \end{smallmatrix} \right]$. The surjectivity of $\Gamma$ now implies that $\left[\begin{smallmatrix}
            \chi_\partial^2 \\ \varepsilon_\partial^2
        \end{smallmatrix} \right] =\Gamma y$ which yields $(f^2,\chi_\partial^2,e^2,\varepsilon_\partial^2) \in \calL$. \par
        Conversely, let $\calL=\calL^{\mperp}$. Following the begin of this proof, the symmetry $\calL \subset \calL^{\mperp}$ immediately shows the abstract Green identity \eqref{eq:Green}. It remains to show that $\Gamma$ is surjective. Let $\left[\begin{smallmatrix}
           -\varepsilon_\partial^2 \\ \chi_\partial^2
        \end{smallmatrix} \right] \in \calG \times \calG$ be orthogonal to $\ran \Gamma$. Thus,
        \begin{equation*}
        \left\langle \hspace{-1.5mm} \left\langle  \left[\begin{smallmatrix}
            \calP \\ \Gamma_0  \\ \calS  \\ \Gamma_1 
        \end{smallmatrix} \right] x, \left[\begin{smallmatrix}
            0 \\ \chi_\partial^2  \\ 0 \\ \varepsilon_\partial^2 
        \end{smallmatrix} \right]  \right\rangle \hspace{-1.5mm}\right\rangle_- =\langle \Gamma_1 x , \chi_\partial^2 \rangle_\calG - \langle \Gamma_0 x, \varepsilon_\partial^2 \rangle_\calG=0
        \end{equation*}
        for all $x \in \dom  \left[\begin{smallmatrix}
            \calP \\ \calS
        \end{smallmatrix} \right]$ and hence $(0,\chi_\partial^2,0,\varepsilon_\partial^2) \in \calL^{\mperp}=\calL$. Lemma \ref{lem:injective_and_closed_range} implies that $\ker \left[\begin{smallmatrix}
            \calP \\ \calS
        \end{smallmatrix} \right]=\{0\}$ holds, and therefore, we conclude $\left[\begin{smallmatrix}
           \chi_\partial^2 \\ \varepsilon_\partial^2
        \end{smallmatrix} \right] = \Gamma 0 = 0$ and consequently, $\ran \Gamma = \calG \times \calG$ which yields that $(\calG, \Gamma_0, \Gamma_1)$ is a boundary triplet for $H^*=\ran \left[\begin{smallmatrix}
        \calP \\ \calS
    \end{smallmatrix}\right]$.
\end{proof}

\subsection{Boundary triplets for skew-symmetric operators}\label{subsec:bdd_triplets_skew}

Our second contribution of this work will be the application of the previously mentioned extension theory to port-Hamiltonian systems of the form \eqref{eq:pHs}.

In addition to the Lagrange subspaces (corresponding to self-adjoint relations as shown in the previous Subsection~\ref{subsec:Lagrange}), port-Hamiltonian systems additionally include Dirac structures, see also \eqref{eq:pHs_dyn}. Thus, we briefly recall boundary triplets for general skew-symmetric operators based on \cite{wegner2017boundary} and \cite[Chap. 2.4]{skrepek2021linear}. 
We begin with the definition of a boundary triplet. Note that, without loss of generality, we assume in the following that the minimal operator $\calJ_0$ is closed. 
\begin{defn}\label{def:operator_bdd_triplet}
    Let $\calJ_0: \calX \subset \dom \calJ_0 \to \calX$ be a densely defined and skew-symmetric operator. We call $(\calG, \Gamma_0, \Gamma_1)$ a boundary triplet for $\calJ_0^*:\calX \subset \dom \calJ_0^* \to \calX$ if $\calG$ is a Hilbert space, $\Gamma = \left[ \begin{smallmatrix}
        \Gamma_0 \\ \Gamma_1
    \end{smallmatrix}\right]: \dom \calJ_0^* \to \calG \times \calG$ is linear and surjective and the \textit{abstract Green identity} 
    \begin{equation}\label{eq:Green_for_skew}  \left\langle \calJ_0^* x , y \right\rangle_\calX + \left\langle x , \calJ_0^* y \right\rangle_\calX= \left\langle \Gamma_1 x, \Gamma_0 y \right\rangle_\calG + \left\langle \Gamma_0 x, \Gamma_1 y\right\rangle_\calG 
    \end{equation}
    holds for all $x,y \in \dom \calJ_0^*$. The space $\calG$ is called the \textit{boundary space}.
\end{defn}
Note that the main difference in comparison to symmetric operators (or relations) is the plus sign instead of the minus sign in \eqref{eq:Green_for_skew}. As for symmetric operators one recovers the minimal operator by restricting the (negative) adjoint to the kernel of the boundary map.
Moreover, it turns out that any operator $\euJ$ with $-\calJ_0 \subset \euJ \subset \calJ_0^*$ is given by $\euJ_\Theta$ for a suitable $\Theta \subset \calG \times \calG$. As for symmetric relations one may characterize any such skew-adjoint (maximal dissipative) operator via the relations in the boundary space. We now state an analogous result to Theorem~\ref{thm:intermediate_extensions} and Theorem~\ref{thm:self-adj_extension}.
\begin{prop}\label{prop:operator_intermediate_extensions}
    Let $\calJ_0: \calX \subset \dom \calJ_0 \to \calX$ be a densely defined and skew-symmetric operator and $(\calG, \Gamma_0, \Gamma_1)$ a boundary triplet for $\calJ_0^*$. Then, \begin{align}\label{eq:claim1}
        \calJ_0=-\calJ_0^* \left. \right|_{\ker \Gamma_0 \cap \ker \Gamma_1}.
    \end{align}
    Further, for $\Theta \subset \calG \times \calG$ define $\euJ_\Theta : \calX \subset \dom \euJ_\Theta \to \calX$ via
    \begin{equation}
        \begin{array}{rcl}
            \dom \euJ_\Theta & \coloneqq & \left\lbrace x \in \dom  \calJ_0^* \,  \middle| \, \Gamma x \in \Theta \right\rbrace = \Gamma^{-1} \Theta \label{eq:skew_intermediate_dom}\\
            \euJ_\Theta x &=& \calJ_0^* x
        \end{array}
    \end{equation}
    Then the following hold:
    \begin{itemize}
        \item[(i)] $\overline{\euJ_\Theta}=\euJ_{\overline{\Theta}}$
        \item[(ii)] $-\euJ_\Theta^*=\euJ_{-\Theta^*}$
        \item [(iii)] $\euJ_\Theta$ is skew-adjoint (maximal dissipative) if and only if $\Theta$ is skew-adjoint (maximal dissipative).
    \end{itemize}
\end{prop}
\begin{proof}
    The identity \eqref{eq:claim1} follows from~\cite[Lem.~2.4.5]{skrepek2021linear} and the remainder is proven in~\cite[Prop. 2.4.10]{skrepek2021linear}. 
\end{proof}
\section{Application to implicit port-Hamiltonian systems on one-dimensional spatial domains}\label{sec:pHs}
In this section, we apply our result to range representations involving operators $\calP$ and $\calS$ that are given as matrix differential operators with constant coefficients as in \eqref{eq:pHs}. The derivation of similar Green formulas in the case of non-constant coefficients will be future work. After recalling results on matrix differential operators, we construct a boundary triplet tailored to the setting of Section~\ref{sec:ext_theory} to characterize the boundary configurations leading to self-adjoint resp.\ maximally dissipative relations. The section concludes with examples such as the biharmonic wave equation and an elastic rod governed by a non-local elasticity condition.

\subsection{Basic properties of matrix differential operators}
On the Hilbert space $\calX \coloneqq L^2(a,b; \C^n)$, i.e.\ the space of square integrable functions defined on an open interval $(a,b)\subseteq\R$ with values in $\C^n$, we consider the minimal matrix differential operator $\euA_0: \calX \supset \dom \euA_0 \to \calX$ given by 
\begin{equation}\label{eq:A} 
    \euA_0f \coloneq \sum\limits_{k = 0}^{M} A_k \dif{k} f, \qquad \dom \euA_0 \coloneqq C_c^\infty(a,b;\C^n) 
\end{equation}
for matrices $A_k \in \C^{n \times n}$, $k=0,\ldots,M$. Unlike the previous studies \cite{augner2016stabilisation,villegas2007port}, we do not assume invertibility of the leading coefficient $A_M$.

To compute the adjoint operator, we briefly recall the notion of the  distributional derivative. Define the space of distributions $\calD'$ as the antidual of the space $\calD \coloneqq C_c^\infty(a,b, \C^n)$. For a distribution $T \in \calD'$ and a test function $\phi \in \calD$, we define
\begin{equation*}
    \langle T , \phi \rangle_{\calD', \calD} \coloneq T(\phi).
\end{equation*}
In particular, we can associate to $f \in L^1_{\operatorname{loc}}(a,b;\C^n)$ a distribution $T_f$, see e.g.\ \cite[p. 48]{yosida2012functional}, 
 by  setting
    \begin{equation}\label{eq:reg_distr}
        \langle T_f , \phi \rangle_{\calD', \calD} = \langle f , \phi \rangle_{\calX}\quad \text{for all $\phi \in \calD$.}
    \end{equation}
    Such distributions $T_f$ are called \emph{regular} and slightly abusing the notation we also denote the corresponding distribution $T_f$ by the same symbol, i.e., $f$. By considering \textit{distributional derivatives}, the operator $\euA_0$ given by  \eqref{eq:A} can be extended to the space of distributions mapping $T \in \calD'$ to  $\euA_0 T \in \calD'$, see e.g.\ \cite{schmudgen2012unbounded}, defined via
\begin{equation}\label{eq:distr_deriv}
    \langle \euA_0 T, \phi \rangle_{\calD', \calD} = \sum\limits_{k = 0}^{M} (-1)^k \left\langle A_kT, \tdif{k}\phi \right\rangle_{\calD', \calD} \quad \text{for all $\phi \in \calD$.}
\end{equation}
\begin{lem}\label{lem:A_star}
Let $\euA_0$ be given by \eqref{eq:A}. Then the adjoint $\euA_0^* \coloneqq \euB$ is the operator
    \begin{equation}\label{eq:A_star}
    \euB f \coloneq \sum\limits_{k=0}^{M} (-1)^k\dif{k}A_k^H  f, \qquad \dom \euB \coloneqq \{ f \in \calX  \mid \euB f \in \calX \}
\end{equation}
where $\euB f$ is meant in the sense of distributions.
\end{lem}
\begin{proof}
    We know that $f \in \dom \euB$ if and only if $f \in \calX$ and $\euB f \in \calX$ in the sense of distributions. Since $f, \euB f \in \calX$ we conclude that $f, \euB f \in L^1_{\operatorname{loc}}(a,b;\C^n)$. Hence, in view of \eqref{eq:reg_distr}, we see
    \begin{equation}
    \langle \euB f, \phi \rangle_{\calX}
    = \left\langle f, \sum\limits_{k = 0}^{M} A_k\dif{k} \phi \right\rangle_{\calX}=\langle f, \euA_0 \phi \rangle_{\calX}
\end{equation}
for all $\phi \in \calD$ and $f \in \dom \euB$. This shows $\dom \euB \subset \dom \euA_0^*$ and $\euB f=\euA_0^* f$ for all $f \in \dom \euB$. Otherwise, let $f \in \dom \euA_0^*$, i.e., there exists a $z \in \calX$ such that 
\begin{equation*}
    \langle z, \phi \rangle_{\calD', \calD} =\langle z, \phi \rangle_\calX = \langle f, \euA_0 \phi \rangle_\calX = \langle f , \euA_0 \phi \rangle_{\calD', \calD}  = \langle \euB f, \phi \rangle_{\calD', \calD} 
\end{equation*}
for all $\phi \in \calD$. As a consequence, $\euB f=z\in \calX$ and therefore $f \in \dom \euB$ holds, which implies $\euB=\euA_0^*$.
\end{proof}

\begin{lem}\label{lem:maximal_closed_densely_def}
The operator $\euB: \calX \supset \dom \euB \to \calX$ defined in \eqref{eq:A_star} is closed and densely defined.
\end{lem}
\begin{proof}
Since $\euB$ is the adjoint operator of $\euA_0$ it is closed, which proves the first claim. For the second claim, we observe that for 
$T \in \calD'$, we have
\[
\langle \euB T , \phi\rangle_{\calD',\calD} 
= \sum_{k=0}^M (-1)^k \langle A_kT, \tdif{k}\phi\rangle_\calX 
\]
for all $\phi \in \calD$. Hence, for the regular distribution $T_f$ corresponding to $f\in L^1_{\rm loc}(a,b;\C^n)$ 
we have $\langle \euB T_f , \phi\rangle_{\calD',\calD} = \sum_{k=0}^M (-1)^k\langle T_{A_kf}, \phi^{(k)}\rangle_{\calD',\calD}$. The fact that $f\in\dom \euB$ thus means that $f\in \calX$ such that $\euB T_f = T_g$ with some $g \in \calX$, that is
\begin{equation}\label{e:samstag}
\sum_{k=0}^N(-1)^k \langle \vphi^{(k)},A_k f\rangle_\calX= \langle \phi ,g\rangle_\calX
\end{equation}
for all $\phi \in \calD$. From here, it is easy to see that $C_c^\infty(a,b;\C^n) \subset\dom \euB$ with $\euB f = \sum_{k=0}^M A_k \tdif{k}f$ for $f \in C_c^\infty(a,b;\C^n)$. In particular, $\euB$ is densely defined.
\end{proof}
\noindent We define the space 
\begin{equation*}
    H(\euB; a,b) \coloneqq \dom \euB = \{ f \in \calX  \mid \euB f \in \calX \},
\end{equation*}
in analogy to the notation in \cite[Sec. 3.1]{skrepek2021linear}, \cite[Eq. 6.56]{villegas2007port}. It follows from Lemma \ref{lem:maximal_closed_densely_def} that $H(\euB; a,b)$ is a Hilbert space when endowed with the graph norm
\begin{equation*}
    \|x\|_\euB \coloneqq \left( \|x\|_\calX^2 + \|\euB x\|_\calX^2 \right)^{1/2}.
\end{equation*}

\subsection{Construction of boundary triplets for  matrix differential operators}
First, we define the matrix differential operators 
\begin{align*}
    \calP_0, \calS_0: \calX \supset C_c^\infty(a,b;\C^n) \to \calX    
\end{align*}
of even order via
\begin{equation}\label{eq:pH_operators}
    \calP_0 x = \sum\limits_{k=0}^N  \dif{k} P_k^H \dif{k}x, \qquad \calS_0 x = \sum\limits_{l=0}^N  \dif{l} S_l^H \dif{l}x 
\end{equation}
for matrices $P_k,S_l \in \C^{n \times n}$, with $k,l=0,\ldots,N$. Note that, without loss of generality and as we do not include an invertibility condition, we assumed both differential operators to have the same degree. 

Following Lemma \ref{lem:A_star} we observe that the Hilbert space adjoints of the minimal operators are given by
\begin{equation}\label{eq:matrix_diff_op}
    \calP x \coloneqq \calP_0^* x = \sum\limits_{k=0}^N  \dif{k} P_k \dif{k} x, \qquad \calS x \coloneqq \calS_0^* x = \sum\limits_{l=0}^N \dif{l} S_l \dif{l} x ,
\end{equation}
where 
\begin{equation*}
    \dom \calP =H(\calP;a,b) 
    ,\qquad \dom \calS = H(\calS;a,b).
\end{equation*}
Recall that $\dom \left[\begin{smallmatrix}
        \calP \\ \calS
    \end{smallmatrix}\right]= H(\calP; a,b) \cap H(\calS; a,b)$. In the terminology of Section~\ref{subsec:bdd_triplets}, we define $H=(\gr \calP_0)^{-1} (\gr \calS_0)$. 

To apply our results from Section~\ref{sec:ext_theory}, we have to verify the central Assumption~\ref{ass:ass1}. 
To this end, we check the coercivity condition \eqref{eq:coercivity}. We first note that \eqref{eq:coercivity} holds if one of the operators $\calS$ or $\calP$ is boundedly invertible or coercive which is satisfied e.g.\ for the biharmonic wave equation or the elasic rod with non-local elasticity as discussed in Section~\ref{sec:applications}. The following result provides an alternative sufficient condition for the coercivity condition~\eqref{eq:coercivity}. 
\begin{prop}\label{prop:suff_coercivity}
Let $\calP$ and $\calS$ be given by \eqref{eq:matrix_diff_op}. Consider $\lambda_k=\left(\frac{k\pi}{b-a}\right)^2$, $k\geq 1$ and assume that for
\[
c_{S,k}:=\sigma_{\min}\left(\sum\limits_{j=0}^N  \lambda_k^{2j}S_j\right)\geq 0,\quad c_{P,k}:=\sigma_{\min}\left(\sum\limits_{j=0}^N  \lambda_k^{2j}P_j\right)\geq 0
\]
there is $c>0$ such that $c_{S,k}^2+c_{P,k}^2\geq c^2$ for all $k\geq 1$.
Then the coercivity condition~\eqref{eq:coercivity} holds for $\left[\begin{smallmatrix}
        \calP \\ \calS
    \end{smallmatrix}\right]$ as in \eqref{eq:matrix_diff_op}. 
    \end{prop}
    \begin{proof}
We examine $\left\lbrace\varphi_k(\cdot)=c_k\sin(\frac{k\pi}{b-a}\cdot)\right\rbrace_{k=1}^\infty$, where $c_k>0$ ensures $\|\varphi_k\|=1$, forming an orthonormal basis for $L^2(a,b)$. 
Hence, for any $x\in L^2(a,b;\C^n)$ we have the decomposition 
$x=\sum_{k=1}^{\infty}x_{k}\varphi_k$ where $x_k=(x_{k,i})_{i=1}^n\in\C^n$ is given by $x_{k,i} = \langle x^{(i)},\varphi_k\rangle_\mathcal{X}\in\mathbb{C}^n$ and $x^{(i)}$ denotes the $i$-th component of the vector-valued function $x$, such that we have the Parseval identity $\|x\|^2=\sum_{k}\|x_{k}\|^2$. For $k\in \N$, we compute
\begin{align}
\nonumber
\langle\calS x_k\varphi_k,\calS x_k\varphi_k\rangle &= 
\left\langle \sum\limits_{j=0}^N \tdif{j} S_j \tdif{j}x_k\varphi_k,\sum\limits_{j=0}^N \tdif{j} S_j \tdif{j}x_k\varphi_k\right\rangle\\
&=\left\langle \sum\limits_{j=0}^N  \lambda_k^{2j}S_jx_k\varphi_k,\sum\limits_{j=0}^N  \lambda_k^{2j}S_jx_k\varphi_k\right\rangle \nonumber\\
&=\left|\hspace{-0.4mm}\left|\sum\limits_{j=0}^N  \lambda_k^{2j}S_jx_k\right|\hspace{-0.4mm}\right|^2\geq c_{S,k}^2\left|\hspace{-0.3mm}\left|x_k\right|\hspace{-0.3mm}\right|^2 \label{eq:x_k_lower}
\end{align}
and by orthonormality and considering a finite sum $x=\sum_{k=1}^Mx_k\varphi_k$ we obtain
\begin{align}
\nonumber
\langle\calS x,\calS x\rangle&=\langle\calS \sum_{k=1}^Mx_k\varphi_k,\calS \sum_{k=1}^Mx_k\varphi_k\rangle\\ &=\left\langle \sum_{k=1}^M\calS x_k\varphi_k, \sum_{k=1}^M\calS x_k\varphi_k\right\rangle \nonumber\\
&=\left\langle \sum\limits_{k=1}^M\sum\limits_{j=0}^N  \lambda_k^{2j}S_jx_k\varphi_k,\sum\limits_{k=1}^M\sum\limits_{j=0}^N  \lambda_k^{2j}S_jx_k\varphi_k\right\rangle\nonumber\\
&=\left\langle \sum\limits_{k=1}^M \varphi_k\sum\limits_{j=0}^N  \lambda_k^{2j}S_jx_k,\sum\limits_{k=1}^M\varphi_k\sum\limits_{j=0}^N  \lambda_k^{2j}S_jx_k\right\rangle\nonumber\\
&=\sum_{k=1}^M\left|\hspace{-0.4mm}\left|\sum\limits_{j=0}^N  \lambda_k^{2j}S_jx_k\right|\hspace{-0.4mm}\right|^2\nonumber\\&\geq \sum_{k=1}^Mc_{S,k}^2\left|\hspace{-0.3mm}\left|x_k\right|\hspace{-0.3mm}\right|^2. \label{eq:s_lower}
\end{align}
Similarly, we estimate for $x=\sum_{k=1}^Nx_k\varphi_k$
\begin{align}
\label{eq:p_lower}
\langle\calP x,\calP x\rangle\geq \sum_{k=1}^M\left|\hspace{-0.4mm}\left|\sum\limits_{j=0}^N  \lambda_k^{2j}P_jx_k\right|\hspace{-0.4mm}\right|^2 \geq  \sum_{k=1}^Mc_{P,k}^2\left|\hspace{-0.3mm}\left|x_k\right|\hspace{-0.3mm}\right|^2.
\end{align}
Hence, the form 
\[
\mathfrak{t}: \dom \left[\begin{smallmatrix}
        \calP \\ \calS
    \end{smallmatrix}\right] \times \dom \left[\begin{smallmatrix}
        \calP \\ \calS
    \end{smallmatrix}\right] \to \C,\quad \mathfrak{t}(x,y):=\langle \calP x, \calP y\rangle_{\calX}+\langle \calS x, \calS y\rangle_{\calX}
\] satisfies
\begin{align}
\label{eq:p_s_combined}    
\mathfrak{t}(x,x)\geq \sum_{k=1}^M(c_{S,k}^2+c_{P,k}^2)\|x_k\|^2\geq c^2\sum_{k=1}^M\|x_k\|^2
\end{align}
for $x=\sum_{k=1}^Mx_k\varphi_k$ by combining the estimates \eqref{eq:s_lower} and \eqref{eq:p_lower}. The result for general $x\in L^2(a,b;\C^n)$ follows by passing to the limit $M\rightarrow\infty$.
    \end{proof}
    
Because in Proposition \ref{prop:suff_coercivity} we only provide a sufficient condition for coercivity, we assume that the coercivity condition~\eqref{eq:coercivity} holds for $\left[\begin{smallmatrix}
        \calP \\ \calS
    \end{smallmatrix}\right]$. Hence, Assumption~\ref{ass:ass1} is satisfied and then by Lemma \ref{lem:injective_and_closed_range} this implies that $\ran \left[\begin{smallmatrix}
        \calP \\ \calS
    \end{smallmatrix}\right]$ is closed in $\calX \times \calX$.  Thus, with Lemma \ref{lem:adjoint_rel}, $H^*=(\gr \calS) (\gr \calP)^{-1}=\ran \left[\begin{smallmatrix}
        \calP \\ \calS
    \end{smallmatrix}\right]$ holds. 
For a construction of a boundary triplet for $\ran \left[\begin{smallmatrix}
        \calP \\ \calS
    \end{smallmatrix}\right]$ in terms of the boundary evaluations we need another assumption.
\begin{ass}\label{ass:mdo}
Let $\calP$ and $\calS$ given by \eqref{eq:matrix_diff_op} satisfy the  following conditions:
    \begin{enumerate}
        \item[(i)] the coefficients $P_k,S_l \in \C^{n \times n}, k,l\hspace{-0.5mm}=\hspace{-0.5mm}0,\ldots,N$ satisfy
        for all $p\hspace{-0.5mm}=\hspace{-0.5mm}0,\ldots,2N$ 
        \begin{equation}\label{eq:symmetric_matrices}
        \sum\limits_{m=0}^p S_m^H P_{p-m} = \left( \sum\limits_{m=0}^p S_m^H P_{p-m} \right)^H, 
    \end{equation}
    where $P_m = S_m =0$ if $m>N$;
        \item[(ii)] there is at least one pair $(P_k,S_l)\neq (0,0)$ with $k\neq l$. 
    \end{enumerate}
\end{ass}
Note that in the case of $\calP=I$, the first assumption in \ref{ass:mdo} is equivalent to symmetry of the operator $\calS_0$. Moreover, the second assumption in \ref{ass:mdo} guarantees that we are able to apply the integration by parts formula.
\begin{rem}
The authors were not able to give a characterization of symmetry of $H$ in terms of the coefficient matrices of $\calP_0, \calS_0$, cf. Proposition \ref{prop:H_symm}. Instead we assume \eqref{eq:symmetric_matrices} that ensures Equation \eqref{eq:symm_prop_2} to hold on $C_c^\infty(a,b;\C^n)$. In \cite{MascvdSc23} it was assumed that the condition 
    \begin{equation}\label{eq:green_without_bdd}
        \left\langle \calS x , \calP y \right\rangle_\calX = \left\langle \calP x , \calS y \right\rangle_\calX
    \end{equation}
    holds for all $x,y \in C_c^\infty(a,b;\C^n)$. Given this, the authors were able to construct the so-called boundary port variables (boundary maps) associated with $H^*$. It can be easily shown by using the integration by parts formula that this so-called \textit{Maxwells reciprocity condition} \eqref{eq:green_without_bdd} is equivalent to the condition \eqref{eq:symmetric_matrices} in Assumption~\ref{ass:mdo}.  We refer the reader to \cite[Prop. 23]{MascvdSc23} where the Lagrangian subspace associated with $\calP, \calS$ is constructed, cf. Section~\ref{subsec:Lagrange}. 
\end{rem}
Now, we are almost in the position to explicitly construct a boundary triplet for $H^*= \ran \left[\begin{smallmatrix}
        \calP \\ \calS
    \end{smallmatrix}\right]$ after providing the following technical auxiliary result.
\begin{lem}
\label{lem:technical}
\allowdisplaybreaks
Consider the operators $\calP$ and $\calS$ given by \eqref{eq:matrix_diff_op} and let 
Assumption~\ref{ass:mdo} be satisfied. Then there exists $B\in\mathbb{C}^{2Nn}$ such that 
\begin{align}\label{eq:Green_with_AB}    
  \left\langle \calS x , \calP y \right\rangle -   \left\langle \calP x , \calS y \right\rangle &=
   \left\langle \begin{bmatrix}
            B-B^H &0 \\ 0& -(B-B^H)
        \end{bmatrix} \gamma x, \gamma y \right\rangle_{\C^{4Nn}} 
\end{align}
holds for all $x,y\in H^{2N}(a,b; \C^n)$, where 
$\gamma: H^{2N}(a,b;\C^n) \to \C^{4Nn}$ denotes the trace operator given by \eqref{eq:trace}. 
\end{lem}
\begin{proof}
     For readability, we let $\left\langle \cdot , \cdot \right\rangle \coloneqq \left\langle \cdot , \cdot \right\rangle_\calX$. Observe that for all $x,y \in H^{2N}(a,b; \C^n)$, we have 
    \begin{equation}
    \label{eq:S_P_hij}
             \left\langle \calS x , \calP y \right\rangle = \sum_{i,j=0}^{N} \underbrace{\left\langle \tdif{i}S_i \tdif{i} x , \tdif{j} P_j \tdif{j} y \right\rangle}_{\eqqcolon h_{ij}(x,y)}.
    \end{equation}
In the following, we apply the integration by parts formula to each of the expressions $h_{ij}(x,y)$ until both sides of the inner product have the same order of differentiation. Precisely, we have
    \begin{equation}
    \label{eq:hij_rewritten}
       h_{ij}(x,y) = 
        (-1)^{|i-j|} \left\langle \tdif{i+j} x , \tdif{i+j} S_i^HP_j y \right\rangle +  B_{S_iP_j}(x,y),
    \end{equation}
    where $B_{S_iP_j}: H^{2N}(a,b; \C^n) \times H^{2N}(a,b; \C^n) \to \C$ is defined by
    \begin{equation*}
        B_{S_iP_j}(x,y)=\scalebox{0.9}{$\begin{cases}
            \left[ \sum\limits_{k=0}^{|i-j|-1}(-1)^k \langle \tdif{2i-1-k} x(\xi)  , \tdif{2j+k} S_i^HP_jy(\xi) \rangle_{\C^n}\right]_{\xi=a}^{\xi=b} ,& \text{if} \ i > j, \\[2.5ex]
            \left[ \sum\limits_{k=0}^{|i-j|-1}(-1)^k \langle \tdif{2i+k} x(\xi)  , \tdif{2j-1-k} S_i^HP_jy(\xi) \rangle_{\C^n}\right]_{\xi=a}^{\xi=b} ,& \text{if} \ i < j, \\
            0, & \text{otherwise.}
        \end{cases}$}
    \end{equation*}
    In particular, 
\begin{equation}\label{eq:symmetric_boundary_form}
       B_{P_jS_i}(x,y) = \overline{B_{S_iP_j}(y,x)}
    \end{equation}
holds. Using $(-1)^{|i-j|}=(-1)^{i+j}$ and combining \eqref{eq:S_P_hij}, \eqref{eq:hij_rewritten}, and \eqref{eq:symmetric_boundary_form} yields 
\allowdisplaybreaks
    \begin{align}
            &\left\langle \calS x , \calP y \right\rangle \\
            &= \sum\limits_{i,j=0}^N 
        (-1)^{|i-j|} \left\langle \tdif{i+j} x , \tdif{i+j} S_i^HP_j y \right\rangle + B_{S_iP_j}(x,y) \nonumber\\
        &= \sum\limits_{p=0}^{2N} (-1)^p\left\langle \tdif{p} x , \tdif{p} \left(\sum\limits_{m=0}^{p} S_m^H P_{p-m}   \right) y \right\rangle + \sum\limits_{i,j=0}^N  B_{S_iP_j}(x,y) \nonumber\\
        &\!\!\!\stackrel{\eqref{eq:symmetric_matrices}}{=} \sum\limits_{p=0}^{2N} (-1)^p\left\langle \tdif{p} x , \tdif{p} \left(\sum\limits_{m=0}^{p} P_{p-m}^H S_{m}   \right) y \right\rangle + \sum\limits_{i,j=0}^N  B_{S_iP_j}(x,y) \nonumber\\
        &= \sum\limits_{i,j=0}^N 
        (-1)^{|j-i|} \left\langle \tdif{j+i} x , \tdif{j+i} P_j^HS_i y \right\rangle + \sum\limits_{i,j=0}^N  B_{S_iP_j}(x,y) \nonumber\\
        &= \sum\limits_{i,j=0}^N  \left\langle \tdif{j}P_j \tdif{j} x , \tdif{i} S_i \tdif{i} y \right\rangle + \sum\limits_{i,j=0}^N  B_{S_iP_j}(x,y) -  \sum\limits_{i,j=0}^N  B_{P_jS_i}(x,y) \nonumber\\
        &= \left\langle \calP x , \calS y \right\rangle + \sum\limits_{i,j=0}^N  B_{S_iP_j}(x,y) -  \sum\limits_{i,j=0}^N  B_{P_jS_i}(x,y),\nonumber\\
    &\!\!\!\stackrel{\eqref{eq:symmetric_boundary_form}}{=}\left\langle \calP x , \calS y \right\rangle + \underbrace{\sum\limits_{i,j=0}^N  B_{S_iP_j}(x,y) -  \sum\limits_{i,j=0}^N \overline{B_{S_iP_j}(y,x)}}_{\eqqcolon \mathfrak{a}(x,y)}. \label{eq:green_with_BPS}
    \end{align}
    We immediately see that $\mathfrak{a}: H^{2N}(a,b; \C^n) \times H^{2N}(a,b; \C^n) \to \C$ is an anti-symmetric sesquilinear form, i.e., $\mathfrak{a}(x,y)=-\overline{\mathfrak{a}(y,x)}$. 

    Next, we aim to rewrite the above sum over the boundary evaluations in terms of the trace operator $\gamma$ given by \eqref{eq:trace}. A direct calculation together with an index change $m=2i-1$ leads to 
    \begin{align}
            & \sum\limits_{i,j=0 }^N  B_{S_iP_j}(x,y) \nonumber  \\
            &= \sum\limits_{l=1}^N \sum\limits_{i=l}^N \sum\limits_{k=0}^{l-1} \left( (-1)^k \left[\langle \tdif{2i-1-k} x(\xi)  , S_i^HP_{i-l} \tdif{2(i-l)+k}  y(\xi) \rangle_{\C^n}\right]_{\xi=a}^{\xi=b} \right. \nonumber  \\
            & \hphantom{\sum\limits_{l=1}^N \sum\limits_{i=l}^N \sum\limits_{k=0}^{l-1}} \left. +(-1)^k \left[\langle \tdif{2(i-l)+k} x(\xi)  , S_{i-l}^HP_i \tdif{2i-1-k}  y(\xi) \rangle_{\C^n}\right]_{\xi=a}^{\xi=b}  \right) \nonumber  \\
            &= \sum\limits_{l=1}^N \sum\limits_{k=0}^{l-1} \sum\limits_{i=l}^N  \left( (-1)^k \left[\langle \tdif{2i-1-k} x(\xi)  , S_i^HP_{i-l} \tdif{2(i-l)+k}  y(\xi) \rangle_{\C^n}\right]_{\xi=a}^{\xi=b} \right. \nonumber \\
            & \hphantom{\sum\limits_{l=1}^N \sum\limits_{i=l}^N \sum\limits_{k=0}^{l-1}} \left. +(-1)^k \left[\langle \tdif{2(i-l)+k} x(\xi)  , S_{i-l}^HP_i \tdif{2i-1-k}  y(\xi) \rangle_{\C^n}\right]_{\xi=a}^{\xi=b}  \right)\nonumber \\
            &= \sum\limits_{l=1}^N \sum\limits_{k=0}^{l-1} \sum\limits_{\substack{m=2l-1 \\ m \text{ odd}}}^{2N-1} \left(\left[\langle \tdif{m-k} x(\xi)  , (-1)^k S_{\tfrac{m+1}{2}}^HP_{\tfrac{m+1-2l}{2}} \tdif{m+1-2l+k}  y(\xi) \rangle_{\C^n}\right]_{\xi=a}^{\xi=b} \right. \nonumber  \\
            & \hphantom{\sum\limits_{l=1}^N \sum\limits_{i=l}^N \sum\limits_{\substack{m=2l-1 \\ m \text{ odd}}}^{2N-1}} \left. +\left[\langle \tdif{m+1-2l+k} x(\xi)  , (-1)^k S_{\tfrac{m+1-2l}{2}}^HP_{\tfrac{m+1}{2}} \tdif{m-k}  y(\xi) \rangle_{\C^n}\right]_{\xi=a}^{\xi=b}  \right) \nonumber \\
            &= \sum\limits_{l=1}^N \sum\limits_{k=0}^{l-1} \sum\limits_{\substack{m=2l-1 \\ m \text{ odd}}}^{2N-1} \left[ 
            \begin{bmatrix} x(\xi)\\ \tdif{}x(\xi)\\ \vdots \\ \tdif{2N-1}x(\xi)\end{bmatrix}^* M_{lkm} \begin{bmatrix} y(\xi)\\ \tdif{}y(\xi)\\ \vdots \\ \tdif{2N-1}y(\xi)\end{bmatrix}\right]_{\xi=a}^{\xi=b}, \label{eq:BPS_almost_trace}
    \end{align}
    and where $M_{lkm}  \in  \C^{2Nn \times 2Nn}$ is a block matrix $M_{lkm}=(\widehat M_{i,j})_{i,j=1}^{2N}$ with $\widehat M_{i,j}\in\mathbb{C}^{n\times n}$ defined by
    \begin{equation*}
    \widehat M_{i,j}\coloneqq \begin{cases}
            (-1)^k  S_{\tfrac{m+1}{2}}^HP_{\tfrac{m+1-2l}{2}}, & \text{if} \ (i,j)=(m-k,m+1-2l+k) ,\\
            (-1)^k S_{\tfrac{m+1-2l}{2}}^HP_{\tfrac{m+1}{2}}, & \text{if} \ (i,j)=(m+1-2l+k, m-k), \\
            0, & \text{elsewhere.}
        \end{cases}
    \end{equation*}
    We define 
    \begin{equation}\label{eq:A_matrix}
        B \coloneqq \left( \sum\limits_{l=1}^N \sum\limits_{k=0}^{l-1} \sum\limits_{\substack{m=2l-1 \\ m \text{ odd}}}^{2N-1} M_{lkm} \right)^H \in \C^{2Nn \times 2Nn}   
    \end{equation}
    and using the trace operators introduced in \eqref{eq:trace}, \eqref{eq:BPS_almost_trace} and \eqref{eq:symmetric_boundary_form}  imply
    \begin{align}
       \mathfrak{a}(x,y)&= \sum\limits_{i,j=0 }^N  B_{S_iP_j}(x,y)-\sum\limits_{i,j=0 }^N  \overline{B_{S_iP_j}(y,x)} \nonumber \\ &= \left\langle \begin{bmatrix}
            B & 0\\0 & - B
        \end{bmatrix} \gamma x, \gamma y \right\rangle_{\C^{4Nn}}-\left\langle \begin{bmatrix}
            B & 0\\0 & - B
        \end{bmatrix}^H \gamma x, \gamma y \right\rangle_{\C^{4Nn}} .
        \label{eq:a_with_A}
    \end{align}
    A combination with \eqref{eq:green_with_BPS} implies \eqref{eq:Green_with_AB}.
\end{proof}

The main theorem concerning the existence of boundary triplets for pairs of matrix differential operators as defined by \eqref{eq:matrix_diff_op} can now be presented.
\begin{thm}\label{thm:boundary_triplet}
Consider the operators $\calP$ and $\calS$ given by \eqref{eq:matrix_diff_op} and let \eqref{eq:coercivity} be satisfied.  Then, Assumption~\ref{ass:ass1} is fulfilled and the following holds:
\begin{itemize}
    \item[\rm (i)] If additionally Assumption~\ref{ass:mdo} is fulfilled and $H \subset H^*$, then there exists a boundary triplet $(\calG, \widetilde{\Gamma}_0, \widetilde{\Gamma}_1)$ for $H^*= \ran \left[\begin{smallmatrix}  \calP \\ \calS  \end{smallmatrix}\right]$.
    \item[\rm (ii)] If $A \coloneqq B-B^H \in\C^{2Nn \times 2Nn}$ where $B$ is given by \eqref{eq:A_matrix} satisfies $\rank A=2Nn$, then we can choose the following  boundary triplet for $H^*$
   \begin{equation}
   \label{def:Gamma_full_rank}
         \Gamma x = \begin{bmatrix}
        \Gamma_0 \\ \Gamma_1
    \end{bmatrix} x = \frac{1}{\sqrt{2}} \begin{bmatrix}
        A & -A \\ -I & -I
    \end{bmatrix} \gamma x,
    \end{equation}
      where $\gamma: H^{2N}(a,b;\C^n) \to \C^{4Nn}$ denotes the trace operator given by \eqref{eq:trace}. Moreover, $\overline{H}=\left[\begin{smallmatrix}
        \calP \\ \calS
    \end{smallmatrix}\right]H_0^{2N}(a,b;\C^n)$ holds.
\end{itemize}
\end{thm}
\begin{proof}
First note that we assume that \eqref{eq:coercivity} holds and thus Assumption~\ref{ass:ass1} is fulfilled.
\begin{itemize}
\item[\rm (i)]
The matrix $A=B-B^H \in\C^{2Nn \times 2Nn}$ where $B$ is given by \eqref{eq:A_matrix} fulfills 
\[
  \left\langle \calS x , \calP y \right\rangle -   \left\langle \calP x , \calS y \right\rangle =\left\langle \begin{bmatrix}
            A &0 \\ 0& -A
        \end{bmatrix} \gamma x, \gamma y \right\rangle_{\C^{4Nn}}.
\]
    The matrix $A$ is skew-Hermitian by definition, and hence, the imaginary multiple $iA$ is Hermitian. 
    Thus, by the spectral theorem for Hermitian matrices, see e.g.\ \cite[Thm.~2.5.6]{HornJohn13}, there exist a unitary matrix $U \in \C^{2Nn \times 2Nn}$ and a diagonal matrix $D \in \C^{2Nn \times 2Nn}$ consisting of the purely imaginary eigenvalues of $A$ such that $A=UDU^H$. We let $U_{\mathrm{r}} \in \C^{2Nn \times \rank A}$ be the matrix whose columns are the eigenvectors corresponding to the nonzero eigenvalues of $A$. Further, we let $D_{\mathrm{r}} \in \C^{\rank A \times \rank A}$ be the diagonal matrix containing only the nonzero eigenvalues of $A$. Thus, we write the reduced factorization 
    \begin{equation*}
        A=U_{\mathrm{r}}D_{\mathrm{r}}U_{\mathrm{r}}^H.
    \end{equation*}
    We now define the boundary map $\widehat\Gamma$ using the trace operator $\gamma$
    \begin{align*}
 \widehat{\Gamma} =  \begin{bmatrix}
        \widehat\Gamma_0 \\ \widehat\Gamma_1
    \end{bmatrix} &: H^{2N}(a,b;\C^n) \to \C^{\rank A} \times \C^{\rank A},\\        \begin{bmatrix}
        \widehat\Gamma_0 \\ \widehat\Gamma_1
    \end{bmatrix} & x = \frac{1}{\sqrt{2}} \begin{bmatrix}
        D_{\mathrm{r}} & -D_{\mathrm{r}} \\ -I & -I
    \end{bmatrix} \begin{bmatrix}
        U_{\mathrm{r}}^H & 0  \\ 0 & U_{\mathrm{r}}^H
    \end{bmatrix} \gamma x.
    \end{align*}
   It is easily verified that
    \begin{equation*}
        \tfrac{1}{\sqrt{2}} \left( \left[ \begin{smallmatrix}
        D_{\mathrm{r}} & -D_{\mathrm{r}} \\ -I & -I
        \end{smallmatrix}\right] \left[\begin{smallmatrix}
        U_{\mathrm{r}}^H &0  \\ 0& U_{\mathrm{r}}^H
    \end{smallmatrix} \right]  \right)^H  \left[ \begin{smallmatrix}
        0 & -I \\ I &0
    \end{smallmatrix}\right] \left[ \begin{smallmatrix}
        D_{\mathrm{r}} & -D_{\mathrm{r}} \\ -I & -I
        \end{smallmatrix}\right] \left[\begin{smallmatrix}
        U_{\mathrm{r}}^H & 0 \\ 0& U_{\mathrm{r}}^H
    \end{smallmatrix} \right] \tfrac{1}{\sqrt{2}} =  \left[ \begin{smallmatrix}
            A^H& 0\\ 0& -A^H
        \end{smallmatrix}\right]
    \end{equation*}
    and therefore
    \begin{align}
            \left\langle  \left[ \begin{smallmatrix}
        \widehat\Gamma_0 \\ \widehat\Gamma_1
    \end{smallmatrix}\right] x, \left[ \begin{smallmatrix}
        0 & -I \\ I &0
    \end{smallmatrix}\right] \left[ \begin{smallmatrix}
        \widehat\Gamma_0 \\ \widehat\Gamma_1
    \end{smallmatrix}\right] y\right\rangle_{\C^{4Nn}}
    =\left\langle \gamma x,  \left[ \begin{smallmatrix}
        A^H &0  \\ 0 & -A^H
    \end{smallmatrix} \right]  \gamma y\right\rangle_{\C^{4Nn}} 
    =\mathfrak{a}(x,y). \label{eq:Green_with_A}
    \end{align}
    Observe that $U_{\mathrm{r}}^H$ has full row rank and $\gamma$ is surjective by Proposition~\ref{prop:trace_surjective}. Hence, $\Gamma$ is surjective as the composition of surjective maps is surjective. Furthermore, along the lines of \cite[Thm. 6.19]{villegas2007port}, it follows that $H^{2N}(a,b;\C^n)$ is dense in $H(\calP; a,b) \cap H(\calS; a,b)$ with respect to the topology induced by \eqref{eq:dom_inner_product}. Let $\widetilde{\Gamma}_0, \widetilde{\Gamma}_1: H(\calP; a,b) \cap H(\calS; a,b) \to \C^{\rank A}$ be the unique extensions of $\Gamma_0,\Gamma_1$. By Lemma \ref{lem:extension_by_density} we conclude that $(\calG, \widetilde{\Gamma}_0, \widetilde{\Gamma}_1)$ where $\calG=\C^{\rank A}$ defines a boundary triplet for $H^*= \ran \left[\begin{smallmatrix}
        \calP \\ \calS
    \end{smallmatrix}\right]$.
    
\item[\rm (ii)]  Since $A$ is invertible, the surjectivity of $\gamma$ that was shown in Proposition~\ref{prop:trace_surjective}  translates to $\Gamma$ that is defined in \eqref{def:Gamma_full_rank}. Moreover, the abstract Green identity is a direct consequence of \eqref{eq:Green_with_A} combined with the abstract Green identity for the boundary triplet $(\calG,\widehat\Gamma_0,\widehat\Gamma_1)$. 
Furthermore, by Proposition \ref{prop:trace_surjective}, the regularity of $A$ and Lemma \ref{lem:Gamma_properties} we conclude
    \begin{equation*}
        \overline{H}=\begin{bmatrix}
        \calP \\ \calS
    \end{bmatrix} \ker \Gamma = \begin{bmatrix}
        \calP \\ \calS
    \end{bmatrix} \ker \gamma = \begin{bmatrix}
        \calP \\ \calS
    \end{bmatrix}H_0^{2N}(a,b;\C^n).
    \end{equation*}
\end{itemize}
\end{proof}
\begin{rem}
    The choice of the boundary maps in Theorem~\ref{thm:boundary_triplet}~(iii) is inspired by \cite[Lem. 7.2.2]{JacoZwar12}, where boundary maps are constructed for a class of linear first-order port-Hamiltonian systems on one-dimensional domains. A similar characterization for higher-order matrix differential operators is given in \cite{MascvdSc23}. 
\end{rem}

In the remainder, we present some general results on intermediate extensions for finite-dimensional boundary spaces for which the parameterizing relations are also finite-dimensional subspaces in kernel representation. 
\begin{cor}\label{cor:kernel_domain}
     Let $(\C^{\rank A}, \widetilde{\Gamma_0}, \widetilde{\Gamma_1})$ be the boundary triplet for $H^*= \ran \left[\begin{smallmatrix}
        \calP \\ \calS
    \end{smallmatrix}\right]$ given by  Theorem~\ref{thm:boundary_triplet}. Moreover, let $K, L \in \C^{\rank A \times \rank A}$. The following assertions apply to the restriction $A_{\ker \left[\begin{smallmatrix}
        K & L
    \end{smallmatrix} \right]}$, cf. \eqref{eq:operator_theta}.
    \begin{itemize}
        \item[(i)] $\ran A_{\ker \left[\begin{smallmatrix}
        K & L
    \end{smallmatrix} \right]}$ is closed;
        \item[(ii)] the adjoint relation is given by
        \begin{equation*}
            (\ran A_{\ker \left[\begin{smallmatrix}
        K & L
    \end{smallmatrix} \right]})^*= \ran A_{\ran \left[\begin{smallmatrix}
                L^H \\ -K^H
            \end{smallmatrix}\right]};
        \end{equation*}
        \item[(iii)] $\ran A_{\ker \left[\begin{smallmatrix}
        K & L
    \end{smallmatrix} \right]}$ is self-adjoint if and only if $KL^H=LK^H$ and $\rank \begin{bmatrix}
        K \hspace{-0.5mm}&\hspace{-0.5mm} L
    \end{bmatrix}\hspace{-0.5mm}=\hspace{-0.5mm}\rank A$;
        \item[(iv)] $\ran A_{\ker \left[\begin{smallmatrix}
        K & L
    \end{smallmatrix} \right]}$ is maximal dissipative if and only if $KL^H+LK^H \geq 0$ and $\rank \begin{bmatrix}
        K & L
    \end{bmatrix}=\rank A$.
    \end{itemize}
\end{cor}
\begin{proof}
\begin{itemize}
    \item[(i)] We define $\Theta= \ker\begin{bmatrix}
        K & L
    \end{bmatrix}$. Since $\begin{bmatrix}
        K & L
    \end{bmatrix}$ is bounded we conclude that $\Theta$ is closed. Theorem \ref{thm:intermediate_extensions} (ii) immediately yields the result.
\item[(ii)] We use Theorem \ref{thm:intermediate_extensions} (iii) to conclude that $(\ran A_\Theta)^*=\ran A_{\Theta^*}$. It remains to show that $\Theta^*=\ran \left[\begin{smallmatrix}
                L^H \\ -K^H
            \end{smallmatrix} \right]$ which follows from \cite[p. 1017]{GernHR21}. 
\item[(iii)] By Theorem~\ref{thm:self-adj_extension} it suffices to show that $\Theta = \ker\begin{bmatrix}
        K & L
    \end{bmatrix}$ self-adjoint if and only if $KL^H=LK^H$ and $\rank \begin{bmatrix}
        K & L
    \end{bmatrix}=\rank A$. This statement is shown in \cite[Lem. 3.2]{GernHR21}.
\item[(iv)] Again, by Theorem~\ref{thm:self-adj_extension} we show that $\Theta = \ker\begin{bmatrix}
        K & L
    \end{bmatrix}$ is maximal dissipative if and only if $KL^H+LK^H\geq 0$ and $\rank \begin{bmatrix}
        K & L
    \end{bmatrix}=\rank A$ which was shown in \cite[Lem. 3.5]{GernHR21}.
    \end{itemize}
\end{proof}    
Note that the characterization of maximal dissipative restrictions of $H^*$ in Corollary~\ref{cor:kernel_domain} (iv)  paves the way to solution theory via generalized Cauchy problems, cf. \cite{arendt2023semigroups}, which is subject to future work.

\subsection{Characterization of skew-adjoint intermediate extensions of matrix differential operators in one spatial dimension}
In this part, we deduce a boundary triplet for (the adjoint of) a skew-symmetric matrix differential operator.
\begin{defn}\label{def:J}
Let coefficient matrices $J_k \in \C^{n \times n}$ for $k=0,\ldots,M$ satisfy
    \begin{equation}\label{eq:skew-symm}
        J_k=(-1)^{k+1}J_k^H.
    \end{equation}
Then we define the operator $\calJ_0 \colon \calX \supset C_c^\infty(a,b; \C^n) \to \calX$ via
     \begin{equation}\label{eq:J}
            \calJ_0 x = \sum\limits_{k=0}^M \dif{k} (-1)^{k}J_k^H x
    \end{equation}
\end{defn}
\noindent Observe that the operator $\calJ$ is skew-symmetric, i.e.,
\begin{equation*}
    \left\langle \calJ_0 x , y \right\rangle_\calX+\left\langle  x , \calJ_0 y \right\rangle_\calX =0
\end{equation*}
for all $x,y \in \dom \calJ_0$. We now define a boundary triplet for $\calJ_0^*$. To this end, by Lemma \ref{lem:A_star} we conclude that the adjoint is given by $\calJ_0^* \colon \calX \supset \dom \calJ_0^* \to \calX$ where
\begin{equation*}
        \calJ_0^* x = \sum\limits_{k=0}^M  J_k \tdif{k} x, \qquad
    \dom \calJ_0^* =H(\calJ_0^*;a,b) .
\end{equation*}
Now, inspired by \cite{le2005dirac} we define the matrix $Q=(Q_{ij})_{i,j=1, \ldots, N} \in \C^{Nn \times Nn}$ where
\begin{equation}\label{eq:Q}
    Q_{ij}=\begin{cases}
        0 , & \text{if} \ i+j >N \\
        (-1)^{j-1}J_k , & \text{if} \ i+j-1 =k. 
    \end{cases}
\end{equation}
It follows from that the symmetry assumption \eqref{eq:skew-symm} that $Q=Q^H$. The spectral theorem for Hermitian matrices now yields that there exist a unitary matrix $S \in \C^{Nn \times Nn}$ and a diagonal matrix $K \in \C^{Nn \times Nn}$ consisting of the eigenvalues of $Q$ such that $Q=SKS^H$. We let $S_{\mathrm{r}} \in \C^{Nn \times \rank Q}$ be the matrix whose columns are the eigenvectors corresponding to the nonzero eigenvalues of $Q$. Further, we let $K_{\mathrm{r}} \in \C^{\rank Q \times \rank Q}$ be the diagonal matrix containing only the nonzero eigenvalues of $Q$. Thus, we write the reduced factorization 
    \begin{equation*}
        Q=S_{\mathrm{r}}K_{\mathrm{r}}S_{\mathrm{r}}^H.
    \end{equation*}
    We now define the boundary maps $\Gamma =  \left[ \begin{smallmatrix}
        \Gamma_0 \\ \Gamma_1
    \end{smallmatrix}\right] : H^{N}(a,b;\C^n) \to \C^{\rank Q} \times \C^{\rank Q}$ via
    \begin{equation*}\begin{bmatrix}
        \Gamma_0 \\ \Gamma_1
    \end{bmatrix} x = \frac{1}{\sqrt{2}} \begin{bmatrix}
        K_{\mathrm{r}} & -K_{\mathrm{r}} \\ I & I
    \end{bmatrix} \begin{bmatrix}
        S_{\mathrm{r}}^H &0  \\ 0& S_{\mathrm{r}}^H
    \end{bmatrix} \gamma x,
    \end{equation*}
    where $\gamma: H^{N}(a,b;\C^n) \to \C^{2Nn}$ denotes the trace operator given by \eqref{eq:trace}.
\begin{thm}\label{thm:skew_sym_bdd_triplet}
    Let $\calJ_0: \calX \supset \dom \calJ_0 \to \calX$ be defined by \eqref{eq:J}. Then, $\left(\C^{\rank Q}, \widehat{\Gamma}_0, \widehat{\Gamma}_1 \right)$ is a boundary triplet for $\calJ_0^*$ where $\widehat{\Gamma}_0, \widehat{\Gamma}_1: H(\calJ_0^*;a,b) \to \C^{\rank Q} \times \C^{\rank Q}$ are the unique extensions of the boundary maps $\Gamma_0,\Gamma_1$.
\end{thm}
\begin{proof}
    We verify Definition \ref{def:operator_bdd_triplet} step by step. Clearly, $\C^{\rank Q}$ is a Hilbert space. Furthermore, by \cite[Thm. 3.1]{le2005dirac} we have 
    \begin{equation*}
            \langle \calJ_0^* x,y \rangle_\calX + \langle  x,\calJ_0^*y \rangle_\calX
            = \left\langle \gamma x,  \left[ \begin{smallmatrix}
        Q & 0\\  0& -Q
    \end{smallmatrix} \right]  \gamma y\right\rangle_{\C^{2Nn}}
    \end{equation*}
    for all $x,y \in \dom H^M(a,b;\C^n)$ where $Q$ is defined by \eqref{eq:Q}. Recall the reduced factorization $Q=S_{\mathrm{r}}K_{\mathrm{r}}S_{\mathrm{r}}^H$ and observe that 
\begin{equation*}
    	\left(\tfrac{1}{\sqrt{2}} \left[ \begin{smallmatrix}
        K_{\mathrm{r}} & -K_{\mathrm{r}} \\ I & I
    \end{smallmatrix} \right]  \left[ \begin{smallmatrix}
        S_{\mathrm{r}}^H & 0 \\ 0& S_{\mathrm{r}}^H
    \end{smallmatrix} \right] \right)^H\left[ \begin{smallmatrix}
    			    	  0& I \\ I &0 
    		    \end{smallmatrix} \right] \left[ \begin{smallmatrix}
        K_{\mathrm{r}} & -K_{\mathrm{r}} \\ I & I
    \end{smallmatrix} \right]  \left[ \begin{smallmatrix}
        S_{\mathrm{r}}^H & 0 \\ 0 & S_{\mathrm{r}}^H
    \end{smallmatrix} \right] \tfrac{1}{\sqrt{2}} =\left[ \begin{smallmatrix}
    			    	Q  & 0 \\ 0 & -Q
    		    \end{smallmatrix} \right] 
    \end{equation*}
and hence
\begin{equation*}
\begin{array}{rcl}
    \langle \calJ_0^* x,y \rangle_\calX + \langle  x,\calJ_0^*y \rangle_\calX
            &=& \left\langle \gamma x,  \left[ \begin{smallmatrix}
        Q & 0  \\ 0 & -Q
    \end{smallmatrix} \right]  \gamma y\right\rangle_{\C^{2Nn}} \\
        &=& 
    \left\langle \Gamma_1 x, \Gamma_0 y\right\rangle_{\C^{\rank Q}} + \left\langle \Gamma_0x, \Gamma_1 y\right\rangle_{\C^{\rank Q}}.
\end{array}
\end{equation*}
Obviously, $S_{\mathrm{r}}^H$ is surjective and since the trace operator $\gamma$ is surjective as well we conclude that $\Gamma$ is surjective. Using the same argument as in the proof of Lemma \ref{lem:extension_by_density} we find that $\Gamma$ uniquely extends to the surjective map $\widehat{\Gamma}$ which shows the claim.
\end{proof}
Note that for invertible leading coefficients $J_M$ we have that the matrix $Q$ defined in \eqref{eq:Q} is also invertible which leads to another choice of the boundary maps without the computation of a spectral factorization of $Q$, cf. \cite{le2005dirac, JacoZwar12}. Since we are able to parametrize the relations in the boundary space via (kernels of) bounded operators, we are yet able to provide a similar characterization of all skew-adjoint intermediate extensions as in \cite[Thm. 2.3]{kurula2015linear}.
\begin{cor}\label{cor:skew_kernel_domain}
    Let $\calJ_0: \calX \supset \dom \calJ_0 \to \calX$ be given by \eqref{eq:J} and let $\left( \C^{\rank Q}, \widehat{\Gamma}_0, \widehat{\Gamma}_1 \right)$ be the boundary triplet for $\calJ_0^*$ from Theorem \ref{thm:skew_sym_bdd_triplet}. Moreover, let $F,E \in \C^{\rank Q \times \rank Q}$. Then  the restriction $\euJ \coloneqq \calJ_0^* \left. \right|_{\dom \euJ}$ to
    \begin{equation}\label{eq:skew_adj_dom}
        \dom \euJ \coloneqq \ker \begin{bmatrix}
       F & E
    \end{bmatrix}\widehat{\Gamma}
    \end{equation}
    has the following properties:
    \begin{itemize}
        \item[(i)] $\euJ$ is closed;
        \item[(ii)] $-\euJ^*=\calJ_0^* \left. \right|_{\dom \euJ^*}$, where
        \begin{equation*}
            \dom \euJ^* = \left\lbrace x \in H(\calJ_0^*;a,b)  \,  \middle| \, \widehat{\Gamma} x \in \ran \begin{bmatrix}
                E^H \\ F^H
            \end{bmatrix}\right\rbrace ;
        \end{equation*}
        \item[(iii)] $\euJ$ is maximal dissipative if and only if 
        \begin{equation}\label{eq:diss_matrices}
            FE^H+EF^H \geq 0
        \end{equation}
        and $\rank \begin{bmatrix}
        F & E
    \end{bmatrix}=\rank Q$ and skew-adjoint if \eqref{eq:diss_matrices} holds with equality.
    \end{itemize}
\end{cor}
\begin{proof}
\begin{itemize}
\item[(i)] Define $\Theta = \ker\begin{bmatrix}
        F & E
    \end{bmatrix}$ and note that $\dom \euJ=\dom \euJ_\Theta$ where $\euJ_\Theta$ is given by \eqref{eq:skew_intermediate_dom}. Hence $\euJ=\euJ_\Theta$ is an intermediate extension of $-\calJ_0$ and thus $\euJ$ is closed if and only of $\Theta$ is closed by Proposition~\ref{prop:operator_intermediate_extensions} $(i)$. This is the case since $\begin{bmatrix}
        F & E
    \end{bmatrix}$ is bounded which yields that $\Theta$ is closed. 
    
    \item[(ii)] We conclude that $-\euJ^*=\euJ_{-\Theta^*}$ where we used Proposition~\ref{prop:operator_intermediate_extensions} $(ii)$. Again, by \cite[p. 1017]{GernHR21} we obtain $\Theta^*=\ran \left[\begin{smallmatrix}
                E^H \\ -F^H
            \end{smallmatrix} \right]$ which yields $-\Theta^*=\ran \left[\begin{smallmatrix}
                E^H \\ F^H
            \end{smallmatrix} \right]$ and hence the claim.
\item[(iii)] By Proposition~\ref{prop:operator_intermediate_extensions} $(iii)$ we have that $\euJ$ is maximal dissipative (skew-adjoint) if and only if $\Theta$ is maximal dissipative (skew-adjoint). The result follows immediately from \cite[Lem. 3.5]{GernHR21}.
\end{itemize}
\end{proof}

\section{Application to generalized implicit port-Hamiltonian systems}\label{sec:applications}
For linear PDEs a commonly used physics-based representation is given by a first-order port-Hamiltonian (pH) system
\begin{align}\label{eq:birgit_hans}
 \frac{\mathrm{d}}{\mathrm{d}t} x(t,\xi)= P_1\tdif{}\mathcal{H}(\xi)x(t,\xi)+P_0\mathcal{H}(\xi)x(t,\xi),\quad  x(0,\xi)=x_0(\xi),
\end{align}
which is considered as an abstract Cauchy problem with states $x(t,\cdot)$ in the Hilbert space of square-integrable functions $\calX=L^2(a,b;\C^n)$ on a bounded one-dimensional spatial interval $(a,b)\subset\mathbb{R}$, where $P_1=P_1^*$ and $P_0=-P_0^*$ and the Hamiltonian density $\mathcal{H}\in L^\infty(a,b;\C^{n\times n})$ is used to describe the total energy of a given state $x(t,\cdot)$. We recall the recent generalization of \eqref{eq:birgit_hans}, the port-Hamiltonian system dynamics \eqref{eq:pHs}
\begin{align}
\label{eq:PJS}
    \frac{\mathrm{d}}{\mathrm{d}t} \underbrace{\sum\limits_{k=0}^N \tdif{k} P_k \tdif{k}  }_{=\calP} x(t, \xi)=\displaystyle \underbrace{\sum\limits_{k=0}^M J_k \tdif{k} }_{=\calJ} \big( \underbrace{\sum\limits_{l=0}^N \tdif{l} S_l \tdif{l} }_{=\calS} x(t,\xi) \big) 
\end{align}
where the coefficient matrices $P_k,S_l$ and $J_k$ satisfy Assumption~\ref{ass:ass1}, \ref{ass:mdo} and \ref{eq:skew-symm}, respectively.  Formally, we can rewrite \eqref{eq:PJS} in the following way
\begin{align*}
    (x(t), \tfrac{\mathrm{d}}{\mathrm{d}t} x(t)) \in \underbrace{\gr \calJ \vphantom{\left[\begin{smallmatrix}
        \calP \\ \calS
    \end{smallmatrix}\right]}}_{\eqqcolon \calD} \underbrace{\ran \left[\begin{smallmatrix}
        \calP \\ \calS
    \end{smallmatrix}\right]}_{\eqqcolon \calL}.
\end{align*}
If the domains of $\calJ$ and $\calL$ are chosen such that $\calD=-\calD^*, \calL=\calL^*$, then we obtain a generalized pH system as e.g.\ studied in finite-dimensional spaces in \cite{MascvdSc18}. Another possibility is to choose the domains such that $\calD$ is maximally dissipative and $\calL$ is maximally nonnegative, i.e.\ $\calL=\calL^*$ and $-\calL$ is dissipative which leads to a dissipative pH framework as studied in finite-dimensional spaces in \cite{GernHR21, GerPPS23, MascvdSc23, MehS23}. 

In the result below, we summarize for which boundary conditions, the implicit pH system \eqref{eq:PJS} can be regarded as a generalized pH system in the above sense. The proposition follows immediately from Corollary \ref{cor:kernel_domain} and \ref{cor:skew_kernel_domain}.
\begin{prop}\label{cor:generalized}
Let $A\in\mathbb{C}^{2Nn}$ and $\left(\C^{\rank A}, \widetilde{\Gamma}_0, \widetilde{\Gamma}_1\right)$ be the boundary triplet for $H^*= \ran \left[\begin{smallmatrix}
        \calP \\ \calS
    \end{smallmatrix}\right]$ given by Theorem~\ref{thm:boundary_triplet} and $\left(\C^{\rank Q}, \widehat{\Gamma}_0, \widehat{\Gamma}_1 \right)$ is a boundary triplet for $\calJ_0^*$ constructed in Theorem \ref{thm:skew_sym_bdd_triplet}.  Moreover, let $K, L \in \C^{\rank A \times \rank A}$ and $F,E \in \C^{\rank Q \times \rank Q}$. Then, 
    \begin{align*}
        \left( \gr \left. -\calJ_0^* \right|_{\ker \begin{bmatrix}
       F & E
    \end{bmatrix}} \right) \left( \ran \left. \left[\begin{smallmatrix}
        \calP \\ \calS
    \end{smallmatrix}\right]\right|_{\ker \begin{bmatrix}
            K & L
    \end{bmatrix}}\right)
    \end{align*}
    is a generalized port-Hamiltonian system if and only if $FE^H=-EF^H$ with $\rank \begin{bmatrix}
        F & E
    \end{bmatrix}=\rank Q$ and $KL^H=LK^H$ and $\rank  \begin{bmatrix}
            K & L
    \end{bmatrix}=\rank A$.
\end{prop}
Eventually we note that, in contrast to the finite-dimensional case, Proposition~\ref{cor:generalized} does not immediately lead to strongly continuous or contraction semigroups to describe the solution to \eqref{eq:PJS}, because we consider products of relations or operators in $\calX$.

\subsection{Dzektser equation revisited}
    Recall the Dzektser equation
    \begin{equation}
        \frac{\mathrm{d}}{\mathrm{d} t} \left(1+ \dif{}1\dif{}\right)x(t,\xi)=\left(\dif{}1\dif{} + \dif{} 2 \dif{}\right)x(t,\xi),
    \end{equation}
    $t > 0, \xi \in (0, \pi)$ where our notion emphasizes the specific form \eqref{eq:matrix_diff_op}. In Section \ref{subsec:dzektser} we defined the associated operators $\calP, \calS$ and we derived a boundary triplet for $\ran \left[\begin{smallmatrix}
        \calP \\ \calS
    \end{smallmatrix}\right]$. However, we will see that the matrix $A \in \C^{4 \times 4}$ given in \eqref{eq:bt_dexter} can also be defined in terms of the coefficient matrices 
        \begin{equation*}
            \begin{array}{rclcrclcrcl}
                P_0&=&1, & \quad P_1&=&1, & \quad P_2&=&0, \\
                S_0&=&0, & \quad S_1&=&1, & \quad S_2&=&2,
            \end{array}
        \end{equation*}
        that obviously satisfy the symmetry conditions \eqref{eq:symmetric_matrices}.
    Moreover, there are three different pairs $(P_k,S_l)\neq (0,0)$ with $k \neq l$, namely $(P_0,S_1), (P_0,S_2)$ and $(P_1, S_2)$. Thus, Assumption \ref{ass:mdo} is satisfied. To derive a boundary triplet for $\ran \left[\begin{smallmatrix}
        \calP \\ \calS
    \end{smallmatrix}\right]$ we apply Theorem~\ref{thm:boundary_triplet} using the matrix $A=B-B^H$ from \eqref{eq:A_matrix} which can be obtained from
    \begin{equation*}
        B=\left( \sum\limits_{l=1}^N \sum\limits_{k=0}^{l-1} \sum\limits_{\substack{m=2l-1 \\ m \text{ odd}}}^{2N-1} M_{lkm} \right)^H = \left( \sum\limits_{l=1}^2 \sum\limits_{k=0}^{l-1} \sum\limits_{\substack{m=2l-1 \\ m \text{ odd}}}^{3} M_{lkm} \right)^H.
    \end{equation*}
    Thus, we find the matrices 
    \begin{equation*}
        \begin{array}{rclcrcl}
            M_{101} &=& \left[\begin{smallmatrix}
                0&0&0&0 \\ 1 &0&0&0 \\ 0&0&0&0 \\ 0&0&0&0
            \end{smallmatrix} \right], & \qquad M_{103} &=& \left[\begin{smallmatrix}
                0&0&0&0 \\ 0&0&0&0 \\ 0&0&0&0 \\ 0&0&2&0
            \end{smallmatrix} \right] \\[2ex]
            M_{203} &=& \left[\begin{smallmatrix}
                0&0&0&0 \\ 0&0&0&0 \\ 0&0&0&0 \\ 2&0&0&0
            \end{smallmatrix} \right], &\qquad M_{213} &=& \left[\begin{smallmatrix}
                0&0&0&0 \\ 0&0&0&0 \\ 0&-2&0&0 \\ 0&0&0&0
            \end{smallmatrix} \right]
        \end{array}
    \end{equation*}
    which yields
    \begin{equation*}
        B=(M_{101}+ M_{103} + M_{203} + M_{213})^H = \left[\begin{smallmatrix}
                0&0&0&0 \\ 1&0&0&0 \\ 0&-2&0&0 \\ 2&0&2&0
            \end{smallmatrix} \right]^H
    \end{equation*}
    and hence
    \begin{equation*}A^H=B^H-B= \left[\begin{smallmatrix}
                0&-1&0&-2 \\ 1&0&2&0 \\ 0&-2&0&-2 \\ 2&0&2&0
            \end{smallmatrix} \right].
    \end{equation*}
    Observe that this is exactly the matrix $A$ as defined in Section~\ref{subsec:dzektser}, where a boundary triplet $(\C^4, \Gamma_0, \Gamma_1)$ for the Dzektser equation was already defined in Theorem~\ref{thm:Dzektser_bdd_triplet}.
\par
Building on this, we want to express the boundary conditions
    \begin{equation}\label{eq:dirichlet_bdd}
        x(0,t)=x(\pi,t)=0, \qquad x''(0,t)=x''(\pi,t)=0, \qquad t > 0
    \end{equation}
    proposed in \cite{jacob2022solvability} in terms of the boundary maps. To this end, observe that
    \begin{equation}\label{eq:Gamma_to_trace}
        \gamma x =\tfrac{\sqrt{2}}{2} \begin{bmatrix}
        A^{-1} & -I \\ -A^{-1} & -I
    \end{bmatrix} \tfrac{1}{\sqrt{2}} \begin{bmatrix}
        A & -A \\ -I & -I
    \end{bmatrix} \gamma x = \tfrac{\sqrt{2}}{2} \begin{bmatrix}
        A^{-1} & -I \\ -A^{-1} & -I
    \end{bmatrix} \Gamma x.
    \end{equation}
    Moreover, we observe that \eqref{eq:dirichlet_bdd} is equivalent to
    \begin{equation*}
         0=\underbrace{\begin{bmatrix}
             1 &0&0&0 \\
             0& 0 & 1 &0 \\
             0&0&0&0 \\
             0&0&0&0
             \end{bmatrix}}_{\eqqcolon Q} \begin{bmatrix}
             x(\pi) \\ x'(\pi) \\ x''(\pi)\\  x^{(3)}(\pi) \end{bmatrix}  +   \underbrace{ \begin{bmatrix}
             0&0&0&0 \\
             0&0&0&0 \\
             1 &0&0&0 \\
             0&0& 1 & 0 
             \end{bmatrix}}_{\eqqcolon R}
         \begin{bmatrix} x(0) \\ x'(0) \\ x''(0)\\  x^{(3)}(0)
         \end{bmatrix} = \begin{bmatrix}
            Q & R
        \end{bmatrix} \gamma x ,
    \end{equation*}
    which can be also written as 
    \begin{equation*}
        \Gamma x \in  \ker \begin{bmatrix}
            Q & R
        \end{bmatrix}  \begin{bmatrix}
        A^{-1} & -I \\ -A^{-1} & -I
    \end{bmatrix} = \ker \begin{bmatrix}
            (Q - R)A^{-1} & -(Q+R)
        \end{bmatrix}
    \end{equation*}
    where we used \eqref{eq:Gamma_to_trace}. It is easily verified that $K\coloneqq (Q - R)A^{-1} \in \C^{4 \times 4}$ and $L=-(Q+R) \in \C^{4 \times 4}$ satisfy
    $KL^H=LK^H$ and $\rank \begin{bmatrix}
            K & L
        \end{bmatrix}=4$. Consequently, by Corollary \ref{cor:kernel_domain} $(iii)$ we have that the range of
        $\left[\begin{smallmatrix}
        \calP \\ \calS
    \end{smallmatrix}\right]: \calX \supset \Gamma^{-1}\ker  \begin{bmatrix}
            K & L
         \end{bmatrix} \to \calX \times \calX$ 
        defines a self-adjoint linear relation in $\calX=L^2(0,\pi)$. Eventually, the boundary conditions \eqref{eq:dirichlet_bdd} lead to a self-adjoint representation of the system dynamics of the Dzektser equation.

\subsection{Biharmonic wave equation as a port-Hamiltonian system}
To illustrate the system class definition introduced herein, we reinterpret the biharmonic wave equation \cite{roetman1967biharmonic} as a port-Hamiltonian system, applying Theorems \ref{thm:boundary_triplet} and \ref{thm:skew_sym_bdd_triplet} to the resulting matrix differential operators. Consider the biharmonic wave equation
\begin{equation}\label{eq:biharmonic_wave}
    \begin{array}{rcll}
        \frac{\mathrm{d}^2}{\mathrm{d} t^2} w(t,\xi)&=&- \dif{4} w(t,\xi), \qquad & \xi \in (a,b), t \geq 0,  \\
        w(0, \xi)&=&w_0(\xi) & \xi \in (a,b), \\
        \frac{\mathrm{d}}{\mathrm{d} t} w(0, \xi)&=&w_1(\xi) & \xi \in (a,b),
    \end{array}
\end{equation}
where $w(t,\xi)$ denotes the elongation of the wave and $w_0, w_1 \in L^2(a,b)$ are initial elongation and momentum.
To derive a system representation as a generalized port-Hamiltonian system we introduce the new state variable
\begin{equation*}
    x(t,\xi) \coloneqq \begin{bmatrix}
        x_1(t,\xi) \\ x_2(t,\xi)
    \end{bmatrix} \coloneqq  \begin{bmatrix}
        \frac{\mathrm{d}}{\mathrm{d} t} w(t,\xi)  \\ \vphantom{\dif{4}} \dif{} w(t,\xi)
    \end{bmatrix}
\end{equation*}
consisting of the momentum $x_1$ and the strain $x_2$. Suppose that $w$ solves~\eqref{eq:biharmonic_wave}. Then
\begin{equation*}
     \frac{\mathrm{d}}{\mathrm{d} t}x(t,\xi) = \begin{bmatrix}
         \frac{\mathrm{d}^2}{\mathrm{d} t^2} w(t,\xi) \\ \dif{} \vphantom{\dif{4}} \frac{\mathrm{d}}{\mathrm{d} t}  w(t,\xi) 
    \end{bmatrix}= \underbrace{\begin{bmatrix}
             0 & \dif{} \\ \dif{} & 0 \vphantom{-\dif{2}}
         \end{bmatrix}}_{\eqqcolon \calJ}\underbrace{ \begin{bmatrix}
             1 & \vphantom{\dif{}}0 \\ 0 & -\dif{2}
         \end{bmatrix}}_{\eqqcolon \calS} \begin{bmatrix}
        x_1(t,\xi) \\ x_2(t,\xi)
    \end{bmatrix}.
\end{equation*}
The involved differential operators are parametrized by matrices via
\begin{equation}\label{eq:waveJQ}
    \calP= I, \qquad \calJ = \begin{bmatrix}
             0 & 1 \\ 1 & 0
         \end{bmatrix} \dif{}, \qquad \calS 
         = \dif{} \begin{bmatrix}
             0 & 0 \\ 0 & -1
         \end{bmatrix}\dif{} + \begin{bmatrix}
             1 & 0 \\ 0 & 0
         \end{bmatrix}
\end{equation}
as well as $\calX \coloneqq L^2(a,b; \C^2)$ and $\calS_0: \calX \supset \dom \calS_0 \to \calX$ via
\begin{align*}
    \calS_0 x = \begin{bmatrix}
             1 & 0 \\ 0 & -\dif{2}
         \end{bmatrix} x , \qquad \dom \calS_0 = C_c^\infty(a,b;\C^2).
\end{align*}
Moreover, define $\calJ_0: \calX \supset \dom \calJ_0 \to \calX$ by
\begin{align*}
    \calJ_0 x = \begin{bmatrix}
             0 & \dif{} \\ \dif{} & 0
         \end{bmatrix} x , \qquad \dom \calJ_0 = C_c^\infty(a,b;\C^2).
\end{align*}
It is easy to see that $\calS_0^*=\calS$ where $\dom \calS = L^2(a,b) \times H^2(a,b)$ and $-\calJ_0^*=\calJ$ where $\dom \calJ = H^1(a,b;\C^2)$. Since $\calP=I$ is trivially invertible, Assumption~\ref{ass:ass1} is fulfilled and we consider $H=\gr \calS_0$ and we aim to use Theorem~\ref{thm:boundary_triplet} to derive a boundary triplet for $H^*=\gr \calS=\ran \left[\begin{smallmatrix}
        I \\ \calS
    \end{smallmatrix}\right]$. To this end, observe that the matrix defined by \eqref{eq:A_matrix} is
\begin{equation*}
    A
    = \begin{bmatrix}
         0&0& 0 &  0 \\
         0&0& 0 & 1 \\
         0 & 0 &0&0 \\
         0 & -1 &0&0
    \end{bmatrix} \in \C^{4 \times 4}.
\end{equation*}
It is easy to verify that the matrices
    \begin{equation*}
        U_{\mathrm{r}} \coloneqq \frac{1}{\sqrt{2}}\begin{bmatrix}
           0 & 0 \\ -i & i \\ 0 & 0 \\ 1 & 1
        \end{bmatrix} \in  \C^{4 \times 2}, \qquad D_{\mathrm{r}} \coloneqq \begin{bmatrix}
          i & 0 \\ 0 & -i
        \end{bmatrix} \in  \C^{2 \times 2}
    \end{equation*}
    satisfy $A = U_{\mathrm{r}}D_{\mathrm{r}}U_{\mathrm{r}}^H$. Define the boundary map $\Gamma =  \left[ \begin{smallmatrix}
        \Gamma_0 \\ \Gamma_1
    \end{smallmatrix}\right] : L^2(a,b) \times  H^{2}(a,b) \to \C^{2} \times \C^{2}$ via
    \begin{equation*}
         \left[ \begin{smallmatrix}
        \Gamma_0 \\ \Gamma_1
    \end{smallmatrix}\right] x = \tfrac{1}{\sqrt{2}} \left[ \begin{smallmatrix}
        D_{\mathrm{r}} & -D_{\mathrm{r}} \\ -I & -I
    \end{smallmatrix} \right] \left[ \begin{smallmatrix}
        U_{\mathrm{r}}^H &0  \\ 0& U_{\mathrm{r}}^H
    \end{smallmatrix} \right] \gamma x =  \tfrac{1}{2} \left[ \begin{smallmatrix}
        i & 0 & - i & 0 \\ 0 & -i & 0 & i \\ -1 & 0 & -1 & 0 \\ 0 & -1 & 0 & -1
        \end{smallmatrix} \right] \left[ \begin{smallmatrix}
         i & 1 & 0 & 0 \\ -i & 1 & 0 & 0 \\ 0 & 0 & i & 1 \\ 0 & 0 & -i & 1
    \end{smallmatrix} \right] \left[
        \begin{smallmatrix}
            x_2'(b) \\  x_2(b) \\  x_2'(a) \\ x_2(a)
        \end{smallmatrix}\right].
    \end{equation*}
    Then, Theorem \ref{thm:boundary_triplet} yields that $(\C^2, \Gamma_0, \Gamma_1)$ is a boundary triplet for $H^*=\gr \calS$. Moreover, we conclude that $(\C^2, \widehat{\Gamma}_0, \widehat{\Gamma}_1)$ where
    \begin{equation*}
    \widehat{\Gamma} =  \left[ \begin{smallmatrix}
        \widehat{\Gamma}_0 \\ \widehat{\Gamma}_1
    \end{smallmatrix}\right] : H^{1}(a,b;\C^2) \to \C^{2} \times \C^{2},\quad 
         \left[ \begin{smallmatrix}
         \widehat{\Gamma}_0 \\ \widehat{\Gamma}_1
    \end{smallmatrix}\right] x = \tfrac{1}{\sqrt{2}} \left[ \begin{smallmatrix}
        0 & 1 & 0 & -1 \\ 1 & 0 & -1 & 0 \\ 1 & 0 & 1 & 0 \\ 0 & 1 & 0 & 1
        \end{smallmatrix}\right]
         \left[ \begin{smallmatrix}
            x_1(b) \\  x_2(b) \\  x_1(a) \\ x_2(a)
        \end{smallmatrix}\right].
    \end{equation*}
    is a boundary triplet for $\calJ_0^*$ where we used \cite[p. 85]{JacoZwar12}. Having defined the boundary triplets for $H^*$ and $\calJ_0^*$, respectively, one may characterize all generalized port-Hamiltonian systems associated with the biharmonic wave equation in the spirit of Proposition~\ref{cor:generalized}.
    
\subsection{Elastic rod with non-local elasticity condition}
We consider the example of an elastic rod with non-local elasticity relation which was provided in \cite{heidari2019port, MascvdSc23}. In the following, $T>0$ denotes the Young's modulus of the rod, $\rho A>0$ is the mass density and $k>0$ is a positive constant. Now let the state be build by  the displacement $u=u(t,\xi)$ of the rod, $\lambda=\lambda(t,\xi)$ the so-called latent strain and $p=p(t,\xi)$ the momentum density. Then, the time evolution is given by
	\begin{align*}\label{eq:elastic_rod}
		\frac{\partial}{\partial t} \underbrace{\begin{bmatrix}
				1 & 0 & 0 \\ 0 & 1 - \tdif{}\mu\tdif{} & 0 \\  0 & 0 & 1 \vphantom{\tdif{}}
		\end{bmatrix} }_{\eqqcolon \calP}
		\begin{bmatrix}
			u \\ \lambda \\ p 
		\end{bmatrix}
		=\underbrace{\begin{bmatrix}
				0 & 0 & 1 \\ 0 & 0 & \tdif{}  \\ -1 & \tdif{} & 0 
		\end{bmatrix}}_{\eqqcolon \calJ}\underbrace{\begin{bmatrix}
				k & 0 & 0 \\ 0 & T & 0 \\ 0 & 0 & \frac{1}{\rho A} 
		\end{bmatrix}}_{\eqqcolon \calS}
		\begin{bmatrix}
			u \\ \lambda \\ p 
		\end{bmatrix},
	\end{align*}
    where $\mu > 0$ is a positive parameter. We define $\calX=L^2(a,b;\C^3)$ and the bounded, coercive and self-adjoint multiplication operator $\calS: \calX \to \calX$. Further, define $\calP_0: \calX \supset C_c^\infty(a,b;\C^3) \to \calX$ via $\calP_0 x = \calP x$ for all $x \in \dom \calP_0$. Clearly, $\calP_0^*=\calP$ with 
    \begin{align*}
        \dom \calP = L^2(a,b) \times H^2(a,b) \times L^2(a,b).
    \end{align*}
    Moreover, it is easy to see that $\calJ_0: \calX \supset C_c^\infty(a,b;\C^3) \to \calX$ with $\calJ_0 x =- \calJ x$ for all $x \in \dom \calJ_0$ satisfies $\calJ_0^*=\calJ$ where $\dom \calJ=L^2(a,b) \times H^1(a,b;\C^2)$. Now, let $H=(\gr \calP_0)^{-1}(\gr \calS)$ and we derive a boundary triplet for $H^*=(\gr \calS)(\gr \calP)^{-1}=\ran \left[\begin{smallmatrix}
        \calP \\ \calS
    \end{smallmatrix}\right]$. Obviously, Assumption \ref{ass:ass1} holds for $H^*$ as $\calS$ is coercive which yields that the methods suggested in Section \ref{subsec:bdd_triplets} are applicable. We compute
    \begin{align*}
        \langle \calS x, \calP y \rangle_\calX 
        &= \langle \calS x,  y \rangle_\calX - \langle T x_2, \tdif{}\mu \tdif{} y_2 \rangle_{L^2(a,b)} \\ 
        &= \langle  x,  \calS y \rangle_\calX + \langle \mu T \tdif{}x_2, \tdif{}y_2 \rangle_{L^2(a,b)}- \left[ \langle  Tx_2(\xi),  \mu\tdif{} y_2(\xi) \rangle_\C \right]_{\xi=a}^{\xi=b} \\
        &= \langle  x,  \calS y \rangle_\calX -\langle \tdif{} \mu \tdif{}x_2, Ty_2 \rangle_{L^2(a,b)} \\
        &\hphantom{=} + \left[ \langle \mu \tdif{} x_2(\xi), T y_2(\xi) \rangle_\C \right]_{\xi=a}^{\xi=b}  - \left[ \langle  Tx_2(\xi),  \mu\tdif{} y_2(\xi) \rangle_\C \right]_{\xi=a}^{\xi=b} 
    \end{align*}
    which yields with $A \coloneqq \left[\begin{smallmatrix}
           0 & -\mu T \\ \mu T & 0
        \end{smallmatrix}\right]$
    \begin{align*}
        \langle \calS x, \calP y \rangle_\calX - \langle \calP x, \calS y \rangle_\calX = \left[ \left\langle A \left( \begin{smallmatrix}
            x_2 (\xi) \\ x_2'(\xi)
        \end{smallmatrix}\right), \left( \begin{smallmatrix}
            y_2 (\xi) \\ y_2'(\xi)
        \end{smallmatrix}\right) \right\rangle_{\C^2}\right]_{\xi=a}^{\xi=b}.
    \end{align*}
    Using Theorem~\ref{thm:boundary_triplet}~(iii), we conclude that $(\C^2, \Gamma_0, \Gamma_1)$ is a boundary triplet for $H^*$ where $\Gamma =  \begin{bmatrix}
        \Gamma_0 \\ \Gamma_1
    \end{bmatrix} : \dom \calP \to \C^2 \times \C^{2}$ is defined by
    \begin{equation*}
         \Gamma x = \begin{bmatrix}
        \Gamma_0 \\ \Gamma_1
    \end{bmatrix} x = \frac{1}{\sqrt{2}} \begin{bmatrix}
        A & -A \\ -I & -I
    \end{bmatrix}
        \begin{bmatrix}
            x_2(b) \\  x_2'(b) \\  x_2(a) \\ x_2'(a)
        \end{bmatrix}=  \frac{1}{\sqrt{2}} \begin{bmatrix}
            \mu T (x_2'(a)- x_2'(b))  \\ \mu T (x_2(b)- x_2(a))  \\  -(x_2(b) + x_2(a)) \\ -(x_2'(b)+ x_2'(a))
        \end{bmatrix}  .
    \end{equation*}
    To derive a boundary triplet for $\calJ_0^*$ we observe that the matrix $Q$ in \eqref{eq:Q} is just $Q=J_1$ and thus
    \begin{align*}
        Q=J_1=\begin{bmatrix}
            0 & 0 & 0 \\ 0 & 0 & 1 \\ 0 & 1 & 0 
        \end{bmatrix} = \frac{1}{\sqrt{2}}\begin{bmatrix}
             0 & 0 \\ 1 & -1 \\ 1 & 1 
        \end{bmatrix}\begin{bmatrix}
             1 & 0 \\ 0 & -1 
        \end{bmatrix} \begin{bmatrix}
             0 & -1 & 1 \\ 0 & 1 & 1 
        \end{bmatrix} \frac{1}{\sqrt{2}} \eqqcolon S_{\mathrm{r}}K_{\mathrm{r}}S_{\mathrm{r}}^H.
    \end{align*}
    Then, Theorem \ref{thm:skew_sym_bdd_triplet} yields that $(\C^2, \widehat{\Gamma}_0, \widehat{\Gamma}_1)$ is a boundary triplet for $\calJ_0^*$ where $\widehat{\Gamma} =  \left[ \begin{smallmatrix}
        \widehat{\Gamma}_0\\ \widehat{\Gamma}_1
    \end{smallmatrix}\right] : \dom \calJ \to \C^{2} \times \C^{2}$ is given by
    \begin{equation*}\begin{bmatrix}
       \widehat{\Gamma}_0\\ \widehat{\Gamma}_1
    \end{bmatrix} x = \frac{1}{\sqrt{2}} \begin{bmatrix}
        1 & 0 & - 1 & 0 \\ 0 & -1 & 0 & 1 \\ 1 & 0 & 1 & 0 \\ 0 & 1 & 0 & 1
        \end{bmatrix}\begin{bmatrix}
         -1 & 1 & 0 & 0 \\ 1 & 1 & 0 & 0 \\ 0 & 0 & -1 & 1 \\ 0 & 0 & 1 & 1
    \end{bmatrix}
        \begin{bmatrix}
            x_2(b) \\  x_3(b) \\  x_2(a) \\ x_3(a)
        \end{bmatrix}.
    \end{equation*}

\section{Conclusion}
We provided a novel notion of boundary triplets for range representations of linear relations which enabled us to characterize the set self-adjoint or maximally dissipative intermediate extensions of such relations via relations in the boundary space. Correspondingly, we gave a one-to-one correspondence between our approach and Lagrangian subspaces having their foundations in the geometrical modeling of port-Hamiltonian systems. As a second contribution of this paper, we characterized suitable boundary maps for a class of implicit port-Hamiltonian systems over one-dimensional spatial domains. We established abstract Green's identities through the coefficient matrices of matrix differential operators. We provided various applications of our results by means of the Dzektser equation, the biharmonic wave equation and an elastic rod. 
\section*{Acknowledgment}
The authors thank Nathanael Skrepek (University of Twente) for pointing out a missing assumption in an earlier version of this work. 


\section*{References}
\bibliographystyle{abbrv}
\bibliography{sample.bib}
\end{document}